\documentclass{amsart}
\usepackage{amssymb,amsmath,amscd,amsthm,enumerate,color,verbatim}
\usepackage[all,cmtip]{xy}
\usepackage{latexsym}
\usepackage{enumitem}
\usepackage{graphicx}
\usepackage{subfigure}
\usepackage{tikz}

\def\SL{{\mathrm {SL}}}
\def\loc{{\mathrm{loc}}}

\def\beq{\begin{equation}}
\def\eeq{\end{equation}}

\newcommand{\widecheck}{\check}

\newcommand{\Z}{{\mathbb Z}}
\newcommand{\R}{{\mathbb R}}

\newcommand{\C}{{\mathbb C}}
\newcommand{\D}{{\mathbb D}}

\newcommand{\N}{{\mathbb N}}

\renewcommand{\H}{{\mathbb H}}

\newcommand{\bbP}{{\mathbb{P}}}

\newcommand{\tr}{{\mathrm{Tr}}}
\newcommand{\inter}{\operatorname{int}}

\newcommand{\CA}{{\mathcal A}}

\newcommand{\DD}{{\mathcal D}}

\newcommand{\CF}{{\mathcal F}}
\newcommand{\CG}{{\mathcal G}}
\newcommand{\CB}{{\mathcal B}}
\newcommand{\CH}{{\mathcal H}}

\newcommand{\CL}{{\mathcal L}}
\newcommand{\CM}{{\mathcal M}}
\newcommand{\CN}{{\mathcal N}}

\newcommand{\CP}{{\mathcal P}}

\newcommand{\CZ}{{\mathcal Z}}
\newcommand{\CU}{{\mathcal U}}

\newcommand{\CO}{{\mathcal O}}

\renewcommand{\vec}[1]{\underline{#1}}

\newtheorem{lemma}{Lemma}[section]
\newtheorem{theorem}[lemma]{Theorem}
\newtheorem{prop}[lemma]{Proposition}
\newtheorem{coro}[lemma]{Corollary}

\newtheorem*{problem}{Problem}

\theoremstyle{definition}
\newtheorem{definition}[lemma]{Definition}
\newtheorem{remark}[lemma]{Remark}

\begin{document}

\title[Positive Lyapunov Exponent]{Schr\"odinger Operators With Potentials Generated by Hyperbolic Transformations:\\I.~Positivity of the Lyapunov Exponent}


\author{Artur Avila, David Damanik, and Zhenghe Zhang}

\address{Institut f\"ur Mathematik, Universit\"at Z\"urich, Winterthurerstrasse 190,
CH-8057 Z\"urich, Switzerland and IMPA, Estrada Dona Castorina, 110, Rio de Janeiro, 22460-320, Brazil}

\email{artur.avila@math.uzh.ch}

\address{Department of Mathematics, Rice University, Houston, TX~77005, USA}

\email{damanik@rice.edu}

\address{Department of Mathematics, University of California, Riverside, CA-92521, USA}
\email{zhenghe.zhang@ucr.edu}

\thanks{D.D.\ was supported in part by NSF grants DMS--1067988, DMS--1301582, DMS--1700131 and by an Alexander von Humboldt Foundation research award}

\thanks{Z.\ Z.\ was supported in part by NSF grant DMS--1764154}

\begin{abstract}
We consider discrete one-dimensional Schr\"odinger operators whose potentials are generated by sampling along the orbits of a general hyperbolic transformation. Specifically, we show that if the sampling function is a non-constant H\"older continuous function defined on a subshift of finite type with an ergodic measure admitting a local product structure and a fixed point, then the Lyapunov exponent is positive away from a discrete set of energies. Moreover, for sampling functions in a residual subset of the space of H\"older continuous functions, the Lyapunov exponent is positive everywhere. If we consider locally constant or globally fiber bunched sampling functions, then the Lyapuonv exponent is positive away from a finite set. Moreover, for sampling functions in an open and dense subset of the space in question, the Lyapunov exponent is uniformly positive. Our results can be applied to any subshift of finite type with ergodic measures that are equilibrium states of H\"older continuous potentials. In particular, we apply our results to Schr\"odinger operators defined over expanding maps on the unit circle, hyperbolic automorphisms of a finite-dimensional torus, and Markov chains.
\end{abstract}

\maketitle

\tableofcontents

\section{Introduction}

\subsection{Statement of Results}\label{ss:results}

In this series of papers, we are mainly concerned with the Anderson localization phenomenon for one-dimensional discrete Schr\"odinger operators $H_\omega$ in $\ell^2(\Z)$ acting by
$$
[H_\omega \psi](n) = \psi(n+1) + \psi(n-1) + V_\omega(n) \psi(n).
$$
Here we assume $\Omega$ to be any compact metric space, $T:\Omega\to\Omega$ a homeomorphism, and $f : \Omega \to \R$ be continuous. We consider potentials $V_\omega : \Z \to \R$ defined by $V_\omega(n) = f(T^n \omega)$ for $\omega \in \Omega$ and $n \in \Z$. For general background on Schr\"odinger operators in $\ell^2(\Z)$ with dynamically generated potentials of this form, we refer the reader to \cite{D17, DF1, DF2}.

Spectral properties of the operators $H_\omega$ can be investigated by studying the behavior of the solutions to the difference equation
\begin{equation}\label{e.deve}
u(n+1) + u(n-1) + V_\omega(n) u(n) = E u(n), \quad n \in \Z
\end{equation}
with $E$ real or complex (depending on the problem in question). These solutions in turn can be described with the help of the Schr\"odinger cocycle $(T,A^E)$ with the cocycle map $A^{E}:\Omega\to\SL(2,\R)$ being defined as
$$
A^E(\omega) = A^{(E-f)}(\omega):= \begin{pmatrix} E-f(\omega) & -1 \\ 1 & 0 \end{pmatrix},
$$
where we often leave the dependence on $f:\Omega \to \R$ implicit as it will be fixed most of the time.

Such cocycles describe the transfer matrices associated with Schr\"odinger operators. Specifically, $u=u(n)$ solves \eqref{e.deve} if and only if
$$
\begin{pmatrix} u(n) \\ u(n-1) \end{pmatrix} = A^{E}_n(\omega) \begin{pmatrix} u(0) \\ u(-1) \end{pmatrix}, \quad n \in \Z,
$$
where
$$
A_n(\omega)=
\begin{cases}
	A(T^{n-1}\omega) \cdots A(\omega),\ & n\ge 1;\\
	[A_{-n}(T^n\omega)]^{-1}, \ &n\le -1,
\end{cases}
$$
and we set $A_0(\omega)$ to be the identity matrix.

The Lyapunov exponent (LE) of the Schr\"odinger cocycle plays a key role in the spectral analysis of the operators. Let $\mu$ be a $T$-ergodic probability measure on $\Omega$. The Lyapunov exponent is given by
\begin{align*}
	L(A^E,\mu)
	 = \lim_{n \to \infty} \frac{1}{n} \int \log \|A^E_n(\omega)\| \, d\mu(\omega)  = \inf_{n \ge 1} \frac{1}{n} \int \log \|A^E_n(\omega)\| \, d\mu(\omega).
\end{align*}
For simplicity, we write $L(E)=L(A^E,\mu)$. By Kingman's subaddive ergodic theorem, we have
$$
\lim_{n \to \infty} \frac{1}{n}\log \|A^E_n(\omega)\|=L(E)
$$
for $\mu$-almost every $\omega \in \Omega$. In particular, certain uniform positivity and uniform large deviation estimates (LDT) for the LE are strong indications of Anderson localization, which in its spectral formulation states that for $\mu$-almost every $\omega \in \Omega$, the operator $H_\omega$ has pure point spectrum with exponentially decaying eigenfunctions.

On the other hand, positivity and LDT estimates for the LE are extensively studied topics in dynamical systems. In general, the more random the base dynamics $(\Omega,T,\mu)$ is, the more likely it is that one has positivity and LDT for the LE. For instance, for the well-known Anderson model, where $V_\omega$ is a realization of independent identically distributed random variables, one does have uniform positivity and uniform LDT on any compact set of energies $E$. These are classic results that go back to the seminal work of F\"urstenberg \cite{furstenburg}. Combining this with a certain elimination of double resonance argument, these two properties indeed lead to a localization result for the Anderson model; see, for example, \cite{bucaj2} for recent proofs of all these results mentioned above.

The Anderson model may be put into the context of the present paper as follows. We consider the Anderson model whose single site measure is an atomic measure supported on a finite number of points, which is the most difficult case. Let $\mathcal{A}=\{1, 2, \ldots, \ell\}$ with $\ell\ge 2$ and let $\tilde\mu$ be a fully supported probability measure on $\CA$. Let $\Omega=\CA^\Z$ be the full shift space and consider the left shift $T : \Omega \to \Omega$ defined by $(T \omega)_n = \omega_{n+1}$ for $\omega \in \mathcal{A}^\Z$ and $n \in \Z$. Let $\mu=\tilde\mu^\Z$, which is strongly mixing with respect to $T$. The Anderson model may be generated by setting $V_\omega=f(T^n\omega)$ where $f:\Omega\to\R$ depends only on $\omega_0$. The potentials generated in this way are the most random ones. A natural question is what if the potentials, or rather the base dynamics $(\Omega,T,\mu)$, are less random. In the language of mathematical physics, what can be said if the $V_\omega$'s are weakly correlated? Or in the language of dynamical systems, what if $(\Omega,T,\mu)$ are general mixing systems such as the Arnold cat map or the doubling map? Or more generally, a subshift of finite type with measure of maximal entropy? It turns out that such systems are much more difficult to analyze.

To further explain what this paper accomplishes, we consider a general framework of the base dynamics that includes most of the systems mentioned above as special classes. Let $(\Omega,T)$ be a subshift of finite type. Let $\mu$ be a $T$-ergodic measure that is fully supported on $\Omega$. Moreover, we further assume that $\mu$ admits a local product structure (a detailed definition may be found in Subsection~\ref{sss:local_product}). Let $f:\Omega\to\R$ be $\alpha$-H\"older for some $0<\alpha\le 1$ and non-constant. We define
$$
\CZ_{f}=\{E: L(E)=0\}.
$$
In the present paper, we address the following question.
\begin{problem}
Let $\Omega$, $T$, $\mu$, and $f$ be described as above. How large is $\CZ_{f}$? In particular, when is it discrete, finite, or even empty?
\end{problem}

Note that the discreteness of $\CZ_f$ can be taken as a starting point to show full spectral localization for the corresponding operators; see, for example, the proof of localization in \cite[Proof of Theorem 1.3]{bucaj2}. We comment on this point in more detail in Remark~\ref{r.localizationproof} below. Earlier partial results along this line may be found for example in \cite{CS, BS, BB, bjerklov, DK, SS1, SS, zhang2}, where either the base dynamics or the choice of $f$ are quite restricted, or $\CZ_f$ is still quite large. The main theorem of this paper is:

\begin{theorem}\label{t:posLEawayDiscreteSet}
Suppose $(\Omega, T)$ is a subshift of finite type and $\mu$ is a $T$-ergodic probability measure that has a local product structure. Suppose $T$ has a fixed point and $f$ is H\"older continuous and non-constant. Then the set $\CZ_f$ is discrete.
\end{theorem}

It is clear that we can add a coupling constant $\lambda$ to $f$ in the statement of Theorem~\ref{t:posLEawayDiscreteSet}.  This further indicates that such systems do behave like the Anderson model, as the Anderson model is always localized as long as $\lambda>0$. If we restrict the choice of $f$ to either locally constant or so that $\|f\|_\infty$ is small, then we can improve the result as follows. Let $C^\alpha(\Omega,\R)$, $0<\alpha\le 1$ be the space of $\alpha$-H\"older continuous functions.

\begin{theorem}\label{t:posLEawayFiniteSet}
Let $(\Omega,T,\mu)$ be as in Theorem~\ref{t:posLEawayDiscreteSet}. Suppose $f \in C^\alpha(\Omega,\R)$ is globally bunched or locally constant. Assume further that $f$ is non-constant and $T$ has a fixed point. Then $\CZ_f$ is finite.
\end{theorem}

A detailed definition of global bunching may be found at Subsection~\ref{ss.glundom}. In particular, $f$ is globally bunched if $\|f\|_\infty$ is small. In fact, the smallness of $\|f\|_\infty$ depends only on $\alpha$ and is not of a perturbative nature. A possible explicit choice of a smallness condition on $\|f\|_\infty$ may be found in \eqref{eq:globalbunching4}. We can again add a coupling constant $\lambda$ to $f$ in Theorem~\ref{t:posLEawayFiniteSet} if $f$ is locally constant. If  $f$ is globally bunched, then $\CZ_{\lambda f}$ might become a discrete but not finite set as $\lambda$ becomes large. This is because we will lose global bunching as $\lambda$ becomes large and we have to apply Theorem~\ref{t:posLEawayDiscreteSet} then. In Section~\ref{s:applications}, we shall show that Theorems~\ref{t:posLEawayDiscreteSet} and ~\ref{t:posLEawayFiniteSet} are sharp in the sense that $\CZ_f$ may indeed be nonempty for a suitable locally constant $f$. Thus another natural question is: when can we remove the discrete or finite set $\CZ_f$? We have the following results.

\begin{theorem}\label{t:fullposLE}
Suppose $(\Omega,T)$ is a subshift of finite type and $\mu$ is $T$-ergodic and has a local product structure. Then there is a residual subset $\CG$ of $C^\alpha(\Omega,\R)$ such that for each $f\in\CG$, $\CZ_f$ is empty.
\end{theorem}

Again, if we restrict the choice of $f$ so that it is either locally constant or globally bunched, then we can obtain a uniform lower bound of the $L(E)$ for a even wider class of choices:

\begin{theorem}\label{t:unifPosLE_1}
Suppose $(\Omega,T)$ is a subshift of finite type and $\mu$ is $T$-ergodic and has a local product structure. Consider the subspaces of $C^\alpha(\Omega,\R)$ consisting of globally bunched or locally constant functions. For each of them, there is an open and dense subset $\CG$ such that for every $f \in \CG$, we have $\inf \{ L(E) : E \in \R \} > 0$.
\end{theorem}

Applications of Theorems~\ref{t:posLEawayDiscreteSet}--\ref{t:unifPosLE_1} to more concrete base dynamics such as the doubling map, Arnold's cat map, Markov shifts may be found in Section~\ref{s:applications}.

\begin{remark}\label{r.localizationproof}
(a) Let us emphasize that from the perspective of a spectral analysis of the operator family $\{ H_\omega \}_{\omega \in \Omega}$, and in particular when seeking a proof of spectral localization for this family, the discreteness of $\CZ_f$ is in general the appropriate first milestone towards the eventual goal. It then needs to be combined with control of the Lyapunov exponent away from $\CZ_f$ (the connected components of $\CZ_f^c$ need to be exhausted by intervals on which the Lyapunov exponent is uniformly bounded away from zero; this is often established by proving the continuity of $L(E)$ in $E$ whenever possible), suitable large deviation estimates, and an argument that rules out the presence of infinitely many double resonances for almost every $\omega$. It then follows for $\mu$-almost every $\omega \in \Omega$ that spectrally almost every energy in $\CZ_f^c$ admits an exponentially decaying eigenfunction for $H_\omega$. As the discrete set $\CZ_f$ almost surely carries no weight with respect to the spectral measures of $H_\omega$, this then shows that for $\mu$-almost every $\omega \in \Omega$, the operator $H_\omega$ admits a basis consisting of exponentially decaying eigenfunctions, and the desired spectral localization statement then follows.
\\[2mm]
(b) One is nevertheless interested in obtaining stronger results on the size of the exceptional set $\CZ_f$, such as finiteness or emptiness, whenever possible, as this leads to stronger versions of the dynamical version of an Anderson localization statement. Here, one is interested in showing that the solutions of the time-dependent Schr\"odinger equation $i \partial_t \psi = H_\omega \psi$ are localized. In other words, one seeks to prove good off-diagonal estimates for the matrix elements of $e^{-itH_\omega}$ relative to the standard basis of $\ell^2(\Z)$, uniformly in the time parameter $t$. Energies $E$ in $\CZ_f$ present an obstacle for proving this and one generally simply projects away from these exceptional energies and considers $\chi_I(H_\omega) e^{-itH_\omega}$ with a set $I \subseteq \CZ_f^c$ that has positive distance from $\CZ_f$. In fact, it has been shown that dynamical localization can actually fail, even when spectral localization holds, if one does not project away from $\CZ_f$; compare, for example, \cite{DT, JSS}. Clearly, it is then desirable to show that $\CZ_f$ is empty whenever this can be expected to be true. Of course, as pointed out earlier, this will not always be the case.
\\[2mm]
(c) Let us emphasize that the road map to spectral localization described in part (a) of this remark is applicable in the general setting of ergodic Schr\"odinger operators, and it has been implemented for special cases ranging from the Anderson model to potentials generated by torus translations, the standard skew-shift, the doubling map, or the Arnold cat map. While the literature is vast, let us just mention a few representative papers, \cite{BoGo, BGS, BS, bucaj2}, and refer the reader to \cite{D17, DF1, DF2} for more information. Regarding the base transformations considered in this paper, the absence of a suitable general and global result showing the discreteness of $\CZ_f$ was the primary obstacle in attempting to implement this road map. Thus, the present paper fills precisely this gap and opens the door to a localization proof, which we intend to work out in detail in the second part of this series \cite{ADZ}.
\end{remark}

\subsection{Strategy of Proofs}

One of the main tools we use to prove our results is the so-called \emph{invariance principle} as coined in \cite{AV}. The first version of the invariance principle goes back to Ledrappier \cite{ledrappier} and is later generalized in \cite{AV}. The version we adopt in this paper is due to Bonatti, G\'omez-Mont,  and Viana \cite{BGV, V}. A detailed statement of the invariance principle may be found in Proposition~\ref{p:invPrinciple}. It says that if the Lyapunov exponent $L(A,\mu)$ of a cocycle $(T,A)$ is $0$ and $A$ depends only on the future or the past, then any $(T,A)$-invariant measure $m$ on $\Omega\times \R\bbP^1$ admits a disintegration $\{m_\omega:\omega\in
\Omega\}$ that depends only on the future or the past, respectively.

Another main tool we use is given by the so-called stable and unstable holonomies, which are defined along the stable or unstable sets of  $\omega$, respectively; see Subsection~\ref{sss.holonomies} for a detailed definition. If $L(A,\mu)=0$, we can define a measurable family of stable and unstable holonomies for $\mu$-almost every $\omega$. Then one can use the stable or unstable holonomies to conjuagate the cocycle $(T,A)$ to one that depends only on the future or the past, respectively.

Combining the two steps above, one can show that the family of invariant measures $\{m_\omega\}$ are invariant with respect to the stable and unstable holomies as well. We call such a family an $su$-state.

It turns out that the existence of $su$-states is a very rare event in the sense that they can be easily perturbed away by modifying the data of the cocycle map $A$ at certain periodic points. Roughly speaking, this is how \cite{BGV,BV,V} show the positivity of the Lyapunov exponent for certain typical $C^\alpha$-cocycles. More precisely, \cite{BGV, BV} did it in case the cocycle is fiber bunched or is locally constant while \cite{V} did it for the general case.

However, to prove Theorems~\ref{t:posLEawayDiscreteSet}--\ref{t:posLEawayFiniteSet}, we need to consider Schr\"odinger cocycles with fixed sampling functions. They are basically fixed cocycle maps parametrized by the energy parameter $E \in \R$. So we are not allowed to perturb the cocycle maps to get typicality. Hence, the above strategy is not sufficient to yield the discreteness or finiteness of $\CZ_{f}$ as stated in Section~\ref{ss:results}. It turns out that in addition we need to deploy certain tools from spectral theory. In particular, we will consider the spectra associated with certain periodic orbits and invoke a result from inverse spectral theory for periodic operators. Moreover, to make use of the periodic data, among other things, we also need to show that periodic orbits with small Lyapunov exponent belong to the topological support of the sets where one can define continuous holonomies.  Finally, to use the periodic data to prove the main results, we have to combine the conformal barycenter concept due to Douady and Earle \cite{DE}, Bowen's specification property \cite{bowen}, and Kalinin's theorem regarding approximating $L(E)$ by the Lyapunov exponent along periodic orbits \cite{kalinin}. In short, the proof is based on a fusion of ideas and results from both dynamical systems and spectral theory.

The structure of the remainder of the paper is as follows. In Section~\ref{s:preliminaries}, we state some necessary preliminaries and lay out our context. In Section~\ref{s:ldt}, we give a proof of an additive version of a large deviation estimate for H\"older continuous functions defined on $\Omega$ and for slightly more restricted measures $\mu$. These large deviation estimates may be of independent interest. Moreover, they will play a key role in the second paper of this series \cite{ADZ}. In Section~\ref{s:invPrinciple_Barycenter}, we introduce our main tools such as the invariance principle and the conformal barycenter, and we also give detailed proofs of certain lemmas. We prove Theorems~\ref{t:posLEawayDiscreteSet}--\ref{t:posLEawayFiniteSet} in Section~\ref{s:PosLya1} and Theorems~\ref{t:fullposLE}--\ref{t:unifPosLE_1} in Section~\ref{ss.uniformpositivity}. In Section~\ref{s:applications}, we apply our general Theorems~\ref{t:posLEawayDiscreteSet}--\ref{t:unifPosLE_1} to several concrete models such as the doubling map, Arnold cat map, and Makov chains. In particular, the class of Markov chains includes general locally constant Schr\"odinger potentials defined on the full shift space as a special case, which yields a generalization of the classical F\"urstenberg theorem. Many of the results are the first of their kind. We also compute an explicit choice of $\lambda_0 > 0$ so that $\|f\|_\infty \le \lambda_0$ is sufficient for $f$ to be globally bunched. Finally, we present an example where we show the finite set $\CZ_f$ appearing in the statement of Theorem~\ref{t:posLEawayFiniteSet} may not be removed in general, so that our results are sharp in a suitable sense.

\section{Preliminaries}\label{s:preliminaries}

\subsection{The Setting}

In this section we describe the setting we will work in. We have chosen subshifts of finite type with appropriate ergodic measures as base transformations as a compromise between concreteness and generality. Other possible choices would have been concrete classes of smooth hyperbolic transformations and expanding maps. For background and discussion of the material presented below, we refer the reader to \cite{BGV, BV, V}.

\subsubsection{The Base Space and the Base Transformation} \label{ss:base_space}

Let $\mathcal{A}=\{1, 2, \ldots, \ell\}$ with $\ell \ge 2$ be equipped with the discrete topology. Consider the product space $\mathcal{A}^\Z$, whose topology is generated by the cylinder sets, which are the sets of the form
$$
[n; j_0,\cdots,j_k] = \{ \omega \in \mathcal{A}^\Z : \omega_{n+i} = j_i , \; 0 \le i \le k \}
$$
with $n \in \Z$ and $j_0, \ldots , j_k \in \mathcal{A}$. The topology is metrizable and for definiteness we fix the following metric on $\mathcal{A}^\Z$. Set
$$
N(\omega,\tilde \omega)=\max\{N\ge 0: \omega_n=\tilde \omega_n \mbox{ for all } |n|<N\},
$$
and equip $\mathcal A^\Z$ with the metric $d$ defined by
$$
d(\omega , \tilde \omega) =e^{-N(\omega,\tilde \omega)}.
$$
We consider the left shift $T : \mathcal{A}^\Z \to \mathcal{A}^\Z$ defined by $(T \omega)_n = \omega_{n+1}$ for $\omega \in \mathcal{A}^\Z$ and $n \in \Z$. Let $\mathrm{Orb}(\omega)=\{T^n\omega:\ n\in\Z\}$ be the orbit of $\omega$ under the dynamics $T$.

\begin{definition}
Let $\Omega \subseteq \mathcal{A}^\Z$ be a subshift of finite type and consider the topological dynamical system $(\Omega,T)$.

We say that a finite word $j_0 j_1 \ldots j_k$, where $j_i \in\{1,\ldots, \ell\}$ for $0 \le i \le k$, is \emph{admissible} if it occurs in some $\omega\in\Omega$, that is, there are $\omega\in\Omega$ and $n\in\Z$ such that $\omega_{n+i}=j_i$ for all $0\le i\le k$.

The \textit{local stable set} of a point $\omega \in \Omega$ is defined by
$$
W^s_\mathrm{loc}(\omega) = \{ \tilde \omega \in \Omega : \omega_n = \tilde \omega_n  \text{ for } n \ge 0 \}
$$
and the \textit{local unstable set} of $\omega$ is defined by
$$
W^u_\mathrm{loc}(\omega) = \{ \tilde \omega \in \Omega : \omega_n = \tilde \omega_n \text{ for } n \le 0 \}.
$$

A set is called $s$-\textit{locally saturated} (resp., $u$-\textit{locally saturated}) if it is a union of local stable (resp., local unstable) sets of the form above.

For each $j \in \mathcal{A}$ and each pair of points $\omega, \tilde \omega\in [0;j]$, we denote the unique point in $W^u_{\mathrm{loc}}(\omega)\cap W^s_{\mathrm{loc}}(\tilde\omega)$ by $\omega\wedge\tilde \omega$.
\end{definition}

\subsubsection{Measures With a Local Product Structure}\label{sss:local_product}

Let the subshift $\Omega$ be equipped with the Borel $\sigma$-algebra and let $\mu$ be a probability measure on $\Omega$ that is ergodic with respect to $T$. We define
\begin{align*}
&\Omega^+=\{(\omega_n)_{n\ge 0}: \omega\in\Omega\}\\
&\Omega^-=\{(\omega_n)_{n\le 0}: \omega\in \Omega\}
\end{align*}
to be the spaces of one-sided right and left infinite sequences, respectively, associated with $\Omega$. Metrics for $\Omega^\pm$ can be defined in a way similar to the definition of the metric for $\CA^\Z$ in Subsection~\ref{ss:base_space}. Abusing notation slightly, we still let $d$ denote their metrics. Let $\pi^{+}$ be the projection from $\Omega$ to $\Omega^+$ and $\mu^+=\pi^+_*(\mu)$ be the pushforward measure of $\mu$ on $\Omega^+$. Similarly, we let  $\pi^-$ be the projection to $\Omega^-$ and $\mu^-$ be the pushforward measure on $\Omega^-$. Let $T_+$ be the left shift operator on $\Omega^u$ and $T_-$ be the right shift on $\Omega^-$. For $n\ge 0$, we let $[n;j_0, \ldots, j_k]^+$ denote the cylinder sets in $\Omega^+$; for $n \le -k$, we let $[n;j_0,\ldots, j_k]^-$ denote the cylinder sets in $\Omega^-$. Let $\omega^\pm$ denote points in $\Omega^\pm$, respectively.

For simplicity, for each $1\le j\le \ell$, we set $\mu_j=\mu|_{[0;j]}$. Similarly, we set $\mu^\pm_j=\mu^\pm|_{[0;j]^\pm}$ and $\Omega^\pm_j=\Omega^\pm\cap[0;j]^\pm$, respectively.

Note that we do not have $\Omega=\Omega^-\times\Omega^+$. However, for each $1\le j\le \ell$ we have a natural homeomorphism
$$
P: \Omega_j\to \Omega^-_j\times\Omega^+_j\mbox{ where } P(\omega)=(\pi^-\omega,\pi^+\omega)
$$
Thus, abusing the notation a bit, we may just write $\Omega_j=\Omega^-_j\times \Omega^+_j$. Moreover, we have for all $\omega \in \Omega$,
\beq\label{eq:local_su_set}
(\pi^+)^{-1}(\pi^+\omega)=W^s_{\mathrm{loc}}(\omega),\ (\pi^-)^{-1}(\pi^-\omega)=W^u_{\mathrm{loc}}(\omega),
\eeq

\begin{definition}
We say $\mu$ has a \emph{local product structure} if there is a $\psi:\Omega\to\R_+$ such that for each $1\le j \le \ell$, $\psi\in L^1(\Omega_j,\mu^-_j\times\mu^+_j)$ and
\beq\label{eq:localproduct}
d\mu_j=\psi \cdot d(\mu^+_j\times \mu^-_j).
\eeq
\end{definition}

The local product structure of $\mu$ amounts to saying that $\mu_j^-\times \mu_j^+$ is equivalent to $\mu_j$.  Indeed, \eqref{eq:localproduct} clearly implies that $\mu_j$ is absolutely continuous with respect to $\mu_j^-\times \mu_j^+$. On the other hand, if $\mu_j(E)=0$, then we must have $(\mu_j^-\times \mu_j^+)(E)=0$ since $\psi(\omega)>0$ for all $\omega\in\Omega$. In particular, we may draw the following conclusion. If $E\subset [0;j]$ is $u$-locally saturated with $\mu(E) > 0$ and $F\subset [0;j]$ is $s$-locally saturated with $\mu(F) > 0$, we have
\begin{align*}
 (\mu^-_j\times  \mu^+_j)(E\cap F)
&= (\mu^-_j\times  \mu^+_j)(\pi^-E\times \pi^+ F)\\
 &=\mu^-_j(\pi^-E)\cdot \mu^+_j(\pi^+F)\\
&=\mu(E)\cdot \mu(F)\\
& > 0,
\end{align*}
which implies that
\beq\label{eq:local_product_structure}
\mu(E\cap F)=\mu_j(E\cap F)>0.
\eeq
Conversely, if $\mu^-_j\times \mu^+_j$ is equivalent to $\mu_j$ for each $1\le j\le \ell$, then $ d\mu_j=\psi\cdot d(\mu^-_j\times \mu^+_j)$, where $\psi \in L^1(\Omega_j,\mu^-_j\times\mu^+_j)$ is the Radon-Nikodym derivative of $\mu_j$ with respect to $\mu^-_j\times\mu^+_j$. Note $1/\psi\in L^1(\Omega_j, \mu_j)$ is the Radon-Nikodym derivative of $\mu^-_j\times\mu^+_j$ with respect to $\mu_j$. Hence we must have that $\psi(\omega)>0$ for all $j$ and for $\mu_j$-a.e. $\omega$. We can of course modify $\psi$ so that it's positive everywhere.

\begin{definition}
A \emph{Jacobian} of the measure $\mu^+$ with respect to $T_+$ on $\Omega^+$ is a measurable function $J_+:\Omega^+\to\R_+$ such for each $i \in \{1,\ldots, \ell\}$, we have
$$
d\mu^+(T_+\omega^+)=J_+(\omega^+)\cdot d((T_+)_*(\mu^+|_{[0;i]}))(T_+\omega^+).
$$
A \emph{Jacobian} of $\mu^-$ with respect to $T_-$ can be defined similarly.
\end{definition}

One consequence of the local product structure of $\mu$ is that  $\mu^\pm$ admit Jacobians with respect to $T_\pm$ on $\Omega^\pm_j$ for each $1\le j\le \ell$, respectively. The following lemma is essentially contained in \cite{BV}. While in \cite[Lemma 2.2]{BV}, $\psi$ is assumed to be continuous, we note that the same proof can be applied to obtain the following lemma.

\begin{lemma}\label{p:jacobian_u}
	The measures $\mu^\pm$ admit positive Jacobians $J_\pm\in L^1(\Omega^\pm_j,d\mu^\pm_j)$ with respect to $T_\pm$ on $\Omega^\pm_j$, respectively,  for each $1\le j\le \ell$.
\end{lemma}

For $\vec l=(l_1,\ldots, l_n)\in \{1,\ldots, \ell\}^n$, we write the cylinder $[0;l_1,\ldots, l_n, j]$ as $[0;\vec l, j]$ and set $|\vec l|:=n$. We use a similar notation for spaces of one-sided sequences. For a cylinder $[0;\vec l, j]^+\subset \Omega^+$, we clearly have a Jacobian for $T^{|\vec l|}_+:[0;\vec l, j]^+\to[0;j]$, which is denoted by $J^{(\vec l,j)}_+:[0;l_1,\ldots, l_n, j]^+\to \R_+$ and is given by the formula
$$
J^{(\vec l, j)}_+(\omega^+)=\prod^{n-1}_{k=0}J_+(T_+^k\omega^+).
$$
By the definition of a Jacobian, we have for any integrable function $f:\Omega^+\to\R$ and any $[0;\vec l,j]^+\subset \Omega^+$ that
\beq\label{eq:jacobian}
\int_{[0;j]^+}f(\eta)d\mu^+(\eta)=\int_{[0;\vec l, j]^+}f(T^{|\vec l|}_+\omega^+)J^{(\vec l, j)}_+(\omega^+)d\mu^+(\omega^+).
\eeq

We first have  the following  immediate consequence of Lemma~\ref{p:jacobian_u}, which will be used in Section~\ref{s:PosLya1}.

\begin{coro}\label{c:jacobian}
	Let $D\subset \Omega^+$ be such that $\mu^{+}(D\cap[0;j]^+)>0$. Then for all $[0;\vec l, j]^+\subset \Omega^+$, we have
	$$
	\mu^+(T^{-|\vec l|}_+(D)\cap[0;\vec l,j]^+)>0.
	$$
	Similarly, if $\mu^-(D\cap[0;j]^{-})>0$ for some $D\subset\Omega^-$, then for all $[-|\vec l|; j,\vec l]^-\subset \Omega^-$, we have
	$$
	\mu^-(T^{-|\vec l|}_-(D)\cap [-|\vec l|; j,\vec l]^-)>0.
	$$
\end{coro}

\begin{proof}
	We only consider the case for $(\Omega^+, T_+, \mu^+)$; the case with $(\Omega^-,T_-, \mu^-)$ can be handled similarly.

Without loss of generality, we may just consider a Borel set $D\subset[0;j]^+$ with positive measure. By \eqref{eq:jacobian}, we have
\begin{align*}
0 & < \mu^+(D) \\
&=\int _{[0;j]}\chi_D(\eta)d\mu^+(\eta)\\&=\int_{[0;\vec l, j]}\chi_D (T^{|\vec l|}_+\omega^+)J^{(\vec l, j)}_+(\omega^+)d\mu^+(\omega)^+\\
&=\int_{[0;\vec l,j]\cap (T_+^{-|\vec l|}D)}J^{(\vec l, j)}_+(\omega^+)d\mu^+(\omega^+),
\end{align*}
which implies that $\mu^+\big([0;\vec l,j]\cap (T_+^{-|\vec l|}D)\big)>0.$
\end{proof}

If $\mu$ has a local product structure, then its topological support $\mathrm{supp} \, \mu$ is a subshift of finite type (see, e.g., \cite[Lemma~1.2]{BGV}) and hence, without loss of generality, we will assume throughout that the measure $\mu$ has full support in $\Omega$. Conversely, given any subshift of finite type, the unique equilibrium state associated with a H\"older continuous potential always has a local product structure, see \cite{bowen2, leplaideur} or \cite[Section 2.2]{BV}. In particular, measures with maximal entropy do have a local product structure.

For some results we will need the measure $\mu$ to obey a quantitative version of local product structure, which is defined as follows.

\begin{definition}
We say that $\mu$ satisfies the \textit{bounded distortion property} if there is $C \ge 1$ such that for all cylinders $[n;j_0,\ldots,j_{k}]\subset\Omega$ and $[l;,i_{0},\ \ldots, j_{m}]\subset \Omega$, where $l> n+k$ and $[n;j_0,\ldots,j_{k}]\cap [l;,i_{0},\ldots, i_{m}]\neq\varnothing$, we have
\beq\label{eq:bdd}
C^{-1} \le \frac{\mu \left( [n;j_0,\ldots,j_{k}]\cap [l;i_{0},\ldots, i_{m}] \right)}{\mu \left( [n;j_0,\ldots,j_{k}] \right) \cdot \mu \left( [l;i_{0},\ldots,i_{m}] \right)} \le C.
\eeq
\end{definition}

It is not difficult to see that every measure satisfying the bounded distortion property has a local product structure. Indeed, for every cylinder  $[-k;j_{-k},\ldots,j_{-1},j_{0},\ldots, j_{k}]\subset \Omega$, we have by \eqref{eq:bdd}
\begin{align*}
(\mu^-_{j_0} \times & \mu^+_{j_0}) \big([-k;j_{-k},\ldots,j_{-1},j_{0},\ldots, j_{m}]\big) \\
& = \mu^-_j\big([-k;j_{-k},\ldots,j_{-1},j_{0}]^-\big)\cdot\mu^+_j\big([0;j_{0},\ldots, j_{m}]^+\big)\\
& = \mu\big([-k;j_{-k},\ldots,j_{-1},j_{0}]\big)\cdot\mu\big([0;j_{0},\ldots, j_{m}]\big)\\
& \le \mu \big([-k;j_{-k},\ldots,j_{-1}]\big)\cdot\mu\big([0;j_{0},\ldots, j_{m}]\big)\\
& \le C \mu \big([-k;j_{-k},\ldots,j_{-1},j_{0},\ldots, j_{m}]\big)\\
& = C \mu_{j_0}\big([-k;j_{-k},\ldots,j_{-1},j_{0},\ldots, j_{m}]\big).
\end{align*}
Similarly, we can obtain such estimates for all other cylinders. Since every Borel set can be approximated by cylinder sets, these estimates clearly imply that $\mu^-_j\times\mu^+_j$ is absolutely continuous with respect to $\mu_j$. Note that by $T$-invariance of $\mu$ and by the definition of $\mu^\pm$, $\mu$ has the bounded distortion property if and only if $\mu^+$ or $\mu^-$ has the bounded distortion property. For instance, the bounded distortion property of $\mu^+$ means that for all $n\ge 0$, $l>n+k$, and $[n;j_0,\ldots,j_{k}]^+\cap [l;,i_{0},\ldots, i_{m}]^+\neq\varnothing$, we have
\beq\label{eq:bdd+}
 C^{-1} \le \frac{\mu^+ \left( [n;j_0,\ldots,j_{k}]^+\cap [l;i_{1},\ldots, i_{m}]^+ \right)}{\mu^+ \left( [n;j_0,\ldots,j_{k}]^+ \right) \cdot \mu^+ \left( [l;i_{1},\ldots,i_{m}]^+ \right)} \le C.
\eeq

In fact, given any subshift of finite type, the unique equilibrium state associated with a H\"older continuous potential always has the bounded distortion property; see Lemma~\ref{l:es_to_bdd}.

\subsubsection{$\mathrm{SL}(2,\R)$-Cocycles and Their Projectivization}

A continuous map $A : \Omega \to \mathrm{SL}(2,\R)$ gives rise to the cocycle $(T,A) : \Omega \times \R^2 \to \Omega \times \R^2$, $(\omega , v) \mapsto (T \omega , A(\omega) v)$. For $n \in \Z$, we let $(T,A)^n = (T^n , A_n)$. In particular, we have
$$
A_n(\omega)=
\begin{cases}
A(T^{n-1}\omega) \cdots A(\omega),\ & n\ge 1;\\
I_2, &n=0;\\
[A_{-n}(T^n\omega)]^{-1}, \ &n\le -1,
\end{cases}
$$
where $I_2$ is the identity matrix. Now let $\mu$ be a $T$-ergodic probability measure with topological support equal to $\Omega$. The Lyapunov exponent is given by
\begin{align*}
L(A,\mu)
& = \lim_{n \to \infty} \frac{1}{n} \int \log \|A_n(\omega)\| \, d\mu(\omega) \\
& = \inf_{n \ge 1} \frac{1}{n} \int \log \|A_n(\omega)\| \, d\mu(\omega).
\end{align*}
By Kingman's subaddive ergodic theorem, we have
$$
\lim_{n \to \infty} \frac{1}{n}\log \|A_n(\omega)\| = L(A,\mu)
$$
for $\mu$-a.e. $\omega$. By linearity and invertibility of each $A(\omega)$, we can projectivize the second component and consider $(T,A) : \Omega \times \R\bbP^1 \to \Omega \times \R\bbP^1$.

\subsubsection{Reduction to a Topologically Mixing Subshift}

We need to reduce to the case where $T:\Omega\to\Omega$ is topologically mixing and collect some standard facts. One may find a detailed discussion of the results stated in this section in \cite{katok}.

One says that $(T,\Omega)$ is \emph{topologically mixing} if for any pair of nonempty open sets $U, V\subset\Omega$, there is an $N\ge 1$ such that $T^n(U)\cap V\neq\varnothing$ for all $n\ge N$.

By the spectral decomposition theorem for hyperbolic basic sets, we may decompose $\Omega$ as $\Omega = \bigsqcup^{s}_{l=1}\Omega_l$ for some $s\ge 1$ and for closed subsets $\Omega_l$, so that the following holds true: $T(\Omega_l)=\Omega_{l+1}$ for $1 \le l < s$ and $T(\Omega_s)=(\Omega_1)$, and $T^s|\Omega_l$ is a topologically mixing subshift of finite type for each $1\le l\le s$. The normalized restriction $\mu_l$ of $\mu$ to $\Omega_l$ is a $T^s$-invariant ergodic, fully supported probability measure with local product structure or bounded distortion property, provided the same property is true for $\mu$ on $\Omega$.

Then for a cocycle map $A : \Omega \to \mathrm{SL}(2,\R)$, we consider $A_s : \Omega_l \to \mathrm{SL}(2,\R)$ as $A_s(\omega)$, which may be considered a cocycle map defined over the base dynamics $T^s : \Omega_l \to \Omega_l$. Clearly, $L(A_s, \mu_l) > 0$ for some $1 \le l \le s$ implies that $L(A,\mu) > 0$. Since the present paper is only concerned with the positivity of the Lyapuonv exponent, we assume from now on that $(\Omega, T)$ is topologically mixing.

Note that $\mathrm{supp}(\mu)=\Omega$ and ergodicity of $\mu$ together already imply that $\overline{\mathrm{Orb}(\omega)}=\Omega$ for $\mu$-almost every $\omega\in\Omega$.

Topological mixing has additional consequences, which are needed in the present paper. First, it implies that the set of periodic orbits is dense in $\Omega$. Moreover, we have the following more quantitative behavior of periodic points, which is called the \emph{specification property}. It concerns shadowing finite pieces of segments of different orbits by a single orbit, in particular, by a periodic orbit. It was first introduced by R.~Bowen \cite{bowen}. The following version for subshifts of finite type is due to Weiss \cite{weiss}. For $a < b \in \Z$, we let $[a,b] \subset \Z$ denote the indicated interval of integers. In other words, $[a,b] = \{ n \in \Z : a \le n \le b \}$.

\begin{prop}\label{p:specification}
	Let $(\Omega,T)$ be a topologically mixing subshift of finite type. For each $\epsilon > 0$, there is an integer $r = r(\epsilon)>0$ such that for any choice of points $\omega_i \in \Omega$ and intervals of integers $I_i = [a_i,b_i]$, $i=1,2$, with $a_2 - b_1 > r$ and any $n > b_2 - a_1 + r$, there exists a periodic point $p$ with period $n$ such that
	$$
	d(T^j p, T^j \omega_i) < \epsilon \mbox{ for } j \in I_i,\ i=1,2.
	$$
\end{prop}

Another fact about topologically mixing subshifts of finite type, which can be obtained from Proposition~\ref{p:specification}, is the following: there is $r_0 \in\Z_+$ so that for all $[k;j_0,\ldots, j_n] \subset \Omega$ and all $[l;i_0,\ldots, i_m] \in \Omega$, where $l-(k+n) \ge r_0$, we have
\beq\label{eq:nonempty_intersection}
[k;j_0,\ldots, j_n]\cap[l;i_0,\ldots, i_m]\neq\varnothing.
\eeq
Throughout this paper, for our $(\Omega,T,\mu)$, we let $r_0$ be a number satisfying \eqref{eq:nonempty_intersection}.

\subsubsection{Stable and Unstable Holonomies}\label{sss.holonomies}

Given $(\Omega,T,\mu)$ as above, consider $A : \Omega \to \mathrm{SL}(2,\R)$ and the projective cocycle $(T,A) : \Omega \times \R\bbP^1 \to \Omega \times \R\bbP^1$. We will denote the fiber $\{ \omega \} \times \R\bbP^1$ by $\mathcal{E}_\omega$.

\begin{definition}
A \textit{stable holonomy} $h^s$ for $A$ is a family of homeomorphisms $h^s_{\omega , \omega'} : \mathcal{E}_\omega \to \mathcal{E}_{\omega'}$, defined whenever $\omega$ and $\omega'$ belong to the same local stable set, satisfying the following properties:
\begin{itemize}

\item[(i)] $h^s_{\omega' , \omega''} \circ h^s_{\omega, \omega'} = h^s_{\omega, \omega''}$ and $h^s_{\omega, \omega} = \mathrm{id}$,

\item[(ii)] $A(\omega') \circ h^s_{\omega, \omega'} = h^s_{T \omega, T \omega'} \circ A(\omega)$,

\item[(iii)] $(\omega, \omega') \mapsto h^s_{\omega, \omega'}(\phi)$ is continuous when $\omega, \omega'$ belong to the same local stable set, uniformly in $\phi$.

\end{itemize}
An \textit{unstable holonomy} $h^u_{\omega , \omega'} : \mathcal{E}_\omega \to \mathcal{E}_{\omega'}$ is defined analogously for pairs of points in the same unstable set.
\end{definition}

By property (i), we have $h^\tau_{\omega,\omega'}=(h^\tau_{\omega',\omega})^{-1}$ for any $\omega'\in W^\tau_\loc(\omega)$, where $\tau\in\{s,u\}$.

These projective holonomies $h^s_{\omega, \omega'}, h^u_{\omega, \omega'}$ typically arise via projectivization of $H^s_{\omega, \omega'}, H^u_{\omega, \omega'} \in \mathrm{SL}(2,\R)$ (for suitable pairs $(\omega,\omega')$) that are obtained as follows,
\begin{equation}\label{e.holonomiesdef}
H^s_{\omega,\omega'} = \lim_{n \to \infty} A_n(\omega')^{-1} A_n(\omega), \quad H^u_{\omega,\omega'} = \lim_{n \to \infty} A_{-n}(\omega')^{-1} A_{-n}(\omega)
\end{equation}
for $\omega , \omega'$ in the same stable (resp., unstable) set. Conditions need to be placed on the cocycle to ensure convergence in \eqref{e.holonomiesdef}; see, for example, the proof of Lemma~\ref{p.bgvprop3}. The analogues of the properties (i)--(iii) for $H^s_{\omega, \omega'}, H^u_{\omega, \omega'}$ follow directly from the construction and this in turn implies (i)--(iii) for $h^s_{\omega, \omega'}, h^u_{\omega, \omega'}$ by projectivization. Holonomies that arise from \eqref{e.holonomiesdef} are called \emph{canonical holonomies} of $A$.

\subsubsection{Invariant Measures of Projective Cocycles}

Consider a projective cocycle $(T,A) : \Omega \times \R\bbP^1 \to \Omega \times \R\bbP^1$ that has stable and unstable holonomies.

\begin{definition}
Suppose we are given a $(T,A)$-invariant probability measure $m$ on $\Omega \times \R\bbP^1$ that projects to $\mu$ in the first component. A \emph{disintegration} of $m$ along the fibers is a measurable family $\{m_\omega: \omega\in\Omega\}$ of conditional probabilities on $\R\bbP^1$ such that $m = \int m_\omega \, d\mu(\omega)$, that is,
$$
m(D)=\int_\Omega m_\omega(\{z\in\R\bbP^1:(\omega,z)\in D\}) \, d\mu(\omega)
$$
for each measurable set $D\subset \Omega\times \R\bbP^1$.
\end{definition}

By Rokhlin's disintegration theorem, such a disintegration exists. Moreover, $\{\tilde m_\omega:\omega\in\Omega\}$ is another disintegration of $m$ if and only if $m_\omega=\tilde m_\omega$ for $\mu$-almost every $\omega\in\Omega$. By a straightforward calculation one checks that $\{ A(\omega)_* m_{\omega} : \omega \in \Omega\}$ is a disintegration of $(T,A)_*m$. In particular, the facts above imply that $m$ is $(T,A)$-invariant if and only if $A(\omega)_* m_\omega = m_{T\omega}$ for $\mu$-almost every $\omega\in\Omega$.

Such a measure $m$ will be called an $s$-\emph{state} (resp., a $u$-\emph{state}) if it is in addition invariant under the stable (resp., unstable) holonomies. That is, the disintegration $\{m_\omega:\omega\in\Omega\}$ satisfies that $(h^s_{\omega,\omega'})_* m_\omega = m_{\omega'}$ for $\mu$-almost every $\omega\in\Omega$ and for every $\omega'\in W^s_\loc(\omega)$ (resp., $(h^u_{\omega,\omega'})_* m_\omega = m_{\omega'}$ for $\mu$-almost every $\omega\in\Omega$ and for every $\omega'\in W^u_\loc(\omega)$). In this case, we say that $\{m_\omega\}$ is $s$-\emph{invariant} (resp. $u$-\emph{invariant}). A measure that is both an $s$-state and a $u$-state is called an $su$-\emph{state}.

\subsubsection{Schr\"odinger Operators and Cocycles}

In this subsection let us initially assume that $\Omega$ is a compact metric space, $T:\Omega\to\Omega$ is a homeomorphism, and $f : \Omega \to \R$ is continuous. We consider potentials $V_\omega : \Z \to \R$ defined by $V_\omega(n) = f(T^n \omega)$ for $\omega \in \Omega$ and $n \in \Z$, and associated Schr\"odinger operators $H_\omega$ in $\ell^2(\Z)$ acting by
$$
[H_\omega \psi](n) = \psi(n+1) + \psi(n-1) + V_\omega(n) \psi(n).
$$
The spectrum $\sigma(H_\omega)$ is defined as
$$
\sigma(H_\omega)=\{E\in\C: H_\omega-E \mbox{ does not have a bounded inverse}\}.
$$

For a subset $S$ of a metric space $(X,d)$ and $\delta > 0$, the open $\delta$-neighborhood of $S$ is given by $B_\delta(S)=\{x\in X: d(x,s)<\delta\mbox{ for some } s\in S\}$. In particular, $B_\delta (x)$ denotes the open ball centered at the point $x\in X$. We need the following uniform estimate that relates the spectrum $\sigma(H_\omega)$ with the orbit $\mathrm{Orb}(\omega)=\{T^n(\omega),\ n\in\Z\}$; see, for example, \cite[Theorem 6]{zhang}.

\begin{prop}\label{p:denseOrb_to_denseSpec}
For each $\varepsilon>0$, there exists a $\delta>0$, depending on $\varepsilon$ only, so that the following holds true. If the orbit $\mathrm{Orb}(\omega_0)$ of some $\omega_0\in \Omega$ satisfies
	$$
	\mathrm{Orb}(\omega_0)\cap B_\delta(\omega)\neq\varnothing
	$$
	for some $\omega\in\Omega$, then
	$$
	\sigma(H_\omega)\subset B_\varepsilon[\sigma(H_{\omega_0})].
	$$
\end{prop}

Proposition~\ref{p:denseOrb_to_denseSpec} implies that if $\mathrm{Orb}(\omega_0)$ is dense in $\Omega$, then $\sigma(H_\omega) \subseteq \sigma(H_{\omega_0})$ for all $\omega \in \Omega$. In this case, we set
$$
\Sigma=\sigma(H_{\omega_0}).
$$

Let us now return to the main scenario of this paper, where $T$ is a topologically mixing shift operator on a subshift of finite type $\Omega$ with an ergodic measure $\mu$ satisfying $\mathrm{supp}(\mu) = \Omega$. Let $\mathrm{Per(T)}$ be the set of periodic points of $T$. Recall that $\overline{\mathrm{Per(T)}} = \Omega$. Recall we have that $\overline{\mathrm{Orb}(\omega)} = \Omega$ for $\mu$-almost every $\omega$. All these facts together with Proposition~\ref{p:denseOrb_to_denseSpec} imply for $\mu$-almost every $\omega$ that
\beq\label{eq:spectrum_variousForms}
\Sigma = \sigma(H_\omega) = \overline{\bigcup_{\omega_p \in \mathrm{Per(T)}}\sigma(H_{\omega_p})}.
\eeq

Spectral properties of the operators $H_\omega$ can be investigated in terms of the behavior of the solutions to the difference equation
\begin{equation}\label{e.eve}
u(n+1) + u(n-1) + V_\omega(n) u(n) = E u(n), \quad n \in \Z,
\end{equation}
with $E$ real or complex (depending on the problem in question). These solutions in turn can be described with the help of the Schr\"odinger cocycle $(T,A^E)$ with the cocycle map $A^{E}:\Omega\to\SL(2,\R)$ (resp., $\SL(2,\C)$ when $E \in \C \setminus \R$) being defined as
$$
A^E(\omega) = A^{(E-f)}(\omega):= \begin{pmatrix} E-f(\omega) & -1 \\ 1 & 0 \end{pmatrix},
$$
where we often leave the dependence on $f:\Omega \to \R$ implicit as it will be fixed most of the time. Such cocycles describe the transfer matrices associated with Schr\"odinger operators with dynamically defined potentials. Specifically, $u$ solves \eqref{e.eve} if and only if
$$
\begin{pmatrix} u(n) \\ u(n-1) \end{pmatrix} = A^{E}_n(\omega) \begin{pmatrix} u(0) \\ u(-1) \end{pmatrix}, \quad n \in \Z.
$$

For the Schr\"odinger cocycle $(T,A^E)$, we set $L(E)=L(A^E,\mu)$. One of the main questions in the spectral analysis of the ergodic family of Schr\"odinger operators $\{H_\omega\}_{\omega\in\Omega}$ (with respect to the ergodic measure $\mu$) is for how many $E \in \Sigma$ we have $L(E) > 0$.

\subsection{Periodic Potentials}\label{s:periodic}

A periodic point $\omega$ of $T$ gives rise to a periodic potential, that is, if $T^p \omega = \omega$, then, $V_{\omega}(n + p) = V_{\omega}(n)$ for every $n \in \Z$. Since much of our work below will involve the study of periodic points and the associated potentials, let us recall some basic properties of Schr\"odinger operators with periodic potentials; see \cite{S} for proofs of the results stated in this subsection.

Consider a Schr\"odinger operator
$$
[H \psi](n) = \psi(n+1) + \psi(n-1) + V(n) \psi(n).
$$
in $\ell^2(\Z)$ with a $p$-periodic potential, $V(n + p) = V(n)$ for every $n \in \Z$. Define, for $E \in \C$, the \textit{monodromy matrix}
$$
M(E) =\prod^{0}_{j=p-1} \begin{pmatrix} E - V(j) & -1 \\ 1 & 0 \end{pmatrix}
$$
and the \textit{discriminant} $\Delta(E) = \tr(M(E))$,  where $\tr(B)$ is the trace of $B$. The function $\Delta(\cdot)$ is a monic polynomial of degree $p$.

\begin{prop}\label{p:spectrumbands}
The set $\Delta^{-1}((-2,2))$ consists of $p$ disjoint open intervals and on each of them, $\Delta$ is strictly monotone. Moreover, $\sigma(H)  = \overline{\Delta^{-1}((-2,2))} = \Delta^{-1}([-2,2])$.
\end{prop}

This shows that the spectrum of $H$ consists of a finite union of closed intervals and, in fact, the number of connected components of the spectrum is bounded by the period of the potential. This suggests an interesting inverse problem. Suppose we are given a set that has such a form, that is, it has finitely many connected components, each being a closed interval. Suppose further that we know that the set is the spectrum of a periodic Schr\"odinger operator. Can we say anything about the period of the potential? \footnote{The more natural inverse problem would lead us too far afield. Namely, one should rather ask, given a finite union of closed intervals, identify within a suitable class of operators those that have the given set as their spectrum. The theory is classical but one needs to pass to the more general class of finite-gap Jacobi matrices to study this question in the appropriate setting.}

\begin{prop}
Suppose $V : \Z \to \R$ is periodic. Denote the spectrum of the associated Schr\"odinger operator by $\sigma$.\\
{\rm (a)} For a probability measure $m$ on $\sigma$, consider its potential energy
$$
\mathcal{E}(m) = \iint \log\left( |E - E'|^{-1} \right) \, dm(E) \, dm(E') \in \R \cup \{ \infty \}.
$$
Then there is a unique measure, $m_\sigma$, which minimizes the potential energy among all probability measures on $\sigma$, and in fact $\mathcal{E} (m_\sigma) = 0$.\\
{\rm (b)} The measure $m_\sigma$ assigns rational weight to each connected component of $\sigma$.\\
{\rm (c)} The potential $V$ is $p$-periodic if and only if the weight of each connected component of $\sigma$ with respect to $m_\sigma$ is an integer multiple of $\frac{1}{p}$.
\end{prop}

This result shows that the shape of the spectrum of a periodic Schr\"odinger operator determines the period of the potential. An immediate consequence is the fact that the spectrum of a periodic Schr\"odinger operator is connected if and only if the period is one, that is, the potential is constant. Another characterization of constant potentials is the following:

\begin{prop}\label{p.perinvmeas}
Suppose $V : \Z \to \R$ is periodic. Then the spectrum $\sigma$ of the associated Schr\"odinger operator has Lebesgue measure at most $4$. Moreover, the Lebesgue measure of $\sigma$ is equal to $4$ if and only if $V$ is constant.
\end{prop}

Finally, we note the following standard facts. For each $E\in\C $ such that $\Delta(E)\neq\pm2$, there are exactly two eigendirections $s(E)$ and $u(E)$ in $\C\bbP^1$ of the monodromy matrix $M(E)$, which are actually the so-called Weyl-Titchmarsh $m$-functions associated with the operator. Moreover, $s(E)\neq u(E)$ are real if and only if $E\in \R\setminus \sigma(H_V)$, and they are the stable and unstable directions of the real hyperbolic matrix $M(E)$. Here we always set $s(E)$ to be the stable direction and $u(E)$ to be the unstable direction. If $E$ is in the upper or lower-half plane or is such that $E\in\R$ and $|\Delta(E)|<2$, then $s(E)$ and $u(E)$ are not real. In the latter case, we have $s(E)=\overline{u(E)}$. For $\Delta(E) = \pm2$, we let $I \subseteq \R$ be the connected component of $\sigma(H_V)$ containing $E$. If $E$ belongs to the boundary of $I$, then $M(E)$ has a unique real invariant direction. We may think of this case as $s(E) = u(E)$. If $E$ is a point at which a spectral gap is collapsed (or, in other words, at which two different components of $\Delta^{-1}(-2,2)$ touch), then $M(E) = \pm I_2$, in which case all directions are invariant.

Based on the description above, we may consider two functions $s$ and $u$ which are holomorphic on the upper or lower half plane $\H$ and $\C\setminus\overline{\H}$, respectively. When restricted to the real line $\R$, they both are continuous functions. Moreover, they are analytic on each spectral gap or in the interior of each connected component of $\sigma(H_V)$. If $E_0$ is on the boundary of some connected component of $I\subseteq \sigma(H_V)$, then $s$ and $u$ are locally like $g\big(\sqrt{\pm(E-E_0)}\big)$ near $E_0$ for some choice of $g$ that is real-analytic near $E_0$. Here the choice of $g$ depends on $s$ or $u$, and the sign of $(E-E_0)$ is determined by whether $E_0$ is the right or left endpoint of $I$. Moreover, $s(E)$ and $u(E)$ are real only when $\sqrt{\pm(E-E_0)}$ is real. Thus, we can find an open disk $D \subseteq \C$ centered at $E_0$ and a ramified (at $E_0$) double cover $\pi : \tilde D \to D$ of $D$ so that $s(\tilde E)$ and $u(\tilde E)$ are holomorphic in $\tilde E \in \tilde D$. Moreover, when $\pi(\tilde E) \in D \cap \R$, $s(\tilde E)$ and $u(\tilde E)$ are real only when $\sqrt{\pm(\pi(\tilde E)-E_0)}$ is real.

\section{Large Deviations}\label{s:ldt}

The main goal of  this section is to prove the following large deviation theorem. Let $C^\alpha(\Omega,\R)$, $0< \alpha\le 1$,  be the space of  $\alpha$-H\"older continuous functions. In other words, $f\in C^\alpha(\Omega,\R)$ if there are $C>0$ such that
$$
|f(\omega)-f(\omega')|<C\cdot d(\omega,\omega')^\alpha \mbox{ for all }\omega,\omega'\in\Omega.
$$
Note that here $C^1(\Omega,\R)$ is the space of Lipschitz continuous functions, not the space of functions with continuous derivatives. Similarly, we can define the space $C^\alpha(\Omega^+,\R)$. Throughout this section, $\mu$, or equivalently $\mu^+$, will be assumed to have the bounded distortion property.

\begin{theorem}\label{t:ldt1}
	Let $(\Omega,T)$ be a topologically mixing subshift of finite type. Let $\mu$ be a $T$-ergodic probability measure that has the bounded distortion property. Let $f \in C^\alpha(\Omega,\R)$ for some $0 < \alpha \le 1$. Then, for each $\varepsilon>0$, there exist $C, c > 0$, depending on $f,\alpha$, and $\varepsilon$, such that
$$
  \mu \bigg\{ \omega \in \Omega : \bigg| \frac{1}{n} \sum^{n-1}_{k=0} f(T^k\omega) - \int_\Omega f \, d\mu \bigg| \ge \varepsilon\bigg\} < Ce^{-cn},\ \forall n\ge 1.
$$
\end{theorem}

Theorem~\ref{t:ldt1} will be a consequence of the following version of large deviations. Recall we have the spaces $(\Omega^\pm,T_\pm,\mu^\pm)$ of one-sided infinite sequences with nonnegative/nonpositive indices.

\begin{theorem}\label{t:ldt2}
Let $(\Omega^+,T_+,\mu^+)$ be a topologically mixing one-sided subshift of finite type and suppose that $\mu^+$ is $T^+$-ergodic and has the bounded distortion property. Let $f\in C^\alpha(\Omega^+,\R)$ for some $0<\alpha\le 1$. Then for each $\varepsilon>0$, there exist $C,c>0$, depending on $f$, $\alpha$, and $\varepsilon$ such that
$$	 \mu^+\bigg\{\omega^+\in\Omega^+:\bigg|\frac{1}{n}\sum^{n-1}_{k=0}f(T_+^k\omega^+)-\int_{\Omega^+} f d\mu^+\bigg|\ge \varepsilon\bigg\}<Ce^{-cn},\ \forall n\ge 1.
$$
\end{theorem}

We first derive Theorem~\ref{t:ldt1} from Theorem~\ref{t:ldt2}. Let us write $S_nf:=\sum^{n-1}_{k=0}f\circ T^k$ for the Birkhoff sums.

\begin{proof}[Proof of Theorem~\ref{t:ldt1}]
Let $f\in C^\alpha(\Omega,\R)$. For each $1\le i\le \ell$, fix a choice of $\omega^{(i)}\in[0;i]$. Define $\varphi(\omega) = \omega^{(\omega_0)}\wedge \omega$, which is continuous and constant on $W^s_\loc(\omega)$ for every $\omega\in\Omega$. Since $f$ is H\"older  continuous and $\varphi(\omega)\in W^s_\loc(\omega)$, a straightforward computation shows that
$$
h^s(\omega) := \sum^{\infty}_{n=0} \big[ f(T^n\omega) - f(T^n\varphi(\omega)) \big]
$$
converges uniformly, and hence is continuous. We define
$$
f^+(\omega):=f(\omega)+h^s(T\omega)-h^s(\omega).
$$
Clearly, $f^+$ is continuous and cohomologous to $f$. Thus, we have
$$
\int_\Omega f \, d\mu = \int_\Omega f^+ \, d\mu \mbox{ and } \big\| S_n f - S_n f^+ \big\|_\infty < \frac{2\|h^s\|_\infty}n,
$$
where $\|\cdot\|_\infty$ denotes the supremum norm. It is straightforward to see that
$$  	
f^+(\omega) = f(\varphi(\omega)) + \sum^{\infty}_{n=0} \big[ f(T^nT\varphi(\omega)) - f(T^n\varphi(T\omega)) \big],
$$
which implies that $f^+$ is constant on $W^s_\loc(\omega)$ for all $\omega \in \Omega$. Moreover, we claim that $f^+ \in C^{\frac\alpha2}(\Omega^+,\R)$. Indeed, take $\omega$ and $\omega' \in \Omega$. Without loss of generality, we may assume $N(\omega,\omega')$ is large and take $k = \lfloor\frac N2\rfloor$. Then we have
\begin{align*}  	
f^+(\omega) - f^+(\omega') = \sum^{k}_{n=0} \big[ f(T^n\varphi(\omega)) - f(T^n\varphi(\omega') \big] + \sum^{k-1}_{n=0} \big[ f(T^n\varphi(T\omega')) - f(T^n\varphi(T\omega') \big] \\
+ \sum^{\infty}_{n=k} \big[ f(T^nT\varphi(\omega)) - f(T^n\varphi(T\omega)) \big] + \sum^{\infty}_{n=k} \big[ f(T^nT\varphi(\omega')) - f(T^n\varphi(T\omega')) \big],
\end{align*}
where the absolute values of the first two terms may be bounded by
$$
C \sum^{k}_{i=1} e^{-\alpha(N-i)} \le C e^{-\alpha \frac N2} = C d(\omega,\omega')^{\frac\alpha2},
$$
and the absolute values of the last two terms may be bounded by
$$
C e^{-\alpha k} \le C e^{-\alpha\frac N2} = C d(\omega,\omega')^{\frac\alpha{2}}.
$$
Thus, $f^+$ descends to a function in $C^{\frac\alpha2}(\Omega^+,\R)$. Abusing notation slightly, let $f^+$ denote its descended function as well. Clearly, we have $\int_\Omega f^+ \, d\mu = \int_{\Omega^+} f^+ \, d\mu^+$ and $S_n f^+(\omega) = S_n f^+(\pi^+\omega)$. Fix any $\varepsilon > 0$ and define
$$  	
\CB^+_n(\varepsilon) := \bigg\{ \omega^+ \in \Omega^+ : \bigg| \frac{1}{n} S_n f^+(\omega^+) - \int_{\Omega^+} f^+ \, d\mu^+ \bigg| > \varepsilon \bigg\}.
$$
By Theorem~\ref{t:ldt2}, there are $C, c > 0$, depending on $f^+$, $\alpha$, and $\varepsilon$, such that
$$ 	
\mu^+(\CB^+_n) < C e^{-cn},\ \forall n\ge 0.
$$
Combining the relations of $f$ and $f^+$ above, there exists a $N = N(\varepsilon)$ such that
$$
\left\{ \omega \in \Omega : \bigg| \frac{1}{n} S_n f(\omega) - \int_\Omega f \, d\mu \bigg| > \varepsilon \right\} \subseteq (\pi^+)^{-1} \CB^+_n(\varepsilon),\ \forall n \ge N.
$$
Changing $C, c$ if necessary, we then have for all $n\ge 1$,
\begin{align*}
\mu \left\{ \omega \in \Omega : \bigg| \frac{1}{n} S_n f(\omega) - \int_\Omega f \, d\mu \bigg| > \varepsilon \right\} & \le \mu[ (\pi^+)^{-1} \CB^+_n(\varepsilon)] \\
& = \mu^+(\CB^+_n(\varepsilon)) \\
& < C e^{-cn},
\end{align*}
as desired.
\end{proof}

To prove Theorem~\ref{t:ldt2}, we first need the following lemma. For $\vec l=(l_1,\ldots, l_n)$ where $l_1\cdots l_n$ is admissible (in this case we also just say that $\vec l$ is admissible), we set $\Omega^+_{\vec l}:=[0;l_1,l_2,\ldots,l_n]^+$, $|\vec l|:=n$, and
$$
\mu^+_{\vec l}=\frac1{\mu^+\big(\Omega^+_{\vec l}\big)}T^{|\vec l|+r_0}_*\mu^+\big|_{\Omega^+_{\vec l}},
$$
where $r_0$ is from \eqref{eq:nonempty_intersection}. In other words, $\mu^+_{\vec l}$ is the normalized push-forward of $\mu^+$ under the homeomorphism $T^{|\vec l|+r_0}:\Omega^+_{\vec l}\to \Omega^+$.  Recall we also write $[n;\vec l]^+=[n;\l_1,\ldots, l_n]^+$.

\begin{lemma}\label{l:bdd+}
Consider a topologically mixing one-sided subshift of finite type $(\Omega^+,T_+,\mu^+)$, where $\mu^+$ has the bounded distortion property. There exists a $C \ge 1$ so that, uniformly for all admissible $\vec l$, we have
\beq\label{eq:bdd+1}
C^{-1}\le \frac{d\mu^+_{\vec l}}{d\mu^+}(\omega^+)\le C \mbox{ for $\mu$-a.e. } \omega^+,
\eeq
where $\frac{d\mu^+_{\vec l}}{d\mu^+}$ is the Radon-Nikodym derivative of $\mu^+_{\vec l}$ with respect to $\mu^+$. In particular, we have for all positive measurable functions $f$ and all admissible $\vec l$,
\beq\label{eq:bdd+2}
C^{-1}\int fd\mu^+\le \int fd\mu^+_{\vec l}\le  C\int fd\mu^+.
\eeq
\end{lemma}

\begin{proof}
Fix an admissible $\vec l=(\l_1,\ldots, \l_n)$. Clearly, \eqref{eq:bdd+1} is equivalent to the existence of a $ C\ge 1$, independent of  $\vec l$, such that for every $[n;\vec i]^+=[n;i_1,\ldots, i_m]^+ \subseteq \Omega^+$ (which implies $n\ge 0$), we have
\beq\label{eq:bdd+3}
C^{-1}\le \frac{\mu^+_{\vec l}([n;\vec i])}{\mu^+([n;\vec i])}\le  C.
\eeq
By  definition of $\mu^+_{\vec l}$, we have
\begin{align*}
\frac{\mu^+_{\vec l}([n;\vec i])}{\mu^+([n;\vec i])}=\frac{\mu^+([0;\vec l]\cap [n+|\vec l|+r_0;\vec i])}{\mu^+([0;\vec l])\mu^+([n;\vec i])}.
\end{align*}
Hence, \eqref{eq:bdd+3} is equivalent to
$$
C^{-1}\le\frac{\mu^+([0;\vec l]\cap [n+|\vec l|+r_0;\vec i])}{\mu^+([0;\vec l])\mu^+([n;\vec i])} \le C.
$$
By $T^+$-invariance of $\mu^+$, the above estimate is then guaranteed by \eqref{eq:bdd+}, as long as $[0;\vec l] \cap [n+|\vec l|+r_0;\vec i] \neq \varnothing$. By definition of $r_0$ from \eqref{eq:nonempty_intersection} and the fact that $n\ge 0$, we indeed have $[0;\vec l]\cap [n+|\vec l|+r_0;\vec i]\neq\varnothing$.
\end{proof}
	
We are now ready to prove Theorem~\ref{t:ldt2}. We adopt the strategy of \cite[Section 6.1]{AD}.

\begin{proof}[Proof of Theorem~\ref{t:ldt2}]
For simplicity, we write $I_n(\omega)=\frac1n S_nf(\omega)$ for the Birkhoff averages. Fix any $\gamma>\int_{\Omega^+} f \, d\mu^+$. By $T^+$ invariance, we have for all $N>0$ that
\beq\label{eq:small_average}
\int_{\Omega^+} \big(I_N(\omega^+)-\gamma \big) \, d\mu^+ = -\kappa' := \int_{\Omega^+} f \, d\mu^+ - \gamma < 0.
\eeq
By the Birkhoff Ergodic Theorem, $I_n(\omega)$ converges to $\int f \, d\mu^+$ pointwise almost everywhere and in $L^1$. Thus we have
$$
\lim_{n\to\infty}\mu^+\{\omega^+: I_n(\omega^+)\ge \gamma\}=0.
$$
By \eqref{eq:bdd+1} of Lemma~\ref{l:bdd+} and by choosing $\vec l$ admissible, we have
\beq\label{eq:small_average1}
\lim_{n\to\infty}\sup_{\vec l}\mu^+_{\vec l}\{\omega^+: I_n(\omega^+)\ge \gamma\}=0.
\eeq
Set $\CB_{n,\vec l}=\{\omega^+: I_n(\omega^+)\ge \gamma\}$. Combine \eqref{eq:bdd+2}, \eqref{eq:small_average} and \eqref{eq:small_average1}, we have for all large $n$
$$
\sup_{\vec l}\int_{\CB^\complement_{n,\vec l}}(I_n(\omega^+)-\gamma)d\mu^+_{\vec l} < \frac1{2C} \int_{\Omega^+}(I_n(\omega^+)-\gamma) \, d\mu.
$$
Since $\int_{\Omega^+} (I_n(\omega^+)-\gamma)d\mu^+_{\vec l}=\int_{\CB_{n,\vec l}}+\int_{\CB^\complement_{n,\vec l}}(I_n(\omega^+)-\gamma)d\mu^+_{\vec l}$, the estimate above and \eqref{eq:small_average} together imply for all large $n$,
\begin{align}\label{eq:bdd_average}
\sup_{\vec l}\int_{\Omega^+} (I_n(\omega^+)-\gamma) \, d\mu^+_{\vec l}\le -(C')^{-1}\kappa',
\end{align}
where $C'$ may be taken as $4C$. Fix a large $N$ so that \eqref{eq:bdd_average} holds true. There clearly exist $\delta'>0$ and $C'>0$ such that for all $|t|<\delta'$ and all $1\le n\le N$, we have
$$
\sup_{\vec l}\int_{\Omega^+}e^{tn(I_nf(\omega^+)-\gamma)} \, d\mu^+_{\vec l}\le \widetilde C.
$$
Hence $\phi_{n, \vec l}(t):=\int e^{t(S_nf(\omega^+)-\gamma)}d\mu^+_{\vec l}$, $1\le n\le N$ and $\vec l$ admissible, are uniformly bounded holomorphic functions on $\{t\in \C: |t|<\delta'\}$. Note that for all $1\le n\le N$,
$$
\phi_{n,\vec l}(0)=1\mbox{ and }\phi'_{n,\vec l}(0)=\int_{\Omega^+} n(I_n(\omega^+)-\gamma) \, d\mu^+_{\vec l}.
$$
By \eqref{eq:bdd_average}, we have for all $|t|<\delta'$ (change $\delta'$ if necessary) that
$$
\sup_{\vec l}\big\{\ln\phi_{N, \vec l}(t)\big\}\le -(C')^{-1}N\kappa't.
$$
Hence, for all $0 < \delta \le \frac{\delta'}2$ so that
\beq\label{eq:small_exp1}
\sup_{\vec l} \int e^{\delta N(I_Nf(\omega^+)-\gamma)} \, d\mu^+_{\vec l} < e^{-(C')^{-1} \delta \kappa' N}.
\eeq
We have obtained the large deviation estimate for all $n=N$. Now we want to extend it to all $n\ge 1$ via the bounded distortion property of $\mu^+$.
		
Since $f \in C^\alpha(\Omega^+,\R)$, it is straightforward to see that
\beq\label{eq:small_dist}
|S_{n+r_0} f(\omega^+) - S_{n+r_0} f(\tilde \omega^+)| \le C \sum^{n+r_0-1}_{k=0} d(T_+^k\omega^+,T_+^k\tilde\omega^+)^\alpha \le C_1 r_0,
\eeq
provided $\omega^+,\tilde\omega^+\in \Omega^+_{\vec l}$, where $|\vec l|=n$. Note here that $C_1$ depends only on $\alpha$ and $f$. We choose $\omega^+_{\vec l,\max}, \omega^+_{\vec l,\min}\in\Omega^+_{\vec l}$ so that
$$
I_{n+r_0}(\omega^+_{\vec l, \max})=\max_{\omega^+\in\Omega^+_{\vec l}}\{I_{n+r_0}(\omega^+)\}\mbox{ and } I_{n+r_0}(\omega^+_{\vec l, \min})=\min_{\omega^+\in\Omega^+_{\vec l}}\{I_{n+r_0}(\omega^+)\}.
$$
In particular we have
\beq\label{eq:small_dist1}
(n+r_0)\big(I_{n+r_0}(\omega^+_{\vec l, \max})-I_{n+r_0}(\omega^+_{\vec l, \min})\big) < C_1 r_0.
\eeq
Since $S_{n+m} f = S_m (f\circ T^n) + S_n f$, we have for all $m, n\in\Z_+$ that
\begin{align*}\label{eq:additivity}
\int_{\Omega^+} & e^{\delta (n+r_0+m)(I_{n+r_0+m}(\omega^+)-\gamma)} \, d\mu^+ \\
& = \int e^{\delta (n+r_0)(I_{n+r_0}(\omega^+)-\gamma)} e^{\delta m(I_{m}(T^{n+r_0}\omega^+)-\gamma)} \, d\mu^+ \\
& = \sum_{|\vec l|=n} \int_{\Omega^+_{\vec l}} e^{\delta (n+r_0)(I_{n+r_0}(\omega^+)-\gamma)} e^{\delta m(I_{m}(T^{n+r_0}\omega^+)-\gamma)} \, d\mu^+ \\
& \le \sum_{|\vec l|=n} e^{\delta (n+r_0)(I_{n+r_0}(\omega^+_{\vec l,\max})-\gamma)} \int_{\Omega^+_{\vec l}} e^{\delta m(I_{m}(T^{n+r_0}\omega^+)-\gamma)} \, d\mu^+ \\
& = \sum_{|\vec l|=n} \mu^+(\Omega^+_{\vec l}) e^{\delta (n+r_0)(I_{n+r_0}(\omega^+_{\vec l,\max})-\gamma)} \int_{\Omega^+_{\vec l}} e^{\delta m(I_{m}(\omega^+)-\gamma)} \, d\mu^+_{\vec l} \\
& \le \bigg(\sup_{\vec l}\int_{\Omega^+_{\vec l}}e^{\delta m(I_{m}(\omega^+)-\gamma)} \, d\mu^+_{\vec l}\bigg) \\
& \qquad \cdot \sum_{|\vec l|=n} \bigg(e^{\delta (n+r_0)[I(\omega^+_{\vec l,\max})-I(\omega^+_{\vec l,\min})]} \cdot \int_{\Omega^+_{\vec l}}e^{\delta (n+r_0)(I_{n+r_0}(\omega^+)-\gamma)} \, d\mu^+\bigg) \\
& \le e^{C_1r_0\delta} \left(\sup_{\vec l} \int_{\Omega^+_{\vec l}} e^{\delta m(I_{m}(\omega^+)-\gamma)} \, d\mu^+_{\vec l}\right) \cdot \int e^{\delta (n+r_0)(I_{n+r_0}(\omega^+)-\gamma)} \, d\mu^+,
\end{align*}
where the last inequality follows from \eqref{eq:small_dist}. Now we choose $N$ large so that $(C')^{-1} N \kappa' > 2 C_1 r_0$ and we set $c = \frac12 (C')^{-1} \delta \kappa'$. By choosing $m = N$ in the above estimate and by \eqref{eq:small_exp1}, we have for all $n\ge 1$,
\beq\label{eq:ldt1}
\int e^{\delta (n+r_0+N)(I_{n+r_0+N}(\omega^+)-\gamma)} \, d\mu^+ \le e^{-cN} \int e^{\delta (n+r_0)(I_{n+r_0}(\omega^+)-\gamma)} \, d\mu^+.
\eeq

Now, given any $n > r_0$, we may apply the Euclidean division $n - r_0 = kN + r$. Using \eqref{eq:ldt1} several times, we obtain for all $n\ge 1$ that
$$
\mu^+ \left\{ \omega^+ : \frac1n S_n f(\omega^+) \ge \gamma \right\} \le \int_{\Omega^+} e^{\delta n(I_n(\omega^+)-\gamma)} \, d\mu^+ \le C e^{-cn}.
$$
Note that the estimate for small $n$ is absorbed into the constant $C$. Similarly, for any $\gamma < \int f \, d\mu^+$, we can apply the same argument to $\gamma-I_n(\omega)$  and obtain
$$
\mu^+ \left\{ \omega^+ : \frac1n S_n f(\omega^+) \le \gamma \right\} \le \int_{\Omega^+} e^{\delta n(\gamma-I_n(\omega^+))} \, d\mu^+ \le C e^{-cn}.
$$
The two inequalities above together clearly imply the desired large deviation estimate as stated in Theorem~\ref{t:ldt2}.
\end{proof}

The fact that an equilibrium state of a H\"older continuous potential has local product structure may be found at \cite{bowen2, BV}. Here we show that equilibrium states of H\"older continuous potentials have the bounded distortion property as defined in \eqref{eq:bdd}, which also implies that such a $\mu$ has a local product structure. In particular this shows that Theorem~\ref{t:ldt1} holds true for such measures. Equivalently, we may consider $(\Omega^+,\mu^+,T_+)$, where $\mu^+$ is an equilibrium state of a H\"older continuous potentials and show that such a $\mu^+$ has the bounded distortion property. Indeed, equilibrium states of H\"older continuous potentials defined over $(\Omega,T)$ are lifts of equilibrium states of H\"older continuous potentials defined over $(\Omega^+, T_+)$; see, for example, \cite{BV}. According to \cite{bowen2, BV}, such a $\mu^+$ has a H\"older continuous Jacobian with respect to $T_+$. So it suffices to prove the following lemma:

\begin{lemma}\label{l:es_to_bdd}
Let $(\Omega^+,T_+,\mu^+)$ be a one-sided subshift of finite type, where $\mu^+$ is a $T_+$-ergodic measure that has a H\"older continuous Jacobian. Then $\mu^+$ satisfies the bounded distortion property as defined in \eqref{eq:bdd+}.
\end{lemma}

\begin{proof}
	To get \eqref{eq:bdd+}, we fix any $[0;\vec l]^+\subset\Omega^+$ and set $n=|\vec l|$.  Choose any $[k;\vec j]^+\subset\Omega^+$ such that $k\ge n$ and $[0;\vec l]^+\cap [k;\vec j]^+\neq\varnothing$.
	
	Let $J_+\in C^\alpha(\Omega,\R_+)$ be the Jacobian of $\mu^+$ with respect to $T_+$. Since it is positive and continuous on $\Omega^+$, we have $\inf_{\omega^+\in\Omega^+} |J_+(\omega^+)| > c > 0$, which implies that $\log J_+\in C^\alpha(\Omega^+,\R_+)$. Consider the map $T^n:[0;\vec l]\to \Omega^+$ and let $J_+^{\vec l}$ be its Jacobian. Then we have
	$$
	J_+^{\vec l}(\omega^+)=\prod^{n-1}_{i=0}J_+(T^i\omega^+).
	$$
	Suppose $\omega^+, \tilde \omega^+\in [0;\vec l]$ for some $|\vec l|=n$. Then we have
	\begin{align*}
	\big|\log J_+^{\vec l}(\omega^+)-\log J_+^{\vec l}(\tilde \omega^+)\big|&=\sum^{n-1}_{i=0}\big|\log J_+(T^i\omega^+)-\log J_+(T^i\tilde \omega^+)\big|\\
	&<C\cdot d(T^i_+\omega^+,T^i_+\tilde\omega^+)^\alpha\\
	&<C,
	\end{align*}
	where $C$ is independent of $\vec l$ , $\omega^+$, and $\tilde\omega^+$. Thus we have
	\beq\label{eq:bdd_jacobian}
	C^{-1}<\bigg|\frac{J_+^{\vec l}(\omega^+)}{J_+^{\vec l}(\tilde\omega^+)}\bigg|<C \mbox{ for all } \omega^+, \tilde\omega^+\in[0;\vec l]^+.
	\eeq
	Now by definition of the Jacobian, we have
	\begin{align*}
	\int_{\Omega^+}\chi_{[k-n;\vec j]^+}(\eta) \, d\mu^+(\eta) = \int_{[0;\vec l]^+}\chi_{[k-n;\vec j]^+}(T^{n}_+\omega^+)J_+^{\vec l}(\omega^+) \, d\mu^+(\omega^+),
	\end{align*}
	which implies that
	$$
	J_+^{\vec l}(\omega^+_{\vec l, \min})\le \frac{\mu\big([k;\vec j]^+\big)}{\mu^+([0;\vec l]^+\cap[k;\vec j]^+)}
	\le J_+^{\vec l}(\omega^+_{\vec l, \max}).
	$$
	Here $\omega^+_{\vec l, \min}$ and $\omega^+_{\vec l, \max}$ are chosen as in the proof of Theorem~\ref{t:ldt2}. Using $1=\int_{\Omega^+}1d\mu^+=\int_{[0;\vec l]}J_+^{\vec l}(\omega^+)d\mu^+$, we obtain
	$$
	\frac1{J_+^{\vec l}(\omega^+_{\vec l, \max})}\le \mu^+([0;\vec l]^+)\le \frac1{J_+^{\vec l}(\omega^+_{\vec l, \min})}.
	$$
	Combining the two estimates above with \eqref{eq:bdd_jacobian}, we clearly get
	$$
	C^{-1}\le \frac{\mu^+\big([0;\vec l]^+\big)\cdot \mu^+\big([k;\vec j]^+\big)}{\mu^+\big([0;\vec l]^+\cap[k;\vec j]^+\big)}\le C
	$$
	which is \eqref{eq:bdd+}.
\end{proof}

\begin{remark}
There are many early works concerning large deviation estimates for functions defined on hyperbolic dynamical systems; see, for example, \cite{young2}. But we could not find a proof that applies in our framework. In Section~\ref{s:invPrinciple_Barycenter}, we shall show that if $\mu$ has the bounded distortion property, then Theorem~\ref{t:ldt1} yields relatively global versions of all the techniques needed. We hope they will be of independent interest. Most importantly, the results, ideas, and proofs given in this section will be used in the second paper \cite{ADZ} of this series.
\end{remark}

\section{Invariance Principle and Conformal Barycenter}\label{s:invPrinciple_Barycenter}

In this section, we develop some tools that are needed to prove our main results in the next two sections. Our main objective is to consider cocycles that have zero Lyapunov exponent. First, we show that a small Lyapunov exponent gives rise to a measurable family of holonomies, which will be integrable if $\mu$ has the bounded distortion property. In the case of a zero Lyapunov exponent, we shall introduce techniques originally developed in \cite{ledrappier}, and then generalized in \cite{BGV, AV, V}, which are referred to as the \emph{invariance principle}. We will use the invariance principle to show the existence of a continuous $su$-state on suitable sets. Concretely, in case of a zero Lyapunov exponent and bounded distortion, we will construct an $su$-state that is continuous on the support of a certain full measure set. In case we have only local product structure, we will construct a local $su$-state that is continuous on the support of some positive measure set. Then we show that periodic points with small Lyapunov exponents belong to the support in question. Finally, we will construct an $su$-invariant family of $\delta$-measures by using the conformal barycenter.

\subsection{Holonomies Resulting from Small Exponents}

For the remainder of this subsection, we fix $0 < \alpha \le 1$ and consider the space of $\alpha$-H\"older continuous cocycles $C^\alpha(\Omega,\SL(2,\R))$, that is, $A \in C^\alpha(\Omega,\SL(2,\R))$ if
$$
\| A(\omega) - A(\tilde \omega) \| < C \cdot d(\omega,\tilde \omega)^\alpha \mbox{ for all }\omega,\omega'\in\Omega.
$$
Thus for every $\omega,\tilde \omega \in \Omega$ and $n \geq 0$, we have
$$
\begin{cases}
\|A(T^n \omega)-A(T^n \tilde \omega)\| \leq C e^{-\alpha n}  & \mbox{ if } \tilde\omega\in W^s_\loc(\omega),\\
\|A(T^{-n} \omega) - A(T^{-n} \tilde \omega)\| \leq C e^{-\alpha n}  & \mbox{ if } \tilde \omega\in W^u_\loc(\omega).
\end{cases}
$$

\begin{lemma}\label{p.bgvprop3}
	Assume $2 L(A,\mu) < \alpha$. For each $\delta > 0$ with $2L(A,\mu) < \alpha - \delta$, there exists an increasing family of compact sets $K_s(N,\delta) \subset \Omega$, $N > 0$, whose union over $N$ has full measure and is $T$-invariant, such that there is a measurable family of stable holonomies defined on $\bigcup_N K_s (N,\delta)$, which are continuous when restricted to any $K_s(N,\delta)$. Moreover, if $\mu$ has bounded distortion, then these holonomies obey the following integrability condition:
	\begin{equation}\label{e.integrability}
		\int_{\Omega} \log \|H^s_{\omega^{(\omega_0)} \wedge \omega, \omega}\| \, d\mu(\omega) < \infty.
	\end{equation}
	Similarly, there are sets $K_u(N,\delta)$, $N > 0$ and unstable holonomies with the same properties. In particular, if $\mu$ has bounded distortion property, then the unstable holonomies satisfies
	\begin{equation}\label{e.integrability1}
		\int_{\Omega} \log \|H^u_{\omega\wedge \omega^{(\omega_0)}, \omega}\| \, d\mu(\omega) < \infty.
	\end{equation}
\end{lemma}

\begin{proof}
	Let $H^{s,n}_{\omega,\tilde \omega} = A_n(\tilde \omega)^{-1} A_n(\omega)$ for $\tilde \omega \in W^s_{\mathrm{loc}}(\omega)$. Define
	$$
	\delta^{s,n}_{\omega,\tilde \omega} = \left( H^{s,n}_{\omega,\tilde \omega} \right)^{-1} \left( H^{s,n+1}_{\omega,\tilde \omega} - H^{s,n}_{\omega,\tilde \omega} \right),
	$$
	so that
	$$
	H^{s,n}_{\omega,\tilde \omega} \left( \mathrm{Id} + \delta^{s,n}_{\omega,\tilde \omega} \right) = H^{s,n+1}_{\omega,\tilde \omega}.
	$$
	We first will estimate $\|\delta^{s,n}_{\omega,\tilde \omega}\|$ as follows. In the case where $\|A_n(\omega)\|^2 \le e^{(\alpha - \delta)n}$, we have
	$$
	\delta^{s,n}_{\omega,\tilde \omega} = A_n(\omega)^{-1} \left( A(T^n \tilde \omega)^{-1} A(T^n \omega) - \mathrm{Id} \right) A_n(\omega)
	$$
	and therefore
	\beq\label{eq:bdd_difference}
	\|\delta^{s,n}_{\omega,\tilde \omega}\| \le C e^{(\alpha - \delta)n} e^{-\alpha n} = C e^{-\delta n}.
	\eeq
	The fact $\lim\limits_{n\to\infty}\frac1n\log\|A_n(\omega)\|=L(A,\mu)$ for $\mu$-a.e. $\omega\in\Omega$ implies that $\frac1n\log\|A_n(\omega)\|$ converges to $L(A,\mu)$ in measure. Thus, if we define for $N > 0$,
	\beq\label{eq:bunching1}
	K_s(N,\delta) = \{ \omega : \|A_n(\omega)\|^2 \le e^{(\alpha - \delta)n} \text{ for every } n \ge N \},
	\eeq
	then the sets $K_s(N,\delta)$ are compact and increasing in $N$, and their union over $N$ has full measure. For $\omega \in K_s(N,\delta)$, we have the uniform summability statement
	$$
	\sum_{n = N}^\infty \|\delta^{s,n}_{\omega,\tilde \omega}\| \le C,
	$$
	where $C$ depends on $\delta$ and $\|A\|_{\infty}$ so that $\|A_n(\omega)\|<e^{Cn}$ for all $\omega$ and all $n$. It then implies that
	\begin{align}\label{eq:bounded_distortion}
	\nonumber \left\|\lim_{n \to \infty} H^{s,n}_{\omega,\tilde \omega}\right\|&=\left\|H^{s,N}_{\omega,\tilde \omega}\prod^{\infty}_{n=N}(H^{s,n}_{\omega,\tilde \omega})^{-1}\cdot H^{s,n+1}_{\omega,\tilde \omega}\right\|\\
	\nonumber
	&\le \|H^{s,N}_{\omega,\tilde \omega}\|\prod^{\infty}_{n=N}\|(H^{s,n}_{\omega,\tilde \omega})^{-1}\cdot H^{s,n+1}_{\omega,\tilde \omega}\|\\
	&\le e^{CN}\exp\left(\sum^{\infty}_{n=N}\log(1+\|\delta^{s,n}_{\omega,\tilde\omega}\|)\right)\\
	\nonumber &\le e^{CN}\exp\left(C\sum^{\infty}_{n=N}\|\delta^{s,n}_{\omega,\tilde\omega}\|\right)\\
	\nonumber &\le e^{CN},
	\end{align}
	where again $C=C(\delta,A)$. Hence we may define $H^s_{\omega,\tilde \omega} := \lim\limits_{n \to \infty} H^{s,n}_{\omega,\tilde \omega}$, where the convergence is uniform on $K_s(N,\delta)$. In particular, $H^s_{\omega,\tilde \omega}$ depends continuously on $\omega \in K_s(N,\delta)$ and $\tilde \omega \in W^s_{\mathrm{loc}}(\omega)$. Note that \eqref{eq:bounded_distortion} implies for the same $\omega$ and $\tilde \omega$ that
	\beq\label{eq:bounded_distortion1}
	\|A_n(\tilde\omega)^{-1}A_n(\omega)\|<e^{CN}
	\eeq
	for all $n\ge 1$.
	
	To get the integrability condition \eqref{e.integrability}, we define $\phi(\omega)=\log\|A(\omega)\|$ and assume without loss of generality that
	\beq\label{eq:small_initial_average}
	2\int_{\Omega}\phi(\omega) \, d\mu < \alpha - \delta,
	\eeq
	because otherwise we may instead consider $\phi(\omega) = \frac1k \log\|A_k(\omega)\|$ for some large $k$, which must satisfy the condition above since $\int_\Omega \frac1k \log\|A_k(\omega)\| \, d\mu$ converges to $L(A,\mu)$. It is straightforward to see that $\phi \in C^\alpha(\Omega,\R)$ since $A$ is $\alpha$-H\"older continuous and $\|A(\omega)\|\ge 1$ for all $\omega\in\Omega$.

	We want to estimate for some $C>0$ the measure of the following set,
	$$
	\CB_N=\{ \omega : \log \|H^s_{\omega,\tilde \omega}\| > CN \text{ for some } \tilde \omega \in W^s_{\mathrm{loc}}(\omega) \}.
	$$
	To this end, we define for $\delta'=\frac\delta2$,
	$$
	\CZ_m=\bigg\{\omega: \frac1n\sum^{n-1}_{j=0}\phi(T^j\omega)<\frac{\alpha-\delta'}2\mbox{ for all }n\ge m\bigg\}
	$$
	If $\omega\in\CZ_1$, then $2\phi(T^n\omega)<\alpha-\delta$ for all $n\ge 0$, which clearly implies that $\|A_n(\omega)\|^2<e^{(\alpha-\delta)n}$ for all $n\ge 0$. Thus by the computation leading to \eqref{eq:bounded_distortion1}, we have for all $\omega\in\CZ_1$,
	$$
	\|A_n(\tilde \omega)^{-1}A_n(\omega)\|\le C \mbox{ for all }\tilde\omega\in W^s_{\mathrm{loc}}(\omega) \mbox{ and for all } n\ge 1.
	$$
	If $\omega\in\CZ_m$ with $m>1$, then we set $0<k<m$ to be the largest integer for which $\omega \notin \CZ_{k}$. A direct computation then shows that $T^k \omega \in\CZ_1$, which implies that for all $\tilde\omega\in W^s_{\mathrm{loc}}(T^k\omega)$,
	$$
	\|A_n(\tilde \omega)^{-1}A_n(T^k\omega)\|\le C \mbox{ for all }n\ge 1.
	$$
	Combining $\|A_n(\omega)\|<e^{Cm}$ for all $1\le n\le m$ and for all $\omega$, we have for all $\omega\in\CZ_m$ and  all $\tilde\omega\in W^s_{\mathrm{loc}}(\omega)$ that
	$$
	\|A_n(\tilde \omega)^{-1}A_n(\omega)\|\le e^{Cm} \mbox{ for all } n\ge 1,
	$$
	which clearly implies that $\log\|H^s_{\omega,\tilde\omega}\|<Cm$. Thus by choosing $C$ appropriately, we have
	$$
	\CB_N\subset\Omega\setminus\CZ_N.
	$$
	However, by \eqref{eq:small_initial_average} it is clear that
	$$
	\Omega \setminus \CZ_N \subseteq \bigcup^{\infty}_{n=N} \bigg\{ \omega : \bigg| \frac1n \sum^{n-1}_{j=0} \phi(T^j\omega) - \int_\Omega \phi \, d\mu \bigg| > \frac\delta4 \bigg\}.
	$$
	Suppose that $\mu$ has bounded distortion. Note that $\phi\in C^\alpha(\Omega,\R)$. Hence, by Theorem~\ref{t:ldt1}, there exist $C>0$ and $\eta>0$ such that
	$$
\mu \bigg\{ \omega : \bigg| \frac1n \sum^{n-1}_{j=0} \phi(T^j\omega) - \int_\Omega \phi \, d\mu \bigg| > \frac\delta4 \bigg\} < C e^{-\eta n} \mbox{ for all } n\ge 1.
	$$
	Clearly, this implies that $\mu(\CB_N) < C e^{-\eta N}$ for all $N \ge 1$, which in turn implies the integrability statement \eqref{e.integrability} since
	\begin{align*}
	\int_\Omega \log \| H^s_{\omega^{(\omega_0)}\wedge \omega,\omega} \| \,  d\mu & = \int_\Omega \log \|H^s_{\omega,\omega^{(\omega_0)}\wedge \omega}\| \, d\mu\\
	& = \sum^{\infty}_{N=1} \int_{\CB_{N} \setminus \CB_{N+1}} \log \|H^s_{\omega,\omega^{(\omega_0)}\wedge \omega}\| \, d\mu \\
	& \le \sum^{\infty}_{N=1} \mu(\CB_N) C N \\
	& \le \sum^{\infty}_{N=1} C e^{-\eta N}N \\
& <\infty.
	\end{align*}	
The case of $K_u(N,\delta)$ can be done similarly after replacing $A_n(\omega)$ by $A_{-n}(\omega)$.
\end{proof}

\begin{definition}
For a periodic point $p$ with period $n$, we let $L(A,p)=\lim\frac1n\log\|A_n(p)\|$ be the individual Lyapunov exponent of $A$ at $p$. We say $p$ is $\gamma$-\emph{bunched} if $2L(A,p) < \gamma \le \alpha$.
\end{definition}

Next, we show the following result, which says that $\alpha$-bunched periodic points are in the support of $K_s(N,\delta)$ for suitable $\delta > 0$ and for large $N$. We note that a similar result has appeared in \cite{butler}.\footnote{We would like to thank Clark Butler for pointing this out to us.}

\begin{lemma}\label{l:PinSuppK}
Suppose $(\Omega, T)$ is a subshift of finite type and $\mu$ has a local product structure. Assume that $2L(A,\mu)<\alpha$. Let $p$ be an $\alpha$-bunched periodic point. Set $\delta = \min \{ \alpha - 2L(A,\mu), \alpha-2L(A,p) \}$. Then for every $0 < \delta_0 < \delta$, there exists $N_0 \in \Z_+$ such that
	$$
	p \in \mathrm{supp}\left(\mu|_{K_{s}(N_0,\delta_0)\cap K_{u}(N_0,\delta_0)}\right).
	$$
\end{lemma}

\begin{proof}
	Assume the period of $p$ is $r$. By definition, we have that $2L(A,p) < \alpha - \delta$ and $2L(A,\mu)<\alpha-\delta$.
	
	Consider the family $K_s(N,\delta)$ as in Proposition \ref{p.bgvprop3}. By \eqref{eq:bounded_distortion1}, we have for all $\omega\in K_s(N,\delta)$, all $\tilde\omega\in W^s_{\mathrm{loc}}(\omega)$, and all $n\ge 1$,
	$$
	\|A_n(\tilde\omega)^{-1}A_n(\omega)\|<CN.
	$$
	Recall \eqref{eq:local_su_set} says $(\pi^+)^{-1}(\pi^+\omega)=W^s_\mathrm{loc}(\omega)$. Thus by definition of $K_s(N,\delta)$ and the estimate above, for each $0<\delta_1<\delta$, there exists $N_1>N$ such that  for all $\omega\in(\pi^+)^{-1}(\pi^+[K_s(N,\delta)])$ and all $n\ge N_1$, we have
	\beq\label{eq:norm_bounded_above1}
	\|A_n(\omega)\|^2<e^{(\alpha-\delta_1)n}.
	\eeq
	Fix such a choice of $\delta_1$ and $N_1$. By choosing $N$ large, we may assume that $\mu( K_s(N,\delta)\cap[0;i])>0$ for each $1\le i\le \ell$, which in turn implies that
	$$
	\mu^+[\pi^+(K_s(N,\delta))\cap[0;i]^+)=\mu( K_s(N,\delta)\cap[0;i])>0.
	$$
	By Corollary~\ref{c:jacobian}, for each $n\ge 0$, we have
	\beq\label{eq:positive_measure1}
	\mu^+(T^{-n}_+[\pi^+(K_s(N,\delta))]\cap[0;p_0,\ldots, p_{n}]^+) > 0.
	\eeq
	
	By the same argument as the one leading to \eqref{eq:bounded_distortion1}, $2L(A,p) < \alpha - \delta$ implies that for each $0 < \delta_2 < \delta$, we can find $m \in \Z_+$ large enough so that $\|A_{rm}(\omega)\|^2 \le e^{(\alpha - \delta_2)rm}$ for $\omega\in T^{-rm}W^u_{\mathrm{loc}}(T^{rm} p)$. By periodicity of $p$, we have for each $l \ge 1$ and each $1\le k < l$,
	$$
	\|A_{rm}(T^{krm}\omega)\|^2 \le e^{(\alpha - \delta_2)rm} \mbox{ for all } \omega \mbox{ with } T^{-lrm} \omega \in W^u_{\mathrm{loc}}(T^{lrm} p),
	$$
	which in turn implies that for each $1\le k\le l$,
	\beq\label{eq:norm_bounded_above2}
	\|A_{krm}(\omega)\|^2<\prod^{k-1}_{j=0}\|A_{rm}(T^{jrm})(\omega)\|^2\le e^{(\alpha-\delta_2)krm}.
	\eeq
	
	For each $l\in\Z_+$, we define the following $s$-\textit{locally saturated} set,
	$$
	\DD^l_+=(\pi^+)^{-1}(T^{-lrm}_+[\pi^+(K_s(N,\delta))]\cap[0;p_0,\ldots, p_{lrm}]^+).
	$$
	By \eqref{eq:positive_measure1}, we have $\mu(\DD^l_+)>0$ and $\DD^l_+\subset [0;p_0,\ldots, p_{lrm}]$. For each $0<\delta_3<\delta_2$,  we can fix a $N'\in\Z_+$ large enough so that the following holds true. For all $l$ large and for each $\omega\in\DD^l_+$, we have for all $N'\le n\le lrm+N_1$ that
	$$
	\|A_n(\omega)\|^2<\|A_{krm}(\omega)\|^2\cdot \|A_{n-krm}(T^{krm}\omega)\|^2<e^{(\alpha-\delta_3)n},
	$$
	where $k$ is so chosen that $0\le n-N_1-krm<rm$. On the other hand, if $n>lrm+N_1$, then we have
	\begin{align*}
	\|A_n(\omega)\|^2
	&\le \|A_{lrm}(\omega)\|^2\cdot \|A_{n-lrm}(T^{krm}\omega)\|^2\\
	&\le \|A_{lrm}(\omega\wedge p)^{-1}\cdot A_{lrm}(\omega)\|^2\cdot \|A_{lrm}(\omega\wedge p)\|^2\cdot \|A_{n-lrm}(T^{krm}\omega)\|^2.
	\end{align*}
	We estimate each factor in the product above. First we consider the last factor. The fact that $\omega\in\DD^l_u$ implies that $T^{krm}(\omega)\in W^s_{\mathrm{loc}}(\omega')$ for some $\omega'\in K_s(N,\delta)$. Thus by \eqref{eq:norm_bounded_above1} and the fact that $n-lrm>N_1$, we have
	$$
	\|A_{n-lrm}(T^{krm}\omega)\|^2<e^{(\alpha-\delta_1)(n-lrm)}.
	$$
	For the second factor, the fact $\omega\in\DD^l_u$ implies that $T^{-lrm}(\omega\wedge p)\in W^u_{\mathrm{loc}}(T^{lrm}p)$. Thus by \eqref{eq:norm_bounded_above2}, we have for each $1\le k\le l$
	$$
	\|A_{krm}(\omega\wedge p)\|^2\le e^{(\alpha-\delta_2)krm}.
	$$
	In particular, the second factor is taken care of by choosing $k=l$. Combining the fact $\omega\in W^s_{\mathrm{loc}}(\omega\wedge p)$ with the estimate above and using
	$$
	A_{lrm}(\omega)=\prod^{0}_{j=l-1}A_{rm}((T^{rm})^j\omega),
	$$
	the same argument getting \eqref{eq:bounded_distortion1} yields
	$$
	\|A_{lrm}(\omega\wedge p)^{-1}\cdot A_{lrm}(\omega)\|^2<C.
	$$
	Thus by setting $\delta'=\min\{\delta_1, \delta_3\}$, we have for all large $l$ and $n\ge lrm+N_1$ that
	$$
	\|A_{n}(\omega)\|^2<e^{(\alpha-\delta')n}.
	$$
	Combining the estimates in the case of $N'\le n\le lrm+N_1$, we obtain for all $l$ large, all $\omega\in\DD^l_u$, and all $n\ge N'$ that
	$$
	\|A_n(\omega)\|^2<e^{(\alpha-\delta')n},
	$$
	which implies that $\DD^l_+\subset K_s(N',\delta')\cap [0;p_0,\ldots p_{lrm}]$ for all large $l$. Note that $0<\delta_i<\delta$, $i=1,2,3$, are arbitrarily chosen, hence $\delta'$ can be any number in $(0,\delta)$. In particular, we have for all $\delta_0\in (0,\delta)$, there is a $N'$ such that $\DD^l_+\subset K_s(N',\delta_0)\cap [0;p_0,\ldots p_{lrm}]$ for all large $l$.
	
	Similarly, for each $0<\delta_0<\delta$, we can find a $N''\in\Z_+$  and a sequence of $u$-\textit{locally saturated}  $\DD^l_-\subset K_u(N'',\delta_0)\cap [-lrm; p_{-lrm},\ldots, p_{-1}, p_0]$ with $\mu(\DD^l_-)>0$. Taking $N_0=\max\{N',N''\}$, we have for all large $l$,
	\begin{align*}
	&\DD^l_-\cap \DD^l_+\subset K_s(N_0,\delta_0)\cap K_u(N_0,\delta_0),\\
	&\DD^l_-\cap \DD^l_+\subset [-lrm; p_{-lrm},\ldots, p_{lrm}],
	\end{align*}
	where the second line implies that $\DD^l_-\cap \DD^l_+$ is contained in arbitrarily small neighborhood of $p$ as $l$ tends to infinity. Finally, combining that $\DD^l_+$ is $s$-\textit{locally saturated} in $[0;p_0]$, $\DD^l_-$ is $u$-\textit{locally saturated} in $[0;p_0]$, and \eqref{eq:local_product_structure}, we have for all $l$ large,
	$$
	\mu(\DD^l_-\cap\DD^l_+)>0,
	$$
	which then implies that $p \in \mathrm{supp}\left(\mu|_{K_{s}(N_0,\delta_0)\cap K_{u}(N_0,\delta_0)}\right).$
\end{proof}

\subsection{Invariance Principle and $su$-States}

Let $\CM$ be the Borel $\sigma$-algebra of the subshift of finite type $(\Omega,T,\mu)$, where $\mu$ has a local product structure. Let $A:\Omega\to\SL(2,\R)$ be a measurable map. Then the following \emph{invariance principle} is due to Ledrappier \cite{ledrappier}, see also \cite{AV, V}:

\begin{prop}\label{p:invPrinciple} Let $\CB \subseteq \CM$ be a $\sigma$-algebra such that
	\begin{enumerate}
		\item $T^{-1}\CB \subseteq \CB $ mod $0$ and $\{T^n\CB: n\in\Z\}$ generates $\CM$ mod $0$.
		\item the $\sigma$-algebra generated by $A$ is contained in $\CB$ mod 0.
		\end{enumerate}
	If $L(A,\mu)=0$, then for any $(T,A)$-invariant measure $m$ on $\Omega\times \R\bbP^1$ that projects to $\mu$ in the first component, the disintegration $\{m_\omega\}_{\omega\in\Omega}$ is $\CB$-measurable mod $0$.
	\end{prop}

\begin{definition}
We say that a function defined on $\Omega$ \textit{only depends on the future} (resp., \emph{past}) if it is constant on every local stable (resp., unstable) set.
\end{definition}

The following consequence of Proposition~\ref{p:invPrinciple} is due to \cite{BGV}. We sketch a proof for the reader's convenience.

\begin{prop}\label{p.bgvprop1}
	Suppose $A$ only depends on the future and $L(A,\mu) = 0$. Then for every $(T,A)$-invariant measure $m$ on $\Omega \times \R\bbP^1$ that projects to $\mu$ in the first component, its disintegration only depends on the future for $\mu$-almost every $\omega\in\Omega$.
\end{prop}

\begin{proof}
	Let $\CB \subseteq \CM$ be the $\sigma$-algebra generated by sets $\{W^s_\loc(\omega): \omega\in \Omega\}$. It is clear that the sets $W^s_\loc(\omega)$ are mutually disjoint. Thus, $D \in \CB$ if and only if for each $\omega \in\Omega$, either $W^s_\loc(\omega)\cap D = \varnothing$ or $W^s_\loc (\omega) \subseteq D$. Since $T \CB$ is the $\sigma$-algebra generated by $\{ T W^s_\loc(\omega) : \omega \in \Omega\}$, it is clear that $\CB \subseteq T\CB$, or equivalently $T^{-1} \CB \subseteq \CB$. More generally, $T^n \CB$ is generated by $\{ T^n W^s_\loc(\omega) : \omega \in \Omega\}$. Now for any cylinder $[n;\vec l]\subset \Omega$, it is clear that it is $T^n \CB$-measurable for some large $n \in \Z_+$. Since $\CM$ is generated by cylinders, we then have that $\{ T^n \CB : n \in \Z \}$ generates $\CM$ mod $0$. The result then follows from Proposition~\ref{p:invPrinciple} and the straightforward fact that $A$ is $\CB$-measurable if and only if $A$ depends on the future.
\end{proof}

An immediate consequence of Proposition~\ref{p.bgvprop1} is that if $A$ is constant along the local stable set and $L(A,\mu)=0$, then for every $(T,A)$-invariant measure $m$ on $\Omega \times \R\bbP^1$ that projects to $\mu$ in the first component, its disintegration is constant on the local stable set $W^s_\loc(\omega)$ for $\mu$-almost every $\omega \in \Omega$. Indeed, we just need to define $\omega' = \varphi(\omega)$ to be the sequence for which $\omega'_n=\omega_{-n}$ for all $n\in\Z$ and set
$$
\Omega' := \{ \omega' = \varphi(\omega) : \omega \in \Omega\}.
$$
Then $\mu$ is again an ergodic measure of $(\Omega',T)$ which has a local product structure. Set $A'(\omega') = A(f(\omega))$ so that $A'$ depends only on the past. Then it is a standard result that $L(A',\mu) = L(A,\mu)=0$ and $m$ is $(T,A')$-invariant if it is $(T,A)$-invariant. Now the conclusion follows from Proposition~\ref{p.bgvprop1}.

We have the following version of $su$-states. It is weaker than the original version where one has canonical holonomies because here the holonomies only exist almost everywhere.

\begin{prop}\label{p.bgvprop2}
	Suppose the cocycle map $A$ is measurable, satisfies  the integrability condition $\int_\Omega \log \|A(\omega)\| \, d\mu < \infty$, and admits stable and unstable holonomies almost everywhere. Suppose that the holonomies satisfy the integrability conditions \eqref{e.integrability} and \eqref{e.integrability1}.
	If $L(A,\mu) = 0$, then every $(T,A)$-invariant measure $m$ on $\Omega \times \R\bbP^1$ that projects to $\mu$ in the first component has a disintegration that is almost surely invariant under the stable and unstable holonomies.
\end{prop}

\begin{proof}
	First we consider the $s$-invariance. For simplicity, we define $\varphi(\omega)=\omega^{(\omega_0)} \wedge \omega$, which depends only on the future. We define a new cocycle map as follows:
	$$
	\tilde A(\omega):=H^s_{T\omega, \varphi(T\omega)}\cdot A(\omega)\cdot H^s_{\varphi(\omega), \omega}.
	$$
	It is clear that $\tilde A$ is conjugate to $A$ via the stable holonomy. By the condition \eqref{e.integrability} and the definition of $\tilde A$, we then obtain $\int_\Omega\log\|\tilde A(\omega)\| \, d\mu < \infty$ and $L(\tilde A,\mu)=0$. On the other hand, by conditions (i)--(ii) of the definition of stable holonomy, we have that
	\begin{align*}
	\tilde A(\omega)
	&=H^s_{T\omega, \varphi(T\omega)}\cdot A(\omega)\cdot H^s_{\varphi(\omega),\omega}\\
	&= H^s_{T\omega, \varphi(T\omega)}\cdot H^s_{T\varphi(\omega),T\omega}\cdot A(\varphi(\omega))\\
	&= H^s_{T\varphi(\omega),\varphi(T\omega)}\cdot A(\varphi(\omega)),
	\end{align*}
	which implies that $\tilde A(\omega)$ depends only on the future. Thus Proposition~\ref{p.bgvprop1} implies that we have for every $(T,\tilde A)$-invariant measure $m$ that projects to $\mu$ in the first component, its disintegration only depends on the future.
	
	Now let $m$ be a $(T,A)$-invariant measure that projects to $\mu$ in the first component. Let $\{m_\omega: \omega\in\Omega\}$ be a disintegration of $m$. Thus $A(\omega)_*m_\omega=m_{T\omega}$ for $\mu$-almost every $\omega$. We define
	$$
	\tilde m_\omega=(H^s_{\omega,\varphi(\omega)})_*m_\omega,\ \omega\in\Omega.
	$$
	One readily checks that $\tilde A(\omega)_*\tilde m_\omega=\tilde m_{T(\omega)}$. Thus the family of conditional measures $\{\tilde m_\omega:\omega\in\Omega\}$ is a disintegration of a $(T,\tilde A)$-invariant measure $\tilde m$. Thus $\tilde m_\omega$ depends only on the future. In other words, for $\mu$-almost every $\omega$, we have for each $\omega'\in W^s_\loc(\omega)$ that
	$$
	(H^s_{\omega,\varphi(\omega)})_*m_\omega=(H^s_{\omega',\varphi(\omega')})_*m_{\omega'}.
	$$
	Since $\varphi(\omega)=\varphi(\omega')$, by condition (i) of the definition of stable holonomy we have
	\beq\label{eq:s-inv}
	m_\omega=(H^s_{\varphi(\omega'),\omega}\cdot H^s_{\omega',\varphi(\omega')})_*m_{\omega'}=(H^s_{\omega',\omega})_*m_{\omega'}.
	\eeq
	In other words, $\{m_\omega:\omega\in\Omega\}$ is $s$-invariant $\mu$-almost everywhere, concluding the proof of the $s$-invariance.
	
	As for the $u$-invariance, we just need to conjugate $A$ to a new $\tilde A$ via the unstable holonomy so that  $\tilde A$ is constant along the local unstable set $W^u_\loc(\omega)$ for $\mu$-almost every $\omega\in\Omega$. Then by repeating the same argument above and using the remark following Proposition~\ref{p.bgvprop1}, we obtain that ${m_\omega}$ is $u$-invariant $\mu$-almost everywhere. This completes the proof.
\end{proof}

\begin{lemma}\label{l.lyap0zp}
	Assume that $L(A,\mu)=0$ and $\mu$ has the bounded distortion property. Then there exists a full measure set $K \subset \Omega$ on which one has stable and unstable holonomies. Moreover, every $(T,A)$-invariant measure $m$ on $\Omega \times \R\bbP^1$ that projects to $\mu$ in the first component has a continuous, $su$-invariant disintegration over $\mathrm{supp}(\mu|_K)\cap K$.
\end{lemma}

\begin{proof}
	Since $L(A,\mu)=0$, Lemma~\ref{p.bgvprop3} applies and yields for each $\delta$ with $0<\delta<\alpha$, the sets $K_s(N,\delta)$, $K_u(N,\delta)$ along with continuous families of holonomies satisfying the estimates required to apply Proposition~\ref{p.bgvprop2}. Thus, applying Proposition~\ref{p.bgvprop2}, choose a $(T,A)$-invariant measure $m$ on $\Omega \times \R\bbP^1$ and consider its disintegration $\{ m_\omega \}$, which is invariant almost everywhere with respect to the stable and unstable holonomies. Recall that both $\bigcup_{N > 0} K_s(N,\delta)$ and $\bigcup_{N > 0} K_u(N,\delta)$ have full measure. We let
	$$
	K_\delta = \left( \bigcup_{N > 0} K_s(N,\delta) \right) \cap \left( \bigcup_{N > 0} K_u(N,\delta) \right).
	$$
	
	As in \cite{BGV}, we can now produce a disintegration $\{ \tilde m_\omega \}$ over $\mathrm{supp} (\mu|_{K_\delta})\cap K_\delta$, which is holonomy-invariant and continuous. For the reader's convenience we give the argument. In the following argument, we will work in the full measure set $K_{\delta}$ so that the stable and unstable holonomies are defined on the local stable and unstable set of $\omega\in K_\delta$, respectively.
	
	For each $\omega\in K_\delta$, if $m_{\omega'}$ already exists for some $\omega'\in W^s_\loc(\omega)$ from the original disintegration of $m$ above, then we may define $m^s_{\omega''}$ via $H^s(\omega', \omega'')_*m_{\omega'}$ for each $\omega''\in W^s_\loc(\omega)$. If $m_{\omega'}$ does not exist from the original disintegration of $m$ for all $\omega'\in W^s_\loc(\omega)$, then $W^s_\loc(\omega)$ is a $\mu$-zero measure set and we may pick any probability measure $m^s_\omega$ and extend via $m^s_{\omega''}=H^s(\omega, \omega'')_*m^s_{\omega}$ for each $\omega''\in W^s_\loc(\omega)$. Clearly, the new family $m^s_\omega$ is invariant under the stable holonomy at every $\omega\in K_{\delta}$. On the other hand, due to the almost sure invariance of the original disintegration $m_\omega$ of $m$, the new family $m^s_\omega$ coincides with $m_\omega$ for $\mu$-almost every $\omega$ in $K_\delta$, and hence it also coincides with $m_\omega$ for $\mu$-almost every $\omega$. In particular, $m^s_\omega$ is again a disintegration of $m$. Similarly, we may construct another disintegration $m^u_\omega$ of $m$ which is invariant under the unstable holonomy at every $\omega\in K_\delta$. Note that the set $\widetilde K=\{\omega\in K_\delta : m^s_{\omega} = m^u_{\omega} \}$ has full $\mu$-measure.
	
	Clearly, for each $1\le j\le \ell$, $[0;j]\cap\widetilde K$ has full $\mu$-measure in $[0;j]$. By the local product structure of $\mu$, $[0;j]\cap\widetilde K$ has full $\mu^-\times\mu^+-$measure in $[0;j]$. Thus by Fubini's theorem, for $\mu^-$-almost every $\omega^-\in [0;j]^-$, we have that $\pi^+[(\{\omega^-\}\times [0;j]^+)\cap \widetilde K]$ has full $\mu^+$-measure in $[0;j]^+$. Note for each $\omega\in[0;j]$ with $\pi^-(\omega)=\omega^-$, we have $\omega^-\times[0;j]^+=W^u_\loc(\omega)$. Thus for each $1\le j\le \ell$, we may choose an $\omega^{(j)}\in [0;j] \cap K_\delta$ such that
	$$
	\mu^+\left(\pi^+\big(W^u_{\mathrm{loc}}(\omega^{(j)})\cap \widetilde K\big)\right)=\mu^+([0;j]^+).
	$$
	By the definition of $\mu^+$, we then have that
	\begin{align}\label{eq:fullmeasure}
	\nonumber \mu\left((\pi^+)^{-1}\bigg[\pi^+\big(\bigcup^{\ell}_{j=1}\big(W^u_{\mathrm{loc}}(\omega^{(j)})\cap \widetilde K\big)\big)\bigg]\right)	& = \mu^+ \left( \pi^+ \bigg( \bigcup^{\ell}_{j=1} \big(W^u_{\mathrm{loc}}(\omega^{(j)}) \cap \widetilde K\big)\bigg) \right)\\
	 & = \sum^{\ell}_{j=1}\mu^+([0;j])^+ \\
\nonumber & = 1.
	\end{align}
	In other words, for $\mu$-almost every $\omega\in K_\delta$, we have $\omega^{(\omega_0)}\wedge\omega\in W^u_\loc(\omega^{(\omega_0)})\cap \widetilde K$. Now for each $\omega\in K_\delta$, we define
	$$
	\tilde m_\omega^s = H^s_{\omega^{(\omega_0)}\wedge \omega,\omega} \cdot m^u_{\omega^{(\omega_0)}\wedge \omega}= H^s_{\omega^{(\omega_0)}\wedge \omega,\omega} \cdot H^u_{\omega^{(\omega_0)},\omega^{(\omega_0)}\wedge \omega}\cdot m^u_{\omega^{(\omega_0)}}.
	$$
	Recall that by the proof of Lemma~\ref{p.bgvprop1}, the stable and unstable holonomies are continuous on each local stable and unstable set, respectively. Thus the equalities in the definition of $\tilde m_\omega$ above imply that for each $1\le j\le \ell$, we have that $\tilde m^s_\omega$ is continuous in $\omega$ at $[0;j]\cap K_\delta$. Thus $\tilde m^s_\omega$ is continuous on $K_\delta$. Clearly, the construction also implies that $\tilde m^s_\omega$ is $s$-invariant. On the other hand, by invariance with respect to stable holonomies, we have that for each $\omega$ such that $\omega^{(\omega_0)}\wedge \omega\in W^u_{\mathrm{loc}}(\omega^{(\omega_0)})\cap \widetilde K_{\omega_0}$, we have
	$$
	\tilde m^s_\omega = H^s_{\omega^{(\omega_0)}\wedge \omega,\omega} \cdot m^u_{\omega^{(\omega_0)}\wedge \omega} = H^s_{\omega^{(\omega_0)}\wedge \omega,\omega} \cdot m^s_{\omega^{(\omega_0)}\wedge \omega} = m^s_\omega.
	$$
	Thus, we have $\tilde m_\omega^s = m_\omega^s$ for $\mu$-almost every $\omega \in K_\delta$, and we obtain an $s$-invariant and continuous disintegration $\{ \tilde m^s_\omega \}$ of $m$. Producing in an analogous fashion a $u$-invariant and continuous disintegration $\{ \tilde m^u_\omega \}$, we find that $\tilde m^s_\omega = \tilde m^u_\omega$ in $\mathrm{supp} (\mu|_{K_\delta})\cap K_\delta$ by continuity and almost everywhere coincidence. This produces an $su$-invariant continuous disintegration $\{ \tilde m_\omega \}$ over $\mathrm{supp}(\mu|_{K_\delta})\cap K_\delta$ by setting $\tilde m_\omega=\tilde m^s_\omega$. By continuity, we also have invariance under $(T,A)$, that is, $A(\omega)_*\tilde m_\omega = \tilde m_{T\omega}$ for every $\omega \in \mathrm{supp}(\mu|_{K_\delta})\cap K_\delta$. Clearly, any $K_\delta$ can be chosen to be the desired $K$.
\end{proof}

\subsection{Application of Conformal Barycenter}

Let $\mathbb H \subseteq \C$ be the upper-half plane, $\D$ the open unit disk, and $S^1=\partial \D$ the unit circle. It is a standard result that the M\"obius transformation associated with an element of the group $\mathrm{SU}(1,1)$ preserves $S^1$ and $\D$. Here $P=\left(\begin{smallmatrix}a & b\\ \bar b & \bar a\end{smallmatrix}\right)\in \mathrm{SU}(1,1)$ if $|a|^2-|b|^2=1$ and the M\"obius transformation associated with it is $P\cdot z=\frac{az+b}{\bar b z+\bar a}$. It is a standard result that $\mathrm{SU}(1,1)$ is conjugate to $\SL(2,\R)$ through the $\SL(2,\C)$-matrix $Q=\frac{-1}{1+i}\left(\begin{smallmatrix}1 & -i\\ 1 &  i\end{smallmatrix}\right)$, that is, $Q^*\mathrm{SU}(1,1)Q=\SL(2,\R)$. In fact, we have the following commutative diagram:
\[\xymatrixcolsep{4pc}
\xymatrix{
	\mathbb H\ar[d]_{Q\cdot } \ar[r]^{(Q^{-1}PQ)\cdot } &\mathbb H\ar[d]^{Q\cdot }\\
	\D \ar[r]^{P\cdot} &\D,}
\]
where all transformations are M\"obius transformations, as well as homeomorphisms. Moreover, $Q$ is a homeomorphism between their boundaries, that is, a homeomorphism from $\R\bbP^1 = \R \cup \{\infty\} = \partial \H$ to $S^1 = \partial \D$. We need the following proposition from \cite[Proposition 1]{DE}.

\begin{prop}\label{p:conformal_barycenter}
For each probability measure $\nu$ on the unit circle $S^1$ containing no atom of mass $\ge \frac12$, there is an unique point $B(\nu)\in\mathbb D$, called the conformal barycenter of $\nu$, so that the map $\nu \to B(\nu)$ is invariant under the M\"obius transformation of $\mathrm{SU}(1,1)$, that is, $B(P_*\nu) = P\cdot B(\nu)$ for each $P\in\mathrm{SU}(1,1)$.
\end{prop}

\begin{lemma}\label{l:inv_from_barycenter}
	Let $(\Omega,T,\mu)$, $A$, and $K\subset \Omega$ be as in Lemma~\ref{l.lyap0zp}. Then there exists a family of $A$-invariant, $su$-invariant measures $\hat m_\omega$ over $\mathrm{supp}(\mu|_K)\cap K$ such that for each $\omega \in\mathrm{supp}(\mu|_K)\cap K$, $\hat m_\omega$ is supported by at most two points of $\C\bbP^1$.
	\end{lemma}
\begin{proof}
	We start with the continuous disintegration $\{\tilde m_\omega:\omega\in\Omega\}$ of $m$ over $\omega \in \mathrm{supp} (\mu|_K)\cap K$ that we constructed as stated in Lemma~\ref{l.lyap0zp}. To produce the family of measures $\{\hat m_{\omega}\}$, we divide it into three different cases.
	
	If $\tilde m_\omega$ has an atom $z(\omega)\in\R\bbP^1$ of mass $> 1/2$, we let $\hat m_\omega=\delta_{z(\omega)}$, that is, the Dirac measure (mass one) supported in this point $z(\omega)$. By invariance of $\tilde m_\omega$ under the holonomies, it is clear that if $\tilde m_\omega$ has such a point $z(\omega)$, then so is $m_{\omega'}$ for each point $\omega'$ in $W^s_\loc(\omega)\cup W^u_\loc(\omega)$.  Moreover, $z(\omega')=H^*_{\omega,\omega'}(z(\omega))$ for $*\in\{s,u\}$, which exactly implies that $\delta_{z(\omega)}$ is invariant under the holonomies. Similarly, by invariance of $\tilde m_\omega$ under $A(\omega)$, we have that $\tilde m_{T^n\omega}$ has such a point mass for all $n\in\Z$ and $A(\omega)_*\delta_{z(\omega)}=\delta_{z(T\omega)}$.
	
	If $\tilde m_\omega$ contains two atoms of mass $1/2$ each, we set $\hat m_\omega = \tilde m_\omega$. Similar to the argument of the case (1) above, we have that $m_{\omega'}$ falls into this case for each point $\omega'$ in $W^s_\loc(\omega)\cup W^u_\loc(\omega)\cup \mathrm{Orb}(\omega)$ and $\hat m_\omega$ is invariant under the holonomies and $A(\omega)$.
	
	In all other cases, by Proposition~\ref{p:conformal_barycenter}, we define $\hat m_\omega$ to be the Dirac measure supported at
	$$
	z(\omega):=Q^{-1}\cdot B(Q_*\tilde m_\omega)\in\H,
	$$
	where $B(Q_*\tilde m_\omega)$ is the conformal barycenter of the measure $Q_*\tilde m_\omega$ of the unit circle $S^1$. Note again by holonomy invariance, if $\omega$ is not in the two cases above, then neither is $m_{\omega'}$ for each point $\omega'$ in $W^s_\loc(\omega)\cup W^u_\loc(\omega)\cup \mathrm{Orb}(\omega)$. Moreover, for $\omega'\in W^s_\loc(\omega)$, we have $QH^s_{\omega,\omega'}Q^{-1}\in \mathrm{SU}(1,1)$ which together with Proposition~\ref{p:conformal_barycenter} and the holonomy invariance of $\tilde m_\omega$ implies that
	\begin{align*}
	H^s_{\omega,\omega'}\cdot z(\omega)&= H^s_{\omega,\omega'}\cdot\big(Q^{-1} \cdot B(Q_*\tilde m_\omega)\big)\\
	&=  Q^{-1}(QH^s_{\omega,\omega'}Q^{-1})\cdot B(Q_*\tilde m_\omega)\\
	&=Q^{-1}\cdot B\left((QH^s_{\omega,\omega'}Q^{-1})_*Q_*\tilde m_\omega\right)\\
	&=Q^{-1}\cdot B\left((QH^s_{\omega,\omega'}Q^{-1}Q)_*\tilde m_\omega\right)\\
	&=Q^{-1}\cdot B\left(Q_*(H^s_{\omega,\omega'})_*\tilde m_\omega\right)\\
	&=Q^{-1}\cdot B\left(Q_*\tilde m_{\omega'}\right)\\
	&=z(\omega'),
	\end{align*}
which in turn implies that $\hat m(\omega)$ is invariant under the stable holonomy. By a similar argument we can establish the invariance under the unstable holonomy and under $A(\omega)$.
	\end{proof}

\subsection{Local $su$-Invariance}

In this subsection, we drop the assumption that $\mu$ has bounded distortion and assume only that it has a local product structure. We adapt the techniques from \cite{V} to produce a certain disintegration of $m$ that has local $su$-invariance. Throughout this subsection, we assume $A \in C^\alpha(\Omega,\SL(2,\R))$ to be such that $L(A,\mu) = 0$ and we fix a $\delta$ with $\frac\alpha2 < \delta < \alpha$.

We start with the following consequence of the proof of Lemma~\ref{p.bgvprop3}. Recall that $K_s(N,\delta)$ was defined in \eqref{eq:bunching1}.

\begin{lemma}\label{l:s-holonomy_inv}
	Let $(\Omega,T,\mu)$, $A$, and $\delta$ be above. Then there exists a $\tilde C=\tilde C(\delta,N)$ so that the following holds true. For all $\omega \in K_s(N,\delta)$, all $\tilde\omega \in W^s_\loc(\omega)$, and all $j \ge 0$, we have
	\beq\label{eq:s-holonomy_inv}
	H^s_{T^j\omega,T^j\tilde\omega}:=\lim_{n\to\infty} H^{s,n}_{T^j\tilde\omega,T^j\omega} \mbox{ exists and }\|H^s_{T^j\omega,T^j\tilde\omega}\|\le\tilde C.
	\eeq
\end{lemma}

\begin{proof}
By \eqref{eq:bounded_distortion}, we have that
$$
\|H^s_{\tilde\omega,\omega}\|=\|\lim_{n\to\infty}H^{s,n}_{\tilde\omega,\omega}\|\le e^{CN}.
$$
A direct computation shows that
$$ H^{s,n}_{T^j\tilde\omega,T^j\omega}=A_j(\tilde\omega)H^{s,n+j}_{\tilde\omega,\omega}A_j(\omega)^{-1},
$$
which implies the existence of
$$
H^{s}_{T^j\tilde\omega,T^j\omega}:=\lim_{n\to\infty}H^{s,n}_{T^j\tilde\omega,T^j\omega}\mbox{ and }\|H^{s}_{T^j\tilde\omega,T^j\omega}\|\le \|H^{s,n}_{T^j\tilde\omega,T^j\omega}\|\cdot\|A_j(\omega)\|\cdot\|A_j(\tilde\omega)\|.
$$
In particular, for all $0\le j\le N$, we have
\beq\label{eq:s-holonomy_inv1}
\|H^{s}_{T^j\tilde\omega,T^j\omega}\|\le e^{3CN}.
\eeq
Fix a $j>N$. Then we have $\|A_j(\omega)\|<e^{(\alpha-\delta)j}$ since $\omega\in K_s(N,\delta)$. Using \eqref{eq:bdd_difference}, a direct computation shows that
	\begin{align*}
	\|\delta^{s,n}_{T^j\tilde\omega,T^j\omega}\|
	& =\|A_j(\omega)\delta^{s,n+j}_{\tilde\omega,\omega}A_j(\omega)^{-1}\|\\
	& \le Ce^{(\alpha-\delta)j}e^{(\alpha-\delta)(n+j)}e^{-\alpha(n+j)}\\
	& = Ce^{(\alpha-2\delta)j}e^{-\delta n} \\
    & < Ce^{-\delta n},
	\end{align*}
	where the last inequality follows from the fact $2\delta>\alpha$. Combining \eqref{eq:s-holonomy_inv} and the proof of \eqref{eq:bounded_distortion}, we obtain for all $j\ge N$
	\beq\label{eq:s-holonomy_inv2}
		\|H^{s}_{T^j\tilde\omega,T^j\omega}\|\le H^{s,N}_{T^j\tilde\omega,T^j\omega} \cdot \exp\left(C\sum^{\infty}_{n=N}\|\delta^{s,n}_{\omega,\tilde\omega}\|\right)\le Ce^{2CN}.
	\eeq
	Combining \eqref{eq:s-holonomy_inv1} and \eqref{eq:s-holonomy_inv2}, we clearly obtain the latter half of \eqref{eq:s-holonomy_inv}, where we may take $\tilde C=\max\{e^{3CN}, Ce^{2CN}\}$.
\end{proof}

So we may choose $N$ large so that $K(N,\delta)=K_s(N,\delta)\cap K_u(N,\delta)$ (which were defined in \eqref{eq:bunching1}) has measure sufficiently close to $1$ and $\mu(K(N,\delta)\cap [0;j])>0$ for all $1\le j\le \ell$. We set $K^j_\tau:=K_\tau(N,\delta)\cap [0;j]$ for $\tau\in\{s,u\}$ and $K^j=K^j_s\cap K^j_u$.

\begin{lemma}\label{l:loc_invPrinciple}
	Let $(\Omega,T,\mu)$, $A$, $\delta$ be as in Lemma~\ref{l:s-holonomy_inv} and let $K(N,\delta)$ be as above. Then for every $(T,A)$-invariant measure $m$ that projects to $\mu$ on the first component, there is a disintegration $\{m_\omega: \omega\in\Omega\}$ of $m$ that is $su$-invariant for $\mu$-almost every $\omega\in K(N,\delta)$.
\end{lemma}

\begin{proof}
	We only consider the case of $s$-invariance, as $u$-invariance can be established in a completely analogous way. We break the argument into three steps.
	\vskip .1cm
	
	\noindent \textit{Step I.} As in the proof of Proposition~\ref{p.bgvprop1}, the first step is to construct a certain $\sigma$-algebra $\CB$ to which we can apply the invariance principle as formulated in Proposition~\ref{p:invPrinciple}. By \eqref{eq:bounded_distortion}, $K_s(N,\delta)$ is $s$-saturated. Similarly, $K_u(N,\delta)$ is $u$-saturated. Fix a $\omega^j\in K^j$ and set $S=W^u_\loc(\omega^j)\cap K^j$. For each $\omega'\in S$, we define $r(\omega')=1$ if $T(W^s_\loc(\omega'))\cap W^s_\loc(\omega'')=\varnothing$ for some $\omega''\in S$; otherwise, we define $2\le r(\omega')\in \Z_+\cup \{\infty\}$ be the largest number such that $T^i(W^s_\loc(\omega'))\cap W^s_\loc(\omega'') = \varnothing$ for all $\omega''\in S$ and for all $0<i< r(\omega')$. Now we define the $\sigma$-algebra $\CB \subseteq \CM$ to be the one generated by the family
	$$
    \{T^i(W^s_\loc(\omega')): \omega'\in S,\ 0\le i< r(\omega')\}.
	$$
	By our definition of $r(\omega')$, it is clear that the sets in the family above are mutually disjoint. Thus, $\CB$ contains all $B\in\CM$ such that for all $\omega'\in S$ and all $0\le i<r(\omega')$, either $B\cap T^i(W^s_\loc(\omega'))=\varnothing$ or $T^j(W^s_\loc(\omega'))\subset B$. First we claim that $\CB$ satisfies condition (1) of Proposition~\ref{p:invPrinciple}. The proof is analogous to the one of Proposition~\ref{p.bgvprop1}. Indeed, $T\CB$ is the $\sigma$-algebra generated by
	$$
	   \{T^{i+1}(W^s_\loc(\omega')): \omega'\in S,\ 0\le i< r(\omega')\},
	$$
	which is again a family of mutually disjoint sets. Since $T^{r(\omega')}(W^s_\loc(\omega')) \subseteq W^s_\loc(\omega'')$ for some $\omega''\in S$, one readily checks that $B\in\CB$ implies $B\in T\CB$. Hence, we have that $T\CB$ contains $\CB$, or equivalently, $T^{-1} \CB \subseteq \CB$. More generally, for all $n\ge 1$, we have that $T^n\CB$ is generated by $\{ T^{i+n} (W^s_\loc(\omega')) : \omega'\in S,\ 0\le i < r(\omega')\}$, which implies that $T^n \CB,\ n\ge 1$ generates $\CM$ mod $0$. Indeed, since $\CM$ is generated by cylinders, we just need to show that all $[k;\vec l]$ are contained in $T^n\CB$ for some large $n$. Taking any $n\ge |k|$, it is clear that $[k;\vec l]\in T^n\CB$.
	\vskip .2cm
	
	\noindent \textit{Step II}. Similarly to the proof of Proposition~\ref{p.bgvprop2}, our second step is to conjugate $A$ to some $\tilde A$, which is measurable with respect to $\CB$. We define $\tilde A$ by
	$$
	\tilde A(\omega):=H^s_{T\omega,T^{i+1}\omega'}A(\omega)H^s_{T^i\omega',\omega}=A(T^i\omega')
	$$
	if $\omega\in T^{i}(W^s_\loc(\omega'))$ for some $\omega'\in S$ (so that $\omega\in W^s_\loc (T^i\omega')$ ) and $0\le i<r(\omega')$;  and
    $$
	\tilde A(\omega):=A(\omega) \mbox{ otherwise}.
	$$
	Clearly, if we set $B(\omega)$ as
	$$
	B(\omega)=
	\begin{cases}H^s_{\omega,T^j\omega'}, &\omega\in T^j(W^s_\loc(\omega')), \omega'\in S, \mbox{ and }0\le j<r(\omega')\\ I_2, &\mbox{otherwise},
	\end{cases}
	$$
	then $\tilde A(\omega)=B(T\omega)A(\omega)B(\omega)^{-1}$. In other words, $\tilde A$ is conjugate to $A$ via $B$. Combining this with the fact that $\omega'\in S\subseteq K_s(N,\delta)$ and Lemma~\ref{l:s-holonomy_inv}, we have $\|B(\omega)\|\le \tilde C$ for all $\omega\in\Omega$. In particular, $\int \log \|\tilde A\| \, d\mu < \infty$ and $L(\tilde A,\mu)=0$. By definition, $\tilde A$ is constant on $T^i(W^s_\loc(\omega'))$ for any $\omega'\in S$ and any $0\le i<r(\omega')$, which clearly implies that $\tilde A$ is $\CB$-measurable.
	\vskip .2cm
	
	\noindent \textit{Step III}. Following the second half of the proof of Proposition~\ref{p.bgvprop2}, for any given disintegration $\{m_\omega:\omega\in\Omega\}$ of $m$, we can set
	$$
	\tilde m_\omega=B(\omega)^{-1}_*m_\omega,\ \omega\in\Omega.
	$$
	Then it becomes a disintegration of a $(T,\tilde A)$-invariant measure $\tilde m$.  We can now apply Proposition~\ref{p:invPrinciple} to $(\CB, \tilde A, \tilde m_\omega)$ and obtain that $\{\tilde m_\omega\}$ is $\CB$-measurable. In particular, $\tilde m_\omega$ is constant on $T^i(W^s_\loc(\omega'))$ for all $\omega'\in S$ and all $0\le i\le r(\omega')$. Taking $i=0$, then  a similar proof to \eqref{eq:s-inv} yields
	 $$
	(H^s_{\omega,\tilde\omega})_*m_\omega=m_{\tilde\omega}\mbox{ for all } \omega,\tilde\omega\in W^s_\loc(\omega') \mbox{ and all }\omega'\in S.
	$$
	Thus we have obtained $s$-invariance for all points in
	$$
	\bigcup_{\omega\in S}\{W^s_\loc(\omega):\omega\in W^u_\loc(\omega^j)\},
	$$
	which contains
	$$
	\bigcup_{\omega\in W^u_\loc(\omega^j)\cap K^j}\{W^s_\loc(\omega):\omega\in W^u_\loc(\omega^j)\cap K^j\}.
	$$
	By local product structure of $\mu$ and following the proof of \eqref{eq:fullmeasure}, we have that the set above is a full measure subset of $K^j$. Since $1\le j\le \ell$ is arbitrarily chosen, we thus obtain $s$-invariance of $\{m_\omega\}$ on a full measure subset of $K(N,\delta)$.
\end{proof}

Now we can apply the proof of  Lemmas~\ref{l.lyap0zp} and ~\ref{l:inv_from_barycenter} (replacing $K_\delta$ by $K(N,\delta)$) to obtain the following corollary.

\begin{coro}\label{c:local_su_inv}
	Using the setup of Lemma~\ref{l:loc_invPrinciple}, there is a disintegration $\{\tilde m_\omega\}$ of $m$ so that $\omega \mapsto \tilde m_\omega$ is continuous and $su$-invariant on $\mathrm{supp}(K(N,\delta))\cap K(N,\delta)$. Moreover, there is family of measures $\{\hat m_\omega\}$ that is $su$-invariant on $\mathrm{supp}(K(N,\delta))\cap K(N,\delta)$ and for each $\omega$, $\mathrm{supp}(\hat m_\omega)$ contains at most two points.
\end{coro}

The main goal of the present Section~\ref{s:invPrinciple_Barycenter} is to obtain the following corollary.

\begin{coro}\label{c.lyap0zp}
Suppose $(\Omega,T)$ is a subshift of finite type and $\mu$ is a $T$-ergodic measure that has a local product structure. Let $A : \Omega \to \SL(2,\R)$ be a cocycle map so that $L(A,\mu)=0$. Then for every periodic point $p$ (of period $n$) such that $2L(A,p) < \frac\alpha2$, there exists a set $Z_p\subseteq \C\bbP^1$, invariant under complex conjugation and under $A_n(p)$, and consisting of either one or two points, with the following property.   Let $q$ be another periodic point such that $2L(A,p)<\frac\alpha2$. If $p_0=q_0$, then
	$$
	H^u_{q,q\wedge p}(Z_q)=H^s_{p, q\wedge p}(Z_p).
	$$
\end{coro}

\begin{proof}
	Since $L(A,\mu)=0$ and $L(A,p)<\frac\alpha2$, by Lemma~\ref{l:PinSuppK}, we clearly have that $p \in \mathrm{supp}(\mu|_{K(N,\delta)})\cap K_{\delta}(N,\delta)$ for some $\frac\alpha2<\delta<\alpha$. Thus we may apply Corollary~\ref{c:local_su_inv} to obtain the measure $\hat m_p$ defined at $p$ which is $A_n(p)$-invariant and $su$-invariant. Hence, if we define
	$$
	Z_p:=\{z(p), \overline{z(p)}:\ z(p)\in\mathrm{supp}(\hat m_{p})\},
	$$
	then $Z_p$ consists of at most two points. Moreover, it is clear that $Z_p$ is $su$-invariant, $A(p)$-invariant, and invariant under complex conjugation.
	
	Since $q$ is also a periodic point such that $L(A,q)<\frac\alpha2$, we can  certainly find $\frac\alpha2<\delta<\alpha$ so that both $p$ and $q$ belong to $\mathrm{supp}(\mu|_{K(N,\delta)})\cap K(N,\delta)$. Thus $Z_p$ and $Z_q$ are both defined. By $su$-invariance of $Z_p$ and $p_0=q_0$, we have
	$$
	H^u_{q,q\wedge p}(Z_q)=H^s_{p, q\wedge p}(Z_p)
	$$
as desired.
\end{proof}

\section{Positivity of the Lyapunov Exponent I}\label{s:PosLya1}

Throughout this section we assume that $\Omega \subseteq \mathcal{A}^\Z$ is a subshift of finite type and $\mu$ is a $T$-ergodic probability measure that is fully supported on $\Omega$ and has a local product structure. We fix a non-constant $f\in C^\alpha(\Omega,\R)$ and consider the one-parameter family of Schr\"odinger cocycles $(T,A^E)$. We shall apply the techniques from Section~\ref{s:invPrinciple_Barycenter} to study the positivity property of the Lyapunov exponent.

In Subsection~\ref{ss.nonundom}, we show under a very general condition that the set of energies with zero Lyapunov exponent is a discrete set.  In Subsection~\ref{ss.glundom}, we apply the same techniques to the scenario where we have global existence of holonomies and obtain a stronger result for the corresponding Schr\"odinger cocycles. Namely, we show that under the same general condition, the set of energies with zero Lyapunov exponent is a finite set. Global existence of the holonomies may be obtained if the $\|\cdot\|_\infty$ norm of the potential is small or if the potentials are locally constant.

\subsection{General Case: Positivity Away from a Discrete Set}\label{ss.nonundom}

Throughout this subsection, we assume that $E_0 \in \R$ is an accumulation point of $\CZ_f = \{ L(E) = 0 \}$. Clearly,
\begin{equation}\label{e.e0insigma}
E_0 \in \Sigma,
\end{equation}
since $\CZ_f \subseteq \Sigma$ and $\Sigma$ is closed.

\begin{definition}
We say that a periodic point $q$ is $\gamma$-\emph{bunched} at $E$ if
$$
2L(A^{E},q) < \gamma \le \alpha.
$$
\end{definition}

We fix $p$ to be a periodic point that is $\frac\alpha2$-bunched at $E_0$. We let $n_p$ denote the period of $p$.

\begin{lemma}\label{l:zeroLE_bunched_p}
$L(A^{E_0},p)=0$. In particular, $E_0\in\sigma(H_{p})$.
\end{lemma}

\begin{proof}
Let $E_n \to E_0$, $E_n \neq E_0$ be a sequence in $\CZ_f$. Recall from \eqref{e.e0insigma} that $E_0\in\Sigma$.

Assume that $p$ is hyperbolic for $E_0$. Then $p$ is still hyperbolic and $\frac\alpha2$-bunched in a small neighborhood $J$ of $E_0$. Let $E_0\in\sigma(H_{\omega'})$ for some $\omega'\in\Omega$. By Proposition~\ref{p:denseOrb_to_denseSpec} and by choosing $\delta>0$ small, we have $J\cap \sigma(H_{\omega})\neq\varnothing$ for any $\omega$ such that $\mathrm{Orb}(\omega)\cap B_{\delta}(\omega')\neq\varnothing$. On the other hand, by Proposition~\ref{p:specification} there is an $r = r(\delta)$ so that for any $I_1=[0,n_1]\subseteq\Z$, there is a periodic orbit $q$ with period $n_q=n_1+r+1$ so that $d(T^jq,T^jp)<\delta$ for all $1\le j\le n_1$ and $d(T^{n_1+r+1}q,\omega')<\delta$. An immediate consequence is that $\sigma(H_q)\cap J\neq \varnothing$. Moreover, as $n_1$ goes to infinity, it clearly holds that $L(A^E,q)$ tends to $L(A^E,p)$ uniformly for all $E\in J$. In particular, by choosing $n_1$ large, we have that $q$ is $\frac\alpha2$-bunched for all $E\in J$ as well. We fix such a periodic point $q$.

Clearly, $p_0=q_0$. Thus we may define
$$
H^E:=H^{u,E}_{q\wedge p, q}\cdot H^{s,E}_{p,q\wedge p}
$$
for each $E\in J$. Here $H^{s,E}_{p,q\wedge p}$ and $H^{u,E}_{q,q\wedge p}$ are the holonomies corresponding to $A^{E}$, which are well-defined since both $p$ and $q$ are $\frac\alpha2$-bunched through $J$. Moreover, they are holomorphic on $J$ since they are limits of uniformly convergent sequences of holomorphic functions $H^{s,n}(E)$ or $H^{u,n}(E)$ on $J$. Thus we have that $E \mapsto H^E$ is analytic. Let $Z_p=Z_p(E_n)$ be as in Corollary~\ref{c.lyap0zp}. By passing to a subsequence, we may assume that
$$
Z_p(E_n)=
\begin{cases}
\{s(E_n)\}\ &\mbox{ for all }n,\\
\{u(E_n)\}\ &\mbox{ for all }n, \mbox{ or }\\
\{s(E_n),u(E_n)\}\ &\mbox{ for all }n.
\end{cases}
$$
Thus, we may extend the definition of $Z_p(E)$ to all $E \in J$ so that $Z_p(E)$ consists of one or two functions that are analytically on $J$. By Corollary~\ref{c.lyap0zp},
$$
Z_q(E):= H^E\cdot Z_p(E)
$$
is invariant under the monodromy matrix of $q$ for infinitely many $E_n$.  By analyticity, it follows that $Z_q(E)$ is invariant by the monodromy matrix of $q$ for every $E \in J$.  Since $Z_p(E)$ is real, so is $Z_q(E)$, which implies that the absolute value of the trace of the monodromy matrix of $q$ cannot become smaller than $2$ throughout $J$.  Since $J$ is open and $E\in \sigma(H_q)$ cannot be an isolated point, we must have $J\cap \sigma(H_q)=\varnothing$, which contradicts our choice of $q$. It follows that $p$ is not hyperbolic for $E_0$. In particular, $E_0 \in \sigma(H_p)$.
\end{proof}

Let $n_p$ be the period of $p$. Choose a small open disk $D \subseteq \C$ around $E_0$ such that $p$ is $\frac\alpha2$-bunched for all energies $E$ in the closed disk $\bar D$.  Recall that by Proposition~\ref{p:spectrumbands}, $\Delta(E)=\tr(A^E_{n_p}(p))$ is monotonic on each connected component of $\Delta^{-1}(-2,2)$. Thus we may also assume that $D$ is small enough so that, through $\bar D \setminus \{E_0\}$, $\Delta(E)$ is different from $-2,2,0$. According to Subsection~\ref{s:periodic}, if $E_0\notin\partial (\sigma(H_p))$, we can then define two holomorphic functions $u,s: D \to \C \bbP^1$, distinct everywhere, such that $u(E)$ and $s(E)$ are eigendirections of $A^{E}_{n_p}(p)$;  otherwise we can still define holomorphic functions $u,s$ on the ramified (at $E_0$) double cover of $\pi: \tilde D\to D$, giving (distinct) eigendirections when $\tilde E \in \tilde D \setminus \{E_0\}$, but taking as value at $E_0$ the single real eigendirection of $A^{E_0}_{n_p}(p)$. Moreover, for $\pi(\tilde E)\in \R$, $s(\tilde E)$ and $u(\tilde E)$ are real if and only if $\pi(\tilde E)$ not in the interior of $\sigma(H_p)$.

\begin{lemma}\label{claim1}
If $q$ is a periodic point that is $\frac\alpha2$-bunched through $E \in \bar D$, then $\sigma(H_p) \cap D = \sigma(H_q) \cap D$.
\end{lemma}

\begin{proof}
Similarly to the proof of Lemma~\ref{l:zeroLE_bunched_p}, we define $H^E=H^{u,E}_{q\wedge p, q}\cdot H^{s,E}_{p,q\wedge p}$ for $E \in D$. Since $E_0\in \sigma(H_p)$, we consider two different cases.

If $E_0\notin\partial (\sigma(H_p))$, then $D\cap \sigma(H_p) = D\cap \R$ by our choice of $D$, and $Z_p(E_n)$ is a nonempty subset of $\{u(E_n),s(E_n)\}$. Following the same argument that showed that $A^{E_0}_{n_q}(q)$ has a real eigendirection in the proof of Lemma~\ref{l:zeroLE_bunched_p}, we obtain that $A^{E}_{n_q}(q)$ of the present lemma has a non-real eigendirection for all $E \in D\cap \R$. This implies that $D\cap \R \subseteq \sigma(H_q)$, and the claim follows in this case.

If $E_0 \in \partial (\sigma(H_p))$, then by our choice of $D$ we have $D\cap \Sigma_p$ is either $[E_0,E_+)$ or $(E_-,E_0]$ where
$$
(E_-,E_+) = D\cap \R.
$$
For simplicity, we will assume that $\inter D \cap \Sigma = [E_0,E_+)$. Recall $\pi:\tilde D\to D$ is the double cover map of $D$ ramified at $E_0$. For each $n$, choose a preimage $\tilde E_n \in \pi^{-1}(E_n)$. Then $Z_p(E_n)$ is a subset of $\{ u(\tilde E_n), s(\tilde E_n) \}$. As in the proof of Lemma~\ref{l:zeroLE_bunched_p} and up to replacing $E_n$ by a subsequence, $Z_p(E_n)$ is always of the form $\tilde Z_p(\tilde E_n)$, where
$$
\tilde Z_p(\tilde E)=
\begin{cases}
\{s(\tilde E)\}\ &\mbox{ for all }\tilde E\in\tilde D,\\
\{u(\tilde E)\}\ &\mbox{ for all }\tilde E\in\tilde D, \mbox{ or }\\
\{s(\tilde E), u(\tilde E)\}\ &\mbox{ for all }\tilde E\in\tilde D.
\end{cases}
$$
Notice that if $\pi(\tilde E) \in (E_-,E_0)$, then $\tilde Z_p(\tilde E)$ consists of real directions; if $\pi(\tilde E) \in (E_0,E_+)$, then $\tilde Z_p(\tilde E)$ consists of non-real directions. We again define
$$
\tilde Z_q(\tilde E) := H^{\pi(\tilde E)} \cdot \tilde Z_p(\tilde E).
$$
Then $ \tilde Z_p(\tilde E)$ is invariant under $A^{\pi(\tilde E)}_{n_q}(q)$ whenever $\tilde E = \tilde E_n$. By the fact that $H^{\pi(\cdot)}$, $u$, and $s$ are all holomorphic on $\tilde D$, it follows that  $\tilde Z_p(\tilde E)$ is invariant under $A^{\pi(\tilde E)}_{n_q}(q)$ for all $\tilde E\in\tilde D$. This implies that $A^{E}_{n_q}(q)$ has at least one real eigendirection for $E \in (E_-,E_0)$ and has at least one non-real eigendirection for $E \in (E_0,E_+)$. This can only happen when $D \cap \sigma(H_q) = [E_0,E_+)$, and the claim follows in this case.
\end{proof}

\begin{lemma}\label{claim2}
If $q$ is any periodic point, then $\sigma(H_q) \cap D = \sigma(H_p)  \cap D$.
\end{lemma}

\begin{proof}
Fix an arbitrary periodic point $q_0$. Let us say that a periodic point $q$ is $(\epsilon, \delta)$-good, $0 < \delta < \epsilon < 1$, if it spends at least a $1 - \epsilon$ proportion of its iterates within distance $\delta$ of $p$, and at least a $\epsilon/2$ proportion of its iterates within distance $\delta$ of $q_0$. By Proposition~\ref{p:specification} and similar to the argument leading to the choice of $q$ in the proof of Lemma~\ref{l:zeroLE_bunched_p}, we see that the set of $(\epsilon,\delta)$-good periodic points is not empty for any choice of $0<\delta<\epsilon<1$. Moreover, if $\epsilon$ is sufficiently small, then an $(\epsilon,\delta)$-good $q$ is $\frac\alpha2$-bunched for energies $E \in \bar D$. By Lemma~\ref{claim1}, it then holds that $\sigma(H_p) \cap D = \sigma(H_q) \cap D$ for all such $q$'s. We fix such a small $\epsilon$ for the remainder of this proof.
	
First we show $\sigma(H_{q_0}) \cap D \subseteq \sigma(H_p)  \cap D$. If this is not true, then there is some $E_0\in(\sigma(H_{q_0})\cap D)\setminus \sigma(H_{p})$. In particular, we have
$$
\varepsilon := \min \{ d(E_0,\sigma(H_{p})), d(E_0,\partial D) \} > 0.
$$
Then for an $(\epsilon, \delta)$-good periodic point $q$, we also have that
$$
\varepsilon = \min \{ d(E_0,\sigma(H_{q})), d(E_0,\partial D) \} > 0.
$$
By Proposition~\ref{p:denseOrb_to_denseSpec} and the fact that $\mathrm{Orb}(q)\cap B_\delta(q_0)\neq \varnothing$ for $(\epsilon, \delta)$-good points, we have for sufficiently small $\delta$ and an $(\epsilon,\delta)$-good periodic point $q$ that
$$
\sigma(H_{q_0}) \subseteq B_{\frac\varepsilon2}(\sigma(H_{q})).
$$
Clearly, this implies that $d(E_0,\sigma(H_{q})) < \frac\varepsilon2$ and we obtain a contradiction. So the first part follows.

Now we show $\sigma(H_p) \cap D \subseteq \sigma(H_{q_0})  \cap D$. Suppose this is not the case. By the first part, there is an $E_0\in (\sigma(H_{p})\cap D) \setminus \sigma(H_{q_0})$. In particular, we have
$$
\tilde\varepsilon:=\min\{d(E_0,\sigma(H_{q_0})), d(E_0,\partial D)\}>0.
$$
Notice that $\mathrm{Orb}(q_0) \cap B_\delta(T^mq)\neq \varnothing$ for an $(\epsilon,\delta)$-good periodic point $p$ and for some $m\in\Z$. Thus by Proposition~\ref{p:denseOrb_to_denseSpec}, and by choosing $\delta$ small (the smallness of which is independent of $q$ or $q_0$), we have
$$
\sigma(H_{T^mq}) \subseteq B_{\frac{\tilde\varepsilon}2}(\sigma(H_{q_0})).
$$
Since $\sigma(H_{T^mq}) = \sigma(H_{q})$ and $E_0 \in \sigma(H_{p})\cap D = \sigma(H_{q})\cap D$, we obtain
$$
d(E_0,\sigma(H_{q_0}))<\frac{\tilde\varepsilon}2,
$$
which is a contradiction and the lemma follows.
\end{proof}

By \eqref{eq:spectrum_variousForms}, the spectrum $\Sigma$ is the closure of the union of the spectra of periodic points. Thus Lemma~\ref{claim2} implies that $\Sigma \cap D = \sigma(H_p) \cap D$. For each $T$-ergodic measure $\nu$ on $\Omega$, we let $\Sigma_\nu$ denote the set such that $\sigma(H_\omega)=\Sigma_\nu$ for $\nu$ almost every $\omega$; see, for example, \cite{pastur}.

\begin{lemma}\label{claim3}
For any $T$-ergodic measure $\nu$ on $\Omega$, we have $\Sigma_\nu\cap D=\sigma(H_p) \cap D$. Moreover,  $L(A^E;\nu)=0$ for all $E \in \sigma(H_p) \cap D$. In particular, $L(E)=0$ for all such $E$'s.
\end{lemma}

\begin{proof}
Since  we have $\sigma(H_\omega)\subseteq\Sigma$ for each $\omega\in\Omega$ and $\Sigma \cap D = \sigma(H_p) \cap D$, it clearly holds that $\Sigma_\nu\cap D\subseteq \sigma(H_p)  \cap D$. On the other hand, if $E\notin\Sigma_\nu$, then the sequence $\{A^E(T^n\omega)\}_{n\in\Z}$ is uniformly hyperbolic for $\nu$-almost every $\omega\in\Omega$, which in turn implies that $L(A^E;\nu)>0$; see, for example, \cite[Theorem 3]{zhang}. Thus $L(A^E;\nu)=0$ implies that $E\in\Sigma_\nu$. So we only need to prove the second part of the lemma.

Assume that the statement is false. In other words, we have $L(A^E;\nu)>0$ for some $E\in \sigma(H_p) \cap D$. By \cite[Theorem 3]{kalinin}, for each $\epsilon>0$, there is a periodic point $q\in\Omega$ so that $|L(A^E;\nu)-L(A^E,q)|<\epsilon$. Thus there is periodic point $q\in\Omega$ so that $L(A^E,q)>0$. In particular, $E\notin \sigma(H_q)$, which contradicts Lemma~\ref{claim2}, concluding the proof.
\end{proof}

\begin{lemma}\label{l.accumulation}
	Let $E_0$ be an accumulation point of $\{L(E)=0\}$.  Assume that there exists a periodic point $p$ that is $\frac\alpha2$-bunched at $E_0$.  Then
	\begin{enumerate}
		\item The connected component $I$ of $E_0$ in the spectrum is isolated,
		\item $L(A^E;\nu)=0$ for all $E\in I$ and all $T$-ergodic measure $\nu$ on $\Omega$.
	\end{enumerate}
\end{lemma}

\begin{proof}
Let $I$ be the connected component of $\sigma(H_p) $ that contains $E_0$. Notice that $p$ is $\frac\alpha2$-bunched for every $E \in I$ since $L(A^E,p)=0$ on $I$. Let $S \subseteq I$ be the set of accumulation points of $\{ L(E) = 0 \}$. It is clearly a closed and non-empty set since $E_0\in S$. Moreover, applying Lemma~\ref{l.accumulation} to all $E \in S$, we see that $S$ is open in $I$.  So $S = I$. Clearly $I\subseteq\Sigma$. Applying the fact $\Sigma\cap D=\sigma(H_p)\cap D$ for some disk $D$ around the boundary points of $I$,  we conclude that $I$ is an isolated component of $\Sigma$. Applying Lemma~\ref{l.accumulation} to all $E \in I$ again, we obtain that $L(A^E;\nu)=0$ for all $E\in I$ and for all $T$-ergodic measure $\nu$ on $\Omega$.
\end{proof}

We have now collected all the tools to prove our main theorem.

\begin{proof}[Proof of Theorem~\ref{t:posLEawayDiscreteSet}]
Suppose to the contrary that there are $E_0 \in \{ E : L(E) = 0 \}$ and $E_n \in \{ E : L(E) = 0 \} \setminus \{ E_0 \}$, $n \in \Z_+$, such that $E_n \to E_0$ as $n \to \infty$. Since $L(E_0) = 0$, we can choose a $\frac\alpha2$-bunched periodic point by \cite[Theorem 3]{kalinin}. It now follows from Lemma~\ref{l.accumulation} that $E_0$ belongs to a non-degenerate compact interval $I$, which is a connected component of $\Sigma$, as well as of all periodic spectra $\sigma(H_p)$. In particular, for the fixed point of $T$, the unique connected component of its spectrum is an interval of length $4$. Since having such a connected component is only possible for constant periodic potentials by Proposition~\ref{p.perinvmeas}, it follows that the potential associated with each periodic point must be constant. This implies that $f$ itself must be constant; contradiction.
\end{proof}

\subsection{Special Cases: Positivity Away from a Finite Set}\label{ss.glundom}

In this subsection we consider sampling functions $f : \Omega \to \R$ for which we have global existence of the holonomies in the sense that the cocycle $A^E$ admits \emph{canonical holonomies} as defined in Subsection~\ref{sss.holonomies} for all $E$ in a complex neighborhood of the convex hull $\CU\subseteq\C$ of the spectrum $\Sigma_f$. Since $A^E$ depends on $E$ holomorphically, we obtain that the holonomies are holomorphic on $\CU$ as well. In this case, we are able to improve the result we obtained in Subsection~\ref{ss.nonundom}.

There are two types of $f$ for which we have such global existence of holonomies. One is the set of $f \in C^\alpha(\Omega,\R)$ for which $A^{E}$ is fiber bunched for every $E$ in the convex hull of the spectrum $\Sigma$. Such an $f$ will be called \textit{globally bunched}. The other is the set of locally constant $f$'s.

\subsubsection{Fiber Bunching and Existence of Holonomies}

Assume that $(\Omega,T,\mu)$ is subject to the same assumptions as before.

\begin{definition}
We say that $A \in C^\alpha (\Omega , \mathrm{SL}(2,\C))$ is \emph{fiber bunched} if there exists $n_0 \ge 1$ such that for every $\omega \in \Omega$, we have
\beq\label{eq:bunching2}
\|A_{n_0}(\omega)\|^2 < e^{\alpha n_0}.
\eeq
Equivalently, there is $\theta < \alpha$ such that $\|A_{n_0}(\omega)\|^2 < e^{\theta n_0}$ for every $\omega \in \Omega$.
\end{definition}

Note that fiber bunching is clearly a $C^0$-open condition. A fiber bunched cocycle has canonical holonomies as defined in Subsection~\ref{sss.holonomies}. In fact, we can run the proof of Lemma~\ref{p.bgvprop3} to show that for $\omega' \in W^u_\loc(\omega)$, $H^{u,n}_{\omega,\omega'} = A_{-n}(\omega')^{-1} \cdot A_{-n}(\omega)$ converges uniformly on $\Omega$ to the stable holonomy $H^u_{\omega,\omega'}$. Similarly for $\omega' \in W^s_\loc(\omega)$, we have that $H^{s,n}_{\omega,\omega'} = A_n(\omega')^{-1} A_n(\omega)$ converges uniformly on $\Omega$ to the unstable holonomy $H^s_{\omega,\omega'}$. Indeed, to obtain the uniform convergence to holonomies,  the only condition we used in the proof of Lemma~\ref{p.bgvprop3} is the condition in \eqref{eq:bunching1}, which is exactly the fiber bunching condition \eqref{eq:bunching2}. We also note the following: If $A^t \in C^\alpha(\Omega,\SL(2,\C))$, $t$ in some domain $U\subseteq \C$, is a continuous family such that $t \mapsto A^t(\omega)$ is holomorphic for every $\omega \in \Omega$ and $A^t$ is fiber bunched for every $t$, then the stable and unstable holonomies depend holomorphically on $t$. Indeed, in this case, the holonomies are limits of uniformly convergent sequences of holomorphic functions. In particular, we may consider Schr\"odinger cocycles $A^E$ with sampling function $f \in C^\alpha(\Omega,\R)$. If $\|f\|_{\infty}$ is sufficiently small, then $A^E$ is fiber bunched in a complex neighborhood of the convex hull of the spectrum $\Sigma$. To see this, we first see that $\left(\begin{smallmatrix} E & -1 \\ 1 & 0 \end{smallmatrix}\right)$ is fiber bunched for all $E\in[-2,2]$ since they are all elliptic or parabolic. By openness of fiber bunching, we then have that $A^{E}$ is fiber bunched for all $E$ in a complex neighborhood of $[-2,2]$ provided $\|f\|_{\infty}$ is sufficiently small. If necessary, we can then choose $\|f\|_{\infty}$ smaller so that the convex hull of $\Sigma_f$ is contained in such an open neighborhood. Thus $f$ is globally bunched.

\subsubsection{Locally Constant Cocycles}

The other class for which the canonical holonomies exist for obvious reasons is defined as follows.

\begin{definition}
We say that $A:\Omega\to\mathrm{SL}(2,\R)$ is \emph{locally constant} if there exists a $n_0$ such that for each $\omega\in\Omega$, $A(\omega)$ depends only on the cylinder set $[-n_0;\omega_{-n_0}, \ldots, \omega_{n_0}]$.
\end{definition}

Evidently, locally constant cocycles are $\alpha$-H\"older continuous for all $\alpha>0$. Locally constant cocycles might not be fiber bunched. However, the holonomies exist trivially. Indeed, if $A$ is locally constant, then there is a $n_0\in\Z_+$ so that for all $\omega$ and all $n>n_0$ we have
$$
H^{\tau,n}_{\omega,\omega^\tau}=H^{\tau,n_0}_{\omega,\omega^\tau},
$$
where $\tau\in\{s,u\}$ and $\omega^\tau\in W^\tau_\loc(\omega)$. Thus $H^{\tau,n_0}_{\omega, \omega^\tau}$ are exactly the holonomies. Now we consider Schr\"odinger cocycles $A^E$ with potential $f : \Omega \to \R$. If there is a $n_0\in\Z_+$ such that $f(\omega)$ depends only on $[-n_0;\omega_{-n_0}, \ldots, \omega_{n_0}]$, then $A^E$ is locally constant for all $E \in \C$. In other words, a locally constant sampling function induces locally constant Schr\"odinger cocycle maps.

\subsubsection{Energies Admitting an $su$-State}

Again, our objective is to study the energies for which $L(E) = L(A^{E},\mu) = 0$. We will point out how the desired statements will follow by simple specialization of the proofs of the lemmas in Section~\ref{s:invPrinciple_Barycenter}.

Assume that $A \in C^0(\Omega, \mathrm{SL}(2,\R))$ has canonical holonomies $H^\tau_{\omega, \omega'}$, where $\tau \in \{s,u\}$. Recall that an $su$-state for $A$ is a $(T,A)$-invariant measure $m$ with a disintegration $\{m_\omega: \omega\in\Omega\}$ that is invariant under the cocycle and the holonomies. In particular, for $\mu$-almost every $\omega\in\Omega$, we have
\begin{enumerate}
	
	\item $A(\omega)_* m_\omega = m_{T \omega}$,
	
	\item $(H^s_{\omega,\omega'})_* m_\omega = m_{\omega'}$ for every $\omega' \in W^s_\loc(\omega)$.
	
	\item $(H^u_{\omega,\omega'})_* m_\omega = m_{\omega'}$ for every $\omega' \in W^u_\loc(\omega)$.
	
\end{enumerate}
Then we have the following invariance principle:

\begin{prop}\label{p:bgv}
If $L(A,\mu) = 0$, then there exists an $su$-state for $A$.
\end{prop}

\begin{proof}
We only need to take Proposition~\ref{p.bgvprop1} as our starting point and run the proof of Proposition~\ref{p.bgvprop2}. Note here that the holonomies take values on $\SL(2,\R)$ and are continuous on $\Omega$. Thus the conditions~\eqref{e.integrability} are automatically satisfied.
\end{proof}

One of the main properties of $su$-states is the following.

\begin{prop}\label{p:CntinuousDisintegration}
	If $m$ is an su-state, then it admits a disintegration for which the conditional measures $m_\omega$ depend continuously on $\omega$ and are both $s$-invariant and $u$-invariant.
\end{prop}

\begin{proof}
	We take an $su$-state $m$. Then we run the proof of Lemma~\ref{l.lyap0zp} where we constructed the disintegration $\tilde m$ which is continuous on $\mathrm{supp}(K_\delta)\cap K_\delta$. In the present setting,  we have $K_\delta=\Omega$ since we have canonical holonomies. The result follows.
\end{proof}

By continuity and almost everywhere coincidence, all the invariance properties in the definition of $su$-states may then hold true for every $\omega\in\Omega$. From now on, we always choose such a disintegration for an $su$-state $m$.

\subsubsection{Finiteness of the Set of Energies Admitting $su$-States}

Now we return to the Schr\"odinger case. Recall that by general principles, $L(E) = 0$ implies $E \in \Sigma \subseteq \R$. Since each of these real energies gives rise to an $su$-state for $A^{E}$, let us consider the following set (whose dependence on $\mu$ and $f$ we leave implicit):
$$
\mathcal{F} = \{ E \in \Sigma : \text{ there is an $su$-state for } A^{E} \}.
$$

\begin{lemma}\label{l.infiniteF}
	Suppose that $0 < \alpha \le 1$ and $f \in C^\alpha(\Omega,\R)$ is globally bunched or locally constant. Assume that $\mathcal{F}$ is infinite. Let $p , q \in \Omega$ be two periodic points of $T$. Then $\sigma(H_{p}) = \sigma(H_{q})$.
\end{lemma}

\begin{proof}
	First, we consider the case that $(p)_i=(q)_j$ for some $i, j\in\Z$. Since $\sigma(H_{\omega})=\sigma(H_{T^n\omega})$ for any $\omega$ and for any $n$, we may assume that $(p)_0=(q)_0$. Recall in this case there is a unique $q\wedge p\in W^u_{\mathrm{loc}}(q)\cap W^s_{\mathrm{loc}}(p)$. Assume that $n_p$ is the period of $p$ and $n_q$ is the one of $q$. By our choice of $f$, we may choose $\tilde \Sigma\subseteq \R$ to be a compact interval containing the spectrum $\Sigma$ and $\CU\subseteq \C$ to be a complex neighborhood of $\tilde\Sigma$ where $A^E$ has canonical holonomies for all $E\in\CU$. Recall that under the conditions of the present lemma, the holonomies are holomorphic functions on $\CU$.
	
	By the arguments in the proof of Lemma~\ref{l.lyap0zp}, Corollary~\ref{c.lyap0zp}, and the existence of canonical holonomies, we can find for each periodic $\omega\in\Omega$, a subset $Z_\omega\subseteq\C\bbP^1$ consisting of at most two points that is invariant under $A(\omega)$ and the holonomies. In particular, for the periodic point $p$ with period $n_p$, $Z_{p}$ is invariant under $A_{n_p}(p)$ and
	$$
	 H^s_{p,q\wedge p}(Z_p)=H^u_{q,q\wedge p}(Z_q)\mbox{ whenever } q_0=p_0.
	$$
Note that if $\tr(A_{n_p}(p))\neq 0$, then $Z_p$ must be a subset of the eigendirections of $A_{n_p}(p)$. In particular, for $A_{n_p}(p)$ with nonzero trace, $A_{n_p}(p)$ is elliptic if and only if $Z_p$ is non-real. We let $\{s(E),u(E)\}$ denote the pair of eigendirections of $A^E_{n_p}(p)$. Note that both $s(E)$ and $u(E)$ are continuous on $\tilde\Sigma$ and analytic on each spectral gap or on the interior of each connected component of $\sigma(H_p)$. We define
$$
H^E=H^{u,E}_{q\wedge p, q}\cdot H^{s,E}_{p,q\wedge p},
$$
which are holomorphic in $E$ on a complex neighborhood $\CU$ of $\tilde\Sigma$.

Let $E_0$ be an accumulation point of $\CF$. Then, similarly to the proof of Lemma~\ref{l:zeroLE_bunched_p} or ~\ref{claim1}, we can find a sequence $\{E_n\}_{n\ge 1}$ in $\CF$ so that $E_n\to E_0$, $E_n\neq E_0$, and
$$
Z_p(E_n)=
\begin{cases}
\{s(E_n)\} &\mbox{ for all } n\ge 1,\\
\{u(E_n)\} &\mbox{ for all } n\ge 1, \mbox{ or }\\
\{s(E_n), u(E_n)\} &\mbox{ for all } n\ge 1.
\end{cases}
$$
Thus we may extend the domain of $Z_p(\cdot )$ from $\{E_n,\ n\ge 1\}$ to $\CU$. Then we define
  $$
  Z_q(E):=H^E(Z_p(E))
  $$
and we get that $Z_q(E_n)$ is invariant under $A^{E_n}_{n_q}(q)$ for all $n\ge 1$. By the continuity and analyticity properties of $H^E$, $s(E)$, and $u(E)$, we obtain the following conclusions: if $E_0$ is in a spectral gap, then $Z_q(E)$ is invariant under $A^{E_n}_{n_q}(q)$ for all $E$ in the closure of that spectral gap; if $E_0$ is the interior of a connected component of $\sigma(H_p)$, then $Z_q(E)$ is invariant under $A^{E_n}_{n_q}(q)$ for all $E$ in that connected component.

Now by the same arguments as in the proof of Lemma~\ref{l:zeroLE_bunched_p}, if $E_0$ is in a spectral gap of $H_p$, then it is away from $\sigma(H_q)$ with a uniform distance for all periodic points $q$. But $E_0\in\Sigma$ since it is an accumulation point of $\CF$. Thus $E_0$ can be approximated by $\sigma(H_q)$ for a certain choice of $q$, a contradiction. We may conclude that $E_0 \in \sigma(H_p)$. So we may let $I\subseteq\sigma(H_p)$ be the connected component containing $E_0$. Now we claim that
$$
Z_q(E)=H^E(Z_p(E)) \mbox{ is invariant under }A^E_{n_q}(q) \mbox{ for all }E\in\tilde\Sigma.
$$
	
If $E_0$ is in the interior of $I$, we have already obtained that $\tilde Z_q(E)$ is invariant under $A^E_{n_q}(q)$ for all $E\in I$. If $E_0$ belongs to the boundary of $I$, then similarly to the proof of Lemma~\ref{claim1}, there is an open disk $D$ centered at $E_0$ with ramified (at $E_0$) double cover $\pi :\tilde D\to D$ so that $s$ and $u$ are holomorphic on $\tilde D$.  Thus we may assume $Z_p(E)=\tilde Z_p(\tilde E)$ where $\tilde E\in\pi^{-1}(E)$ and
 $$
 \tilde Z_p(\tilde E)=
 \begin{cases}
 \{s(\tilde E)\}\ &\mbox{ for all }\tilde E\in\tilde D,\\
 \{u(\tilde E)\}\ &\mbox{ for all }\tilde E\in\tilde D, \mbox{ or }\\
 \{s(\tilde E), u(\tilde E)\}\ &\mbox{ for all }\tilde E\in\tilde D.
 \end{cases}
 $$
Then we define
$$
\tilde Z_q(\tilde E):=H^{\pi(\tilde E)}(\tilde Z_p(\tilde E)),
$$
so that $\tilde Z_p(\tilde E)$ is invariant under $A^{\pi(\tilde E)}_{n_q}(q)$ for infinitely many $\tilde E_n\in \tilde D$. By the fact that $H^{\pi(\cdot)}$, $s$, and $u$ are holomorphic on $\tilde D$, we obtain that $\tilde Z_{q}(\tilde E)$ is invariant under $A^{\pi(\tilde E)}_{n_q}(q)$ for all $\tilde E\in \tilde D$. Descending to $D$, we obtain that
$$
Z_q(E):=H^E(Z_p(E))
$$
is invariant under $A^{E}_{n_q}(q)$ for all $E\in D$. In particular, $Z_q(E)$ is invariant under $A^{E}_{n_q}(q)$ for all $E\in (E_0-\rho,E_0+\rho)$, where $\rho>0$ is the radius of $D$.

By the analysis above, we obtain that no matter whether $E_0$ belongs to the boundary or to the interior of $I$, after a finite number of continuations, we get that $Z_q(E)$ is invariant under $A^{E}_{n_q}(q)$ for all $E \in \tilde\Sigma$, as claimed. As in the proof of Lemma~\ref{claim1}, and by the fact that $H^E$ is real for $E$ real, we obtain that $Z_p(E)$ and $Z_q(E)$ are simultaneously real or non-real for all $E \in \tilde \Sigma \supseteq \Sigma$. This clearly implies that
$$
\sigma(H_p) = \sigma(H_q) \mbox{ whenever } p_0=q_0.
$$

Now we remove the condition $(p)_i=(q)_j$ for some $i,j\in\Z$. As in the proof of Lemma~\ref{claim2}, we can find a periodic point $p'$ with some iterates very close to $p$ and some very close to $q$. In particular, $p'_i=p_j$ for some $i, j\in\Z$ and $p'_k=q_m$ for some $k, m\in\Z$. Thus by the first case we consider above, we have
	$$
	\sigma(H_{p})=\sigma(H_{p'})=\sigma(H_{q}).
	$$
This concludes the proof.
\end{proof}

\subsubsection{Proof of Theorem~\ref{t:posLEawayFiniteSet}}

Theorem~\ref{t:posLEawayFiniteSet} is an immediate consequence of  the following theorem.

\begin{theorem}\label{t:posLEawayFiniteSet1}
Suppose $0<\alpha \le 1$ and let $f \in C^\alpha(\Omega,\R)$ be globally bunched or locally constant. If the periodic spectra associated with periodic points of $T$ in $\Omega$ are not all identical, then $\{ E : L(E) = 0 \}$ is finite.
\end{theorem}

\begin{proof}
As $\{ E : L(E) = 0 \} \subseteq \mathcal{F}$, the statement follows from Lemma~\ref{l.infiniteF}.
\end{proof}

\begin{remark}\label{r.nofixedpoint}
Theorem~\ref{t:posLEawayFiniteSet1} is particularly easy to apply when $T$ has a fixed point, as the latter property ensures the presence of a constant potential and all one needs to do in order to show that not all periodic spectra are the same is to use the non-constancy of the sampling function to produce a non-constant periodic potential. However, there are certainly cases of interest where the base dynamics given by $T$ is fixed-point-free. In this case Theorem~\ref{t:posLEawayFiniteSet1} still provides a direct tool for proving that $\{ E : L(E) = 0 \}$ is finite for many globally bunched or locally constant $f \in C^\alpha(\Omega,\R)$, one just needs to take a closer look at the resulting periodic spectra.
\end{remark}

\begin{remark} \label{r:application_to_expanding_maps}
Consider $(\Omega^+,T_+, \mu^+)$ and assume that we can lift $\mu^+$ to an ergodic measure $\mu$ on $(\Omega,T)$ that has a local product structure. Then all our main results of this section, in particular Theorem~\ref{t:posLEawayDiscreteSet} and Theorem~\ref{t:posLEawayFiniteSet}, can be applied to $f \in C^\alpha(\Omega^+,\R)$. Indeed, such an $f$ can be lifted to an $\bar f\in C^\alpha(\Omega, \R)$ that depends only on the future. Then all our results follow since $L(\mu, A^{(E-\bar f)})=L(\mu^+, A^{(E-f)})$.
	\end{remark}

\section{Positivity of the Lyapunov Exponent II}\label{ss.uniformpositivity}

We first show that in the scenario of Subsection~\ref{ss.glundom}, we may remove the finite exceptional set for an open and dense subset of sampling functions. Then we apply similar arguments to the general case discussed in Subsection~\ref{ss.nonundom} and obtain that for a residual set of sampling functions, the discrete exceptional set can be removed. Throughout this section, we again assume that $(\Omega,T)$ is a subshift of finite type with an ergodic measure $\mu$ that has a local product structure. Note that the space $C^\alpha(\Omega,\SL(2,\R))$ is a Banach space with a $C^\alpha$ norm defined by
$$
\|A\|_{0,\alpha}=\|A\|_{0}+\sup_{\omega\neq\omega'}\frac{\|A(\omega)-A(\omega')\|}{d(\omega,\omega')^\alpha},
$$
where $\|A\|_0$ is standard $C^0$ norm $\|A\|_0=\sup_{\omega\in\Omega}\|A(\omega)\|$. Similarly, the space $C^\alpha(\Omega,\R)$ is a Banach space with a $C^\alpha$ norm that can be defined analogously. We say that a subset of $C^\alpha(\Omega,\SL(2,\R))$ has \emph{codimension infinity} if it is locally contained in finite unions of closed submanifolds with arbitrary codimension. The same notion can be defined when we consider a subspace or an open subset of $C^\alpha(\Omega,\SL(2,\R))$.

\subsection{Special Cases: Uniform Positivity in a Dense Open Set}

In this subsection, we assume that $A\in C^\alpha(\Omega,\SL(2,\R))$ is fiber bunched or locally constant, and hence admits canonical holonomies by our earlier discussion.

We first introduce the follow notion of typical cocycles.

\begin{definition}
We say $A$ is \emph{typical} if there are two periodic points $p$ and $q$ with periods $n_p$ and $n_q$ such that $p_0=q_0$ and the following properties hold:
\begin{enumerate}
	\item $A_{n_p}(p)\neq I_2$ and $\tr(A_{n_p}(p))\neq 0$.
	\item Let $\{s(p),u(p)\}\subseteq\C\bbP^1$ be the set of eigendirections of $A_{n_p}(p)$. Then there is no $Z_p\subseteq\{s(p),u(p)\}$ so that $H^u_{q\wedge p, q}\cdot H^s_{p,q\wedge p}\cdot Z_p$ is invariant under $A_{n_q}(q)$.
\end{enumerate}
\end{definition}

Since the definition involves two periodic points $p$ and $q$, we may more precisely say that $A$ is \emph{typical with respect to} $(p,q)$. Note that $A$ might be typical with respect to many other pairs of periodic points as well. Clearly, the defining conditions of a typical cocycle are open in the $C^0$ topology. Thus they are open in the $C^\alpha$ topology as well.

	The notion of a typical cocycle in the present scenario was first introduced in \cite{BGV,BV}. Our version is slightly different from theirs. It is adapted for the proof of Theorem~\ref{t:unifPosLE_1} below. In particular, employing the arguments from \cite{BGV,BV}, one can show the following result. We only sketch the proof for the convenience of the reader.
	
\begin{prop}\label{p:opendense_typical_1}
		The set of typical cocycles as defined above forms a $C^\alpha$-open and dense subset in the set of fiber bunched {\rm (}resp., locally constant{\rm )} cocycles. Moreover, the complement of the set of typical cocycles has codimension infinity.
\end{prop}

\begin{proof}
	 	Following the arguments from \cite{BGV,BV}, for each fixed pair of periodic points $p$ and $q$ with $p_0=q_0$, the complement of the set cocycles satisfying conditions (1) and (2), denoted by $\CB_{p,q}$, is seen to be contained in the union of a finite number of sets of the form
	 	$$
	 	\{A:\CH(A)=0\},
	 	$$
	 	where each $A\mapsto\CH(A)$ is a $C^1$ submersion when restricted to suitable sets of $C^\alpha(\Omega,\SL(2,\R))$. Thus for each fixed pair $(p,q)$, one can show that $\CB_{p,q}$ is a submanifold of $C^\alpha(\Omega,\SL(2,\R))$ with positive codimension. Note that the complement of the set of typical cocycles is
	 	$$
	 	\bigcap_{p,q\in\mathrm{Per}(T):\ p_0=q_0}\CB_{p,q}.
	 	$$
	  Since there are infinitely many such pairs $(p,q)$, the set above is contained in a subset of $C^\alpha(\Omega,\SL(2,\R))$ with codimension infinity. Thus, the complement of the set of typical cocycles has codimension infinity and the set of typical cocycles is open and dense.
\end{proof}

\begin{remark}\label{r:PerturbToTypical}
In particular, we want to mention that one can have the following type of perturbation from \cite{BGV,BV}: for each fixed pair of periodic points $p$ and $q$ with $p_0=q_0$, one can modify the values of $A$ at other points without changing its values at $p$ and $q$ as well as without changing its holonomies on the local stable and unstable sets of these two points.
\end{remark}

We first note the following consequence of our proof of Lemma~\ref{l.infiniteF}, which also recovers one of the results in \cite{BGV}:

\begin{lemma}\label{l:posLEtypical_1}
	Assume that the fiber bunched or locally constant $A\in C^\alpha(\Omega,\SL(2,\R))$ is typical. Then $L(A,\mu)>0$. In particular, there is an open and dense subset $\CG$ of fiber bunched or locally constant cocycles whose complement has codimension infinity and $L(A,\mu)>0$ for all $A\in\CG$.
\end{lemma}

\begin{proof}
	Assume that $L(A,\mu)=0$. Let $p$ and $q$ be two periodic points satisfying the conditions in the definition of typical cocycles. Then by the proof of Lemma~\ref{l.infiniteF}, we know that there is a set $Z_p\subseteq\C\bbP^1$ consisting of at most two points with the following properties:
	\begin{enumerate}
		\item $A_{n_p}(p)\cdot Z_p=Z_p$
		\item 	$H^u_{q\wedge p, q}\cdot H^s_{p,q\wedge p}\cdot Z_p$ is invariant under $A_{n_q}(q)$.
	\end{enumerate}
	Since $p$ and $q$ satisfy the conditions stated in the definition of typical cocycles, we have that $A_{n_p}(p)\neq I_2$ and $\tr(A_{n_p}(p))\neq 0$. Thus property (1) implies that $Z_p$ is a subset of $\{s(p),u(p)\}$. As a consequence, property (2) contradicts condition (2) of the definition of typical cocycles, concluding the proof.
\end{proof}

We note the following consequence of \cite[Theorem 2.8]{backes}.

\begin{prop}\label{p:continuityLE}
	Suppose $(\Omega,T)$ is a subshift of finite type and $\mu$ is $T$-ergodic with a local product structure. Let $f \in C^\alpha(\Omega,\R)$ be globally fiber bunched or locally constant. Then $E \mapsto L(E)$ is continuous on $\R$.
\end{prop}

Indeed, \cite[Theorem 2.8]{backes} implies that the Lyapunov exponent is continuous on the subspace of $C^\alpha(\Omega,\SL(2,\R))$ of globally fiber bunched or locally constant cocycles. If $f\in C^\alpha(\Omega,\R)$ is globally fiber bunched or locally constant, then there is a connected compact interval $\hat \Sigma$ that contains the spectrum $\Sigma = \Sigma_f$ so that $A^{E}$ is fiber bunched or locally constant for all $E \in \hat\Sigma$. Thus $L(E)$ is continuous on $\hat \Sigma$. On the other hand, $L(E)$ is smooth outside of the spectrum as $(T,A^E)$ is uniformly hyperbolic for $E\notin\Sigma$ and the Lyapunov exponent is pluriharmonic on the set of uniformly hyperbolic cocycles. Thus $L(E)$ is continuous on $\R$.

\begin{proof}[Proof of Theorem~\ref{t:unifPosLE_1}]
We focus on the case where $f$ is globally fiber bunched as the proof in the locally constant case is completely analogous.

Fix an $f \in C^\alpha(\Omega,\R)$ that is non-constant and globally fiber bunched. Thus we may find a compact connected interval $\hat\Sigma$ whose interior contains the spectrum $\Sigma_f$ so that $A^{(E-f)}$ is fiber bunched for each $E\in\hat\Sigma$. Note that fiber bunching is a $C^0$ open condition and
	\beq\label{eq:unifclose}
	\sup_{E \in \hat\Sigma}\left\{\|A^{(E-f_1)}-A^{(E-f_2)}\|_{0}\right\} < C \|f_1-f_2\|_{0}.
	\eeq
	Thus, for any open neighborhood $\CU_f \subseteq C^\alpha(\Omega,\R)$ of $f$ that is sufficiently small, we have for each $g \in \CU_f$ that $\Sigma_g \subseteq \hat\Sigma$ and $A^{(E-g)}$ is fiber bunched for all $E \in \hat\Sigma$. In the remaining part of the proof, we fix such a sufficiently small $\CU_f$ and work inside it.
	
	If $\sigma(H_{p,f})=\sigma(H_{q,f})$ for all periodic points $p$ and $q$, then by Remark~\ref{r:PerturbToTypical} we can modify the value of $f$ at $q$ without changing its value along the orbit of $p$. On the other hand, if we choose $E$ on the boundary of the spectrum of $\sigma(H_{q,f})$, we can certainly perturb $f$ to $g$ so that $L(A^{(E-g)},q)>0$.  Thus we may perturb $f$ to a $g$ that is arbitrarily close to $f$ with the property $\sigma(H_{p,g})\neq \sigma(H_{q,g})$. Then we can instead work with $g$.
	
	Thus, we may assume without loss of generality that $f$ is such that $\sigma(H_{p,f})\neq \sigma(H_{q,f})$ for suitably chosen periodic points $p$ and $q$. By the proof of Lemma~\ref{l.infiniteF}, we may further assume that $p_0=q_0$. As described in Subsection~\ref{s:periodic}, we again let $\{s(E),u(E)\}_{E\in\Sigma}$ be the pair of functions associated with the eigendirections of $A^{(E-f)}_{n_p}(p)$. Define $H^E=H^{u,E}_{q\wedge p, q}\cdot H^{s,E}_{p,q\wedge p}$. Then by the proof of Lemma~\ref{l.infiniteF}, if we define $Z_p(E)$ to be
	$$
	Z_p(E)=
	\begin{cases}
	\{s(E)\}\ &\mbox{ for all }E\in\hat\Sigma,\\
	\{u(E)\}\ &\mbox{ for all }E\in\hat\Sigma, \mbox{ or }\\
	\{s(E),u(E)\}\ &\mbox{ for all }E\in\hat\Sigma,
	\end{cases}
	$$
	then the set
	$$
	\left\{E\in\hat\Sigma:\ A^{(E-f)}_{n_q}(q)\cdot H^E\cdot Z_p(E)=H^E\cdot Z_p(E)\right\}
	$$
	is finite. On the other hand, the set
	$$
	\left\{E\in\hat\Sigma:\ A^{(E-f)}_{n_p}(p)=\pm I_2\mbox{ or } \tr(A^{(E-f)}_{n_p}(p))=0\right\}
	$$
	is finite as well. Combining the facts above, we then have that
	$$
	\CB_f:=\left\{E\in\hat\Sigma:\ A^{(E-f)} \mbox{ is not typical}\right\}
	$$
	is finite. Note that for all $E\notin\CB_f$, $A^{(E-f)}$ is typical with respect to $(p,q)$. By Remark~\ref{r:PerturbToTypical} we can modify the values of $A^{(E-f)}$ at different points and keep its values at $p$ and $q$, as well as their holonomies. In particular, after a finite number of perturbations, we can perturb $f$ to $g$ with the following properties. There is a pair of periodic points $(p',q')$ with $p'_0=q'_0$ and $A^{(E-g)}$ is typical with respect to $(p,q)$ for all $E\notin\CB_f$ and typical with respect to $(p',q')$ for all $E\in\CB_f$. Thus we have that $A^{(E-g)}$ is typical for all $E\in\hat\Sigma$. By the fact that the defining properties of typical cocycles are open conditions with respect to the $C^0$ topology, property \eqref{eq:unifclose}, and the compactness of $\hat\Sigma$, we obtain a neighborhood $\CU_g\subseteq\CU_f$ of $g$ so that for each $h\in\CU_g$, we have
$$
L(A^{(E-h)},\mu)>0\mbox{ for all }E\in\hat\Sigma.
$$
By Proposition~\ref{p:continuityLE}, $L(A^{(E-h)})$ is continuous on $\R$. On the other hand, it is well known that $(T,A^{(E-h)})$ is uniformly hyperbolic outside of $\Sigma$ and $L(A^{(E-h)},\mu)$ tends to $\infty$ as $|E|$ tends to $\infty$. Combining all these statements, we find that for each $h\in\CU_g$, we have
$$
\inf_{E\in\R}L(A^{(E-h)},\mu)>0.
$$
This concludes the proof.
\end{proof}

\subsection{General Case: Full Positivity for Generic Sampling Functions}

In this subsection, we return to the general setting of Theorem~\ref{t:posLEawayDiscreteSet}. Note that in this case we have neither the canonical holonomies, nor global existence of holonomies. Moreover, the discrete set can in principle be infinite. To remove the discrete exceptional set, the price we need to pay is that we can only do it for $C^\alpha$-generic sampling functions. For the remaining part of the section, we fix $0 < \alpha \le 1$ and consider the space $C^\alpha(\Omega,\SL(2,\R))$.

We start with a new definition of typical cocycles that is adapted for the purpose of this section.

\begin{definition}
We say $A\in C^\alpha(\Omega,\SL(2,\R))$ is \emph{typical} if there are two periodic points $p$ and $q$ with periods $n_p$ and $n_q$, respectively, such that $p_0=q_0$ and the following properties hold:
\begin{enumerate}
	\item $p$ and $q$ are $\frac\alpha2$-bunched, that is, $2L(A,p)<\frac\alpha2$ and $2L(A,q)<\frac\alpha2$.
	\item $A_{n_p}(p)\neq I_2$ and $\tr(A_{n_p}(p))\neq 0$.
	\item Let $\{s(p),u(p)\}\subseteq\C\bbP^1$ be the set of eigendirections of $A_{n_p}(p)$. Then there is no $Z_p\subseteq\{s(p),u(p)\}$ so that $H^u_{q\wedge p, q}\cdot H^s_{p,q\wedge p}\cdot Z_p$ is invariant under $A_{n_q}(q)$.
\end{enumerate}
\end{definition}

Note that the existence of the holomomies of $p$ and $q$  in condition (3) is guaranteed by condition (1). As in the previous subsection, we may also say that $A$ is \emph{typical with respect to} $(p,q)$, as the definition involves $p$ and $q$.

Define
$$
\mathcal T_\alpha:=\left\{A\in C^\alpha(\Omega,\SL(2,\R)): A \mbox{ is a typical cocycle} \right\}.
$$

It is a standard fact that $A \mapsto L(A,\mu)$ is upper-semicontinuous on $C^\alpha(\Omega,\SL(2,\R))$. In particular, the set
$$
\CL_\alpha = \left\{ A \in C^\alpha(\Omega,\SL(2,\R)) : 2L(A,\mu) < \frac\alpha2 \right\}
$$
is open in $C^\alpha(\Omega,\SL(2,\R))$. Again by \cite[Theorem 3]{kalinin}, if $2L(A,\alpha) < \frac\alpha2$, there exists a periodic point $p$ such that $2L(A,p) < \frac\alpha2$, that is, $p$ is $\frac\alpha2$-bunched. Then, as in the proof of Lemma~\ref{l:zeroLE_bunched_p}, we may use the specification property to produce infinitely many pairs of $\frac\alpha2$-bunched periodic points $(p,q)$ so that $p_0=q_0$. In particular, similarly to Proposition~\ref{p:opendense_typical_1}, we have the following:

\begin{prop}\label{p:opendense_typical_2}
Suppose $(\Omega,T)$ is a subshift of finite type and $\mu$ is a $T$-ergodic measure that has a local product structure. Consider the space $C^\alpha(\Omega,\SL(2,\R))$ for $\alpha>0$ and let $\mathcal T_\alpha$ and $\CL_\alpha$  be defined as above. $\mathcal T_\alpha\cap \CL_\alpha$ forms an open and dense subset of $\CL_\alpha$. Moreover, $\CL_\alpha\setminus\mathcal T_\alpha$ has codimension infinity in $\CL_\alpha$.
\end{prop}

Similarly to Lemma~\ref{l:posLEtypical_1}, Proposition~\ref{p:opendense_typical_2} has the following consequence, which has appeared in \cite{V}. For simplicity, we define
$$
\CP_\alpha=\left\{A\in C^\alpha(\Omega,\SL(2,\R)): L(A,\mu)>0\right\}.
$$

\begin{lemma}\label{l:posLEtypicaol_2}
We have $\mathcal T_\alpha\subseteq \CP_\alpha$. In other words, $L(A,\mu) > 0$ for each $A$ that is typical. Moreover, the set $\CP_\alpha$ contains an open and dense subset of $C^\alpha(\Omega,\SL(2,\R))$ and the complement of $\CP_\alpha$ has codimension infinity.
\end{lemma}

\begin{proof}
	If $A\notin\CL_\alpha$, then $L(A,\mu)\ge \frac\alpha4>0$. If $A\in\CL_\alpha$ is typical, then we may apply the proof of Lemma~\ref{l:posLEtypical_1} to get $L(A,\mu)>0$. However, here we have to use the full strength of Subsection~\ref{ss.nonundom}. Specifically, $\frac\alpha2$-bunching of $p$ and $q$ and the proof of Lemma~\ref{p.bgvprop3} guarantee the existence of the holonomies associated with $p$ and $q$. Then Lemmas~\ref{l:PinSuppK} and \ref{l.lyap0zp} and Corollary~\ref{c.lyap0zp} can be used to guarantee the existence and holonomy-invariance of $Z_p$ and $Z_q$. Once we have all these tools, the proof of $L(A,\mu)>0$ is then identical to the proof of Lemma~\ref{l:posLEtypical_1}.
	
Next, we want to show that the set $\CP_\alpha$ contains an open and dense set. To this end, we fix any $A\in C^\alpha(\Omega,\SL(2,\R))$. If there is an open neighborhood $\CU_A$ of $A$ such that for each $B\in\CU_A$, $L(B,\mu)\ge \frac\alpha4$, then there is nothing we need to say. Otherwise, in any open neighborhood $\CU$ of $A$, we can find a $B\in\CL_\alpha$. Then by Proposition~\ref{p:opendense_typical_2} and the proof above, we can find an open set $\mathcal V\subseteq\CU\cap \mathcal T_\alpha$, which implies that $L(B,\mu)>0$ for each $B\in\mathcal V$.

Finally, it is clear that the complement of $\CP_\alpha$ is contained in $\CL_\alpha\setminus\mathcal T_\alpha$, which has codimension infinity in $\CL_\alpha$. Hence, the complement of $\CP_\alpha$ has codimension infinity in $C^\alpha(\Omega,\SL(2,\R))$ as well.
\end{proof}

Note that this is an improved version of Lemma~\ref{l:posLEtypical_1}, as here we remove the assumption of global bunching or local constancy of $f$.

Now we are ready to generically remove the discrete set that appeared in Theorem~\ref{t:posLEawayDiscreteSet}.

\begin{proof}[Proof of Theorem~\ref{t:fullposLE}]
Via the arguments from the proof of Lemma~\ref{l:posLEtypicaol_2} we can show that the set
	$$
	\CZ_\alpha:=\{f\in C^\alpha(\Omega,\R): L(A^{(f)},\mu)=0\}
	$$
has codimension infinity in $C^\alpha(\Omega,\R)$. In other words, $\CZ_\alpha$ is locally contained in finite unions of closed submanifolds with arbitrary codimension. More precisely, for each $k\in\Z_+$ and each $f\in\CZ_\alpha$, we can find an open neighborhood $\CU_f$ of $f$ and submanifolds $\CM_j$, $1\le j\le m$, each with codimension $k$, so that
	 \beq\label{eq:codimensionk}
	 \left(\CZ_\alpha\cap\CU_f\right)\subseteq \bigcup^m_{j=1}\CM_j.
	 \eeq
On the other hand, if we define the set $\CB_\alpha$ to be
	$$
	\CB_\alpha:=\{f\in C^\alpha(\Omega,\R): L(A^{(E-f)},\mu)\in\CZ_\alpha \mbox{ for some }E\in\R\},
	$$
then for each $g\in\CB_\alpha$, we can find $f\in\CZ_\alpha$ and $E\in\R$ so that $g=E+f$. Thus $\CB_\alpha$ is locally contained in finite unions of submanifolds of arbitrary codimension as well. Indeed, for the $g$ and $f$ above, we may just assume that $f$ is the one in \eqref{eq:codimensionk}. Thus for a fixed $k\in\Z_+$,  for each $\CM_j$ in \eqref{eq:codimensionk}, the set
$$
\CN_j:=\{h\in C^\alpha(\Omega,\R):h-E\in\CM_j\mbox{ for some }E\in\R\}
$$
may be viewed as a submanifold of $C^\alpha(\Omega,\R)$ with codimension $k-1$ whose local charts can be obtained from those of $\CM_j$ and $E\in\R$. In particular, it is nowhere dense if $k\ge 2$. On the other hand, \eqref{eq:codimensionk} clearly implies that the open neighborhood
 $$
 \CU_g=\CU_f+E:=\{h\in C^\alpha(\Omega,\R): h-E\in\CU_f\}
 $$
 of $g$ satisfies
$$
(\CB_\alpha\cap\CU_g)\subseteq \bigcup^m_{j=1}\CN_j.
$$
Since $g\in\CB_\alpha$ and $k\in\Z_+$ can be arbitrarily chosen, we obtain that $\CB_\alpha$ is nowhere dense. Equivalently, we may say that the complement $\CB^c_\alpha$ of $\CB_\alpha$ is residual in $C^\alpha(\Omega,\R)$. By definition of $\CB_\alpha$, we have for each $f\in \CB^c_\alpha$ that
	 $$
	 L(A^{(E-f)},\mu)>0\mbox{ for all }E\in\R,
	 $$
concluding the proof.
\end{proof}

\begin{remark}
	Similarly to Remark~\ref{r:application_to_expanding_maps}, all the main results in this section can be applied to H\"older continuous sampling functions defined on $(\Omega^+,T_+,\mu^+)$, where the lift $\mu$ of $\mu^+$ has a local product structure. Indeed, in this case, $C^\alpha(\Omega^+,\R)$ can be considered as a closed subspace of $C^\alpha(\Omega,\R)$ whose elements depend only on the future. All the perturbations can then be performed within this subspace.
\end{remark}

\section{Applications}\label{s:applications}

All of the results of this paper may be applied to H\"older continuous cocycles defined over any  transitive Anosov diffeomorphism (or transitive, uniformly expanding differentiable map), where $\mu$ is taken to be the equilibrium state of a H\"older continuous potential. By a standard technique one can reduce the cocycles in question to H\"older continuous cocycles over a subshift of finite type via a Markov partition; see, for example, \cite{bowen2, katok}. Although the applicability is much wider, we will focus on a particular case as follows. It is standard result that if an invariant measure $\mu$ of a $C^2$ transitive Anosov diffeomorphism (or a $C^2$ transitive, uniformly expanding map) is absolutely continuous with respect to the volume measure, then it is an equilibrium state of a H\"older continuous potential; see, for example, \cite{bowen2}.


To illustrate this, we choose three differential models that have been widely studied in both the dynamical systems and mathematical physics communities.  The first type of model is given by linear expanding maps of the circle,
$$
T : \R/\Z \to \R/\Z, \quad Tx=kx, \; k\ge 2,
$$
and the measure is taken to be the Lebesgue measure $m$ on $\R/\Z$. One may find some existing results for this case in \cite{CS, BS, BB, bjerklov, DK, SS, vianayang, young, zhang2}. In particular, the case $k=2$ corresponds to the doubling map, which is the most difficult map to study within this family of maps, as it is the least mixing among them. The second type is given by hyperbolic automorphisms of $\R^d/\Z^d$, where $\mu$ is taken to be the Lebesgue measure $m$ on $\R^d/\Z^d$. The most intensively studied case is the famous Arnold cat map, where
$$
T : \R^2/\Z^2 \to \R^2/\Z^2, \quad T = \begin{pmatrix} 2&1\\ 1&1 \end{pmatrix}
$$
and $\mu$ is taken to be the Lebesgue measure $m$ on $\R^2/\Z^2$. One may find earlier results for this case in \cite{CS, BS, SS, young}. It is clear that both linear expanding maps of the circle and hyperbolic toral automorphisms meet all the conditions necessary to apply our main theorems in Sections~\ref{s:PosLya1} and \ref{ss.uniformpositivity}. In particular, they all have a fixed point.

Our theorems then yield the following results. To unify the statements, we let $(\Omega,T, \mu)$ be any of the following: $(\R/\Z, T_k, m)$, where $T_kx=kx$ and $k\ge 2$ is an integer;  $(\R^d/\Z^d, T_A, m)$ where $d\ge 2$ and $T_A$ is the hyperbolic toral automorphism generated by some hyperbolic $A\in\SL(d,\Z)$. Recall that for a sampling function $f$, we set $L(E)=L(\mu,A^{(E-f)})$ and define
$$
\CZ_f:=\{E: L(E)=0\}\subseteq \R.
$$
For $0 < \alpha \le 1$ and $\lambda > 0$, we $C^\alpha_\lambda(\Omega,\R)=\{f\in C^\alpha(\Omega,\R): \|f\|_\infty<\lambda\}$.

\begin{theorem}\label{t:posLEhyp}
	Let $(\Omega,T,\mu)$ be as above and let $0<\alpha\le 1$. For all non-constant $f\in C^\alpha(\Omega, \R)$, $\CZ_f$ is a discrete set. Moreover, $\CZ_f=\varnothing$ for $f$'s in a residual subset of $C^\alpha(\Omega,\R)$. There is $\lambda_0=\lambda_0(\alpha)>0$ such that $\CZ_f$ is a finite set for all non-constant $f\in C^\alpha_{\lambda_0}(\Omega,\R)$. Finally, there is an open and dense subset $\CO^\alpha$ of $C^\alpha_{\lambda_0}(\Omega,\R)$ such that for all $f\in \CO^\alpha$, $\inf_{E\in\R}L(E) > 0$.
\end{theorem}

If we introduce a coupling constant $\lambda$ , then we have the following immediate consequence of Theorem~\ref{t:posLEhyp}.

\begin{coro}\label{c:no_transition}
	Let $(\Omega,T,\mu)$ and $\alpha$ be as in Theorem~\ref{t:posLEhyp}. Fix a non-constant $f\in C^\alpha(\Omega, \R)$. Then  $\CZ_{\lambda f}$ is a discrete set for all $\lambda>0$. Moreover, there is a $\lambda_0=\lambda_0(\|f\|_\infty,\alpha)>0$ such that $\CZ_{\lambda f}$ is finite for all $0<\lambda<\lambda_0$.
\end{coro}

\begin{remark}\label{r:posLEhyp}
	To the best of our knowledge, if we take $T$ to be the doubling map for $d=1$ or the Arnold cat map for $d\ge 2$, then the results we stated in Theorem~\ref{t:posLEhyp} and Corollary~\ref{c:no_transition} are the first global results that do away with smallness or largeness assumptions for the coupling constant. 
In the large coupling regime, Herman's subharmonicity trick \cite{herman} can be applied (for trigonometric polynomials), and in the (perturbatively!) small coupling regime, the perturbative analysis of Chulaevsky-Spencer \cite{CS} and Sadel-Schulz-Baldes \cite{SS1, SS} can be applied. Other methods get around changing the coupling constant by changing the base dynamics instead, specifically to increase its hyperbolicity; compare Bourgain-Bourgain-Chang \cite{BB} and Bjerkl\"ov \cite{bjerklov}.
\end{remark}

\begin{remark}\label{r:sample}
	Taking the doubling map as an example, we give two sample computations. First, we show how to reduce a H\"older continuous cocycle on $\R/\Z\times\R^2$ to one on $\Omega\times\R^2$ for some subshift of finite type. Let $\Omega^+=\{0,1\}^\N$ and $(\Omega^+,T_+,\mu^+)$ be the one-sided Bernoulli shift. Here we choose $\mu^+=\tilde \mu^\N$ where $\tilde\mu(0)=\tilde\mu(1)=\frac12$. Then it is well know that the map
	$$
	\pi:\Omega^+\to\R/\Z, \; \omega^+ \mapsto \sum^\infty_{n=0}\frac{\omega^+_n}{2^{n+1}}
	$$
	codes the dynamics of doubling map $(\R/\Z, T_2, m)$ to that of $(\Omega^+, T_+, \mu^+)$ since $T_2\circ\pi=\pi\circ T_+$ and $\pi_*\mu^+=m$. In particular, for any cocycle map $A:\R/\Z\to\SL(2,\R)$, we set $A^+:\Omega^+\to\SL(2,\R)$ where $A^+=A\circ \pi$. Then the dynamics of the cocycle $(T,A)$ is the same as the one of  $(T_+, A^+)$. In particular, $L(A,m)=L(A^+,\mu^+)$. Now we consider the full shift space $(\Omega,T,\mu)$ whose one-sided shift is $(\Omega^+,T_+,\mu^+)$, as described above. By setting $\bar A(\omega)=A^+(\pi^+\omega)$, we obtain a cocycle $(T,\bar A)$ that shares the dynamics with $(T_+,A^+)$. It is clear that $\bar A$ is $\alpha$-H\"older continuous as long as $A^+$ is, since $d(\pi^+\omega,\pi^+\tilde \omega)\le d(\omega,\tilde\omega)$. So we just need to show that the H\"older continuity can be carried over from $A$ to $A^+$. This in turn follows from the following straightforward estimate:
	$$
	|\pi(\omega^+)-\pi(\tilde \omega^+)|\le d(\omega^+,\tilde\omega^+)^{\ln 2}.
	$$
	In particular, $\alpha$-H\"older continuity of $A$ implies $(\alpha\ln 2)$-H\"older continuity of $A^+$ since
	\begin{align*}
		\|A^+(\omega^+)-A^+(\tilde\omega^+)\|&=\|A(\pi \omega^+)-A(\pi \tilde\omega^+)\|\\
		&\le C|\pi\omega^+-\pi\tilde \omega^+|^\alpha\\
		&\le Cd(\omega^+,\tilde \omega^+)^{\alpha\ln2}.
		\end{align*}
	
Next, we compute some explicit choices for the value of $\lambda_0$ appearing in Theorem~\ref{t:posLEhyp} and Corollary~\ref{c:no_transition}, when the base dynamics in question is the doubling map. Clearly, the above process still works if we replace $A:\R/\Z\to\SL(2,\R)$ by $f:\R/\Z\to\R$. Given $f\in C^\alpha(\R/\Z,\R)$, we may instead consider the corresponding $\bar f\in C^{\alpha\ln2}(\Omega,\R)$. In particular, $\|f\|_\infty=\|\bar f\|_\infty$. We want to find a $\lambda_0$ so that $f$ is globally bunched if $\|f\|_\infty<\lambda_0$. In other words,
$$
A^E(\omega)=\begin{pmatrix}E-\bar f(\omega) &-1\\ 1 &0\end{pmatrix}
$$
is fiber bunched for all $E\in [-2-\|f\|_\infty,2+\|f\|_\infty]$. To simplify the computation, we simply impose fiber bunching at step $1$. That it, we want for all  $E\in[-2-\|\bar f\|_\infty, 2+\|\bar f\|_\infty]$ that
$$
\|A^E(\cdot)\|_\infty<e^{\frac{\ln2}{2}\alpha}=2^{\frac\alpha2}.
$$
Recall the fiber bunching condition is only assumed to ensure the existence of stable and unstable holonomies. Thus, by the construction of the holonomies from the proof of Lemma~\ref{p.bgvprop3}, it is clear that we may reduce the condition above to the following condition. For each $E\in [-2-\|\bar f\|_\infty, 2+\|\bar f\|_\infty]$, there is a $P(E)\in\mathrm{SL}(2,\R)$ so that
\beq\label{eq:globalbunching}
\|P(E)^{-1}A^E(\cdot)P(E)\|_\infty<2^{\frac\alpha2}.
\eeq
First, we take care of the $E$'s that are away from $\pm 2$. For each $E\in (-2,2)$, a direct computation shows that
$$
P(E)^{-1}\begin{pmatrix}E&-1\\ 1 &0\end{pmatrix} P(E)\in \mathrm{SO}(2,\R)
$$
which has norm one and where
$$
P(E)=\begin{pmatrix}\frac{\sqrt2}{(4-E^2)^{\frac14}}&0\\ \frac{E}{\sqrt2(4-E^2)^{\frac14}}&\frac{(4-E^2)^{\frac14}}{\sqrt2}\end{pmatrix}.
$$
If we choose $\lambda_0$ so that for all $E\in [-2+\lambda_0,2-\lambda_0]$ and all $|\lambda|<\lambda_0$, we have
$$
\left\|P(E)^{-1}\begin{pmatrix}\lambda&0\\ 0 &0\end{pmatrix} P(E)\right\|< 2^{\frac\alpha2}-1,
$$
then we have \eqref{eq:globalbunching} for any $\|\bar f\|_\infty=\|f\|_\infty<\lambda_0$ and all $E\in[-2+\lambda_0, 2-\lambda_0]$. It is straightforward to see that
$$
P(E)^{-1}\begin{pmatrix}\lambda&0\\ 0 &0\end{pmatrix} P(E)=\begin{pmatrix}\lambda&0\\ -\frac{E\lambda}{\sqrt{4-E^2} }&0\end{pmatrix}.
$$
Thus we have fiber bunching for all $E\in[-2+\lambda_0, 2-\lambda_0]$ if for all such $E$'s and for all $|\lambda|<\lambda_0$, we have
$$
|\lambda|+\left|\frac{E\lambda}{\sqrt{4-E^2}}\right|<2^\frac\alpha2-1.
$$
Clearly, it suffices to have
$$
\lambda_0+\frac{\lambda_0}{\sqrt{\lambda_0-\lambda_0^2}}<2^{\frac\alpha2}-1,
$$
which in turn can be guaranteed, for example, by the condition $3\sqrt{\lambda_0}\le 2^\frac\alpha2-1$. In particular, if we choose any
\beq\label{eq:globalbunching2}
0<\lambda_0\le \frac{(2^\frac\alpha2-1)^2}9,
\eeq
then we have fiber bunching for all $E\in[-2+\lambda_0, 2-\lambda_0]$ and for all $\|f\|_\infty<\lambda_0$.

Now we take care of the energies $E\in [-2-\lambda_0, -2+\lambda_0]\cup [2-\lambda_0, 2+\lambda_0]$. Take $E=2$ for example. Then we have
$$
P_a^{-1}\begin{pmatrix}2 &-1\\ 1& 0\end{pmatrix}P_a=\begin{pmatrix}1 &a\\ 0& 1\end{pmatrix},
$$
where $a>0$ and
$$
P_a=\begin{pmatrix}\frac1{\sqrt{a}} &-\sqrt{a}\\ \frac1{\sqrt{a}}& 0\end{pmatrix}.
$$
It is easy to see that for $0<a<1$, we have
$$
\left\|\begin{pmatrix}1 &a\\ 0& 1\end{pmatrix}\right\|<1+2a.
$$
On the other hand, we can see via a straightforward computation that
$$
P_a^{-1}\begin{pmatrix}\lambda&0\\ 0 &0\end{pmatrix}P_a=\begin{pmatrix}0&0\\ -\frac\lambda{a}&\lambda \end{pmatrix}.
$$
Thus it suffices to choose $a<1$ and $\lambda_0 > 0$ so that for all $|\lambda| \le 2\lambda_0$, we have
$$
1+2a+\frac{|\lambda|}{a}+|\lambda|<2^\frac\alpha2,
$$
which may be guaranteed by
$$
a+\frac{\lambda_0}a+\lambda_0<\frac{2^\frac\alpha2-1}{2}.
$$
Clearly, we may choose $a=\frac14(2^\frac\alpha2-1)<\frac14$. It is then easy to see that if we choose any $\lambda_0$ such that
\beq\label{eq:globalbunching3}
0<\lambda_0\le \frac{(2^\frac\alpha2-1)^2}{20},
\eeq
then we have fiber bunching for all $E\in [2-\lambda_0, 2+\lambda_0]$ and for all $f$ with $\|f\|_\infty<\lambda_0$. A similar computation shows that the $\lambda_0$ in \eqref{eq:globalbunching3} works for $E\in [-2-\lambda_0, -2+\lambda_0]$ as well. Combining \eqref{eq:globalbunching2} and $\eqref{eq:globalbunching3}$, we see that in the statement of Theorem~\ref{t:posLEhyp} and Corollary~\ref{c:no_transition} for the doubling map, we may choose
$$
\lambda_0=\frac{(2^\frac\alpha2-1)^2}{20}.
$$
\end{remark}

\begin{remark}\label{r:globalbunching}
	The computation of $\lambda_0$ in Remark~\ref{r:sample} actually works for $A^E$ defined on any subshift of  finite type $(\Omega,T,\mu)$. Moreover, since we do not have the coding process as in Remark~\ref{r:sample}, we have that $f\in C^\alpha(\Omega,\R)$ is globally bunched if
	\beq\label{eq:globalbunching4}
	\|f\|_\infty\le \lambda_0=\frac{(e^\frac\alpha2-1)^2}{20}.
	\eeq
	In particular, this value of $\lambda_0$ works for Theorem~\ref{t:posLEmarkov} below.
	\end{remark}

Let us now apply our results to Markov chains. We consider the full shift $(\CA^\Z,T)$, where $\CA=\{1,\ldots, \ell\}$. Let $P=(P_{ij})_{1\le i,j\le \ell}$ be a stochastic matrix, in other words, $P_{ij}\ge 0$ and $\sum^{\ell}_{j=1}P_{ij}=1$. Then there is a unique probability vector $\vec p=(p_1,\ldots p_\ell)$ (i.e., $p_i > 0$ and $\sum^{\ell}_{i=1} p_i=1$) such that $\sum^{\ell}_{i=1} p_i P_{ij} = p_j$. Assume that $P$ is irreducible, that is, for all $i,j \in \CA$, there is $n \in \Z_+$ such that the $(i,j)$-entry of $P^n$ is positive. Now we define the measure $\mu$ on $\CA^\Z$ via
$$
\mu([0;k_0, \ldots, k_n])=p_{k_0}\prod^{n-1}_{i=0}P_{k_ik_{i+1}}.
$$
Such a measure $\mu$ is called a Markov measure. By a standard result, the topological support of $\mu$ is a subshift of finite type $\Omega$ with the adjacency matrix $A=(a_{ij})$ such that $a_{ij}=1$ whenever $p_{ij}>0$ and $a_{ij}=0$ otherwise. Thus we may instead consider the space $(\Omega, T,\mu)$. Moreover, $\mu$ is $T$-ergodic if and only if $P$ is irreducible. Consider its associated one-sided space $(\Omega^+,T_+,\mu^+)$. It is a standard result that $\mu^+$ is the unique equilibrium state of the potential $\phi(\omega^+)=-\log P_{\omega^+_0\omega^+_1}$, which is locally constant; see, for example, \cite{young2}. Thus by Lemma~\ref{l:es_to_bdd}, $\mu$ has the bounded distortion property, and hence a local product structure as well.

\begin{theorem}\label{t:posLEmarkov}
Let $(\Omega,T,\mu)$ be a Markov chain as described above. Fix $0 < \alpha \le 1$. Then there is a residual set $\CG^\alpha\subseteq C^\alpha(\Omega,\R)$ such that $\CZ_f=\varnothing$ for all $f\in\CG^\alpha$. There are $\lambda_0=\lambda_0(\alpha)>0$ and an open dense subset $\CO^\alpha\subseteq C^\alpha_{\lambda_0}(\Omega,\R)$ with the following property. For each $f\in \CO^\alpha$, we have $\inf _{E\in\R}L(E)>0$. If in addition $(\Omega,T)$ has a fixed point {\rm (}which happens if and only if  $P_{ii}>0$ for some $1\le i\le \ell${\rm )}, $\CZ_f$ is a discrete set for all non-constant $f\in C^\alpha(\Omega,\R)$ and it is a finite set for all non-constant $f\in C^\alpha_{\lambda_0}(\Omega,\R)$ or for all non-constant $f$ that are locally constant. In particular, $\CZ_{\lambda f}$ is discrete for all $\lambda>0$ and finite for all $0<\lambda<\lambda_0$ for all non-constant $f\in C^\alpha(\Omega,\R)$. If $f$ is locally constant and non-constant, then $\CZ_{\lambda f}$ is a finite set for all $\lambda>0$.
\end{theorem}

\begin{remark}
Reiterating what we said in Remark~\ref{r.nofixedpoint}, even if $(\Omega,T)$ does not have a fixed point (i.e., when $P_{ii} = 0$ for every $1\le i\le \ell$), we can work with periodic spectra of higher periods and test for non-coincidence of two of them. In concrete cases this procedure is easy to implement and will in many cases lead to the desired result. For instance, we can apply it to the last example we present in the end of this section.
\end{remark}

Note that the Anderson model is a special case of the Markov chains described above, provided that the single-site measure is supported on a finite set. Indeed, such models may be generated as follows. Let $\mu$ be a probability measure on the full shift space $\CA^\Z$ that is generated by a single site measure $\bar \mu\{i\}=p_i$ where $\vec p=(p_1,\cdots p_\ell)$ is a probability vector. It is clearly a Markov chain with the same probability vector and with the stochastic matrix $p_{ij}=p_j$. Thus, we have the following corollary of Theorem~\ref{t:posLEmarkov}.

\begin{coro}\label{c:PosLE_Anderson}
Consider the full shift space $(\CA^\Z,T,\mu)$, where $\mu=\tilde\mu^\Z$ and $\tilde\mu$ is a probability measure on $\CA = \{1,\ldots \ell\}$ that has full support. Then all the conclusions that we stated in Theorem~\ref{t:posLEmarkov} hold true. In particular, if $f$ is locally constant and non-constant, then $\CZ_{\lambda f}$ is a finite set for all $\lambda>0$.
\end{coro}

In particular, the Anderson model is generated by a sampling function $f : \CA^\Z \to \R$ that depends only the $0$-th position. Note that such a function is in particular locally constant. Corollary~\ref{c:PosLE_Anderson} implies the finiteness of $\CZ_{\lambda f}$ for all such $f$'s that are non-constant. Of course, in this case, the celebrated F\"urstenberg's Theorem yields uniform positivity of the Lyapuonv exponent. However, the finiteness of $\CZ_f$ for all non-constant locally constant $f:\CA^\Z\to\R$ already may not be directly obtained from F\"urstenberg's Theorem. Moreover, our result is basically sharp. Indeed, there are plenty of examples where $\CZ_f$ is not empty for locally constant and non-constant $f:\CA^\Z\to\R$, see \cite{bucaj}. Nevertheless, the finiteness of $\CZ_f$ can already be a starting point to prove full spectral localization.

For the reader's convenience, we provide an example with the property $\CZ_f\neq\varnothing$, where $f$ is a non-constant locally constant function defined over a Markov chain. To give such an example, let us show that the well-known random dimer model (cf., e.g., \cite{BG, DWP}) is covered by our framework. The random dimer model arises from the standard Bernoulli-Anderson model by doubling up the sites. That is, with $\{ \omega_n \}_{n \in \Z}$ i.i.d.\ random variables taking two different values, say $0$ and $\lambda$ with probability $0<p<1$ and $1-p$, the potentials are given by $V_\omega(2n) = V_\omega(2n+1) = \omega_n$. To realize these potentials in our framework, consider the subshift of finite type $\Omega$ over the alphabet $\{ 1, 2, 3, 4 \}$ with the adjacency matrix
	$$
	A = \begin{pmatrix} 0 & 0 & 1 & 0 \\ 0 & 0 & 0 & 1 \\ 1 & 1& 0 & 0 \\ 1 & 1 & 0 & 0 \end{pmatrix}.
	$$
	 The measure $\mu$ is the Markov measure generated by the following probability vector and the stochastic matrix
	 $$
	 \vec p=\left(\frac p2, \frac{1-p}2,\frac p2, \frac{1-p}2\right),\
	 P = \begin{pmatrix} 0 & 0 & 1 & 0 \\ 0 & 0 & 0 & 1 \\ p & 1-p& 0 & 0 \\ p & 1-p & 0 & 0 \end{pmatrix}.
	 $$
	 The sampling function $f:\Omega\to\R$ is  generated by $\bar f : \{ 1,2,3,4\} \mapsto \{ 0 , \lambda \}$, $\bar f(1) = \bar f(3) = 0$, $\bar f(2) = \bar f(4) = \lambda$ via $f(\omega)=\bar f(\omega_0)$ which is locally constant. It is readily checked that the resulting model is indeed the random dimer model. It is well known, and in fact easy to see, that for $-2 < \lambda < 2$, $A^{(E-f)}_n(\omega)$ is bounded for all $n$ at energies $0$ and $\lambda$. Thus $\{0,\lambda\}\subseteq \CZ_f$. Although this system has no fixed point, we do have that $f$ is constant on the orbit of $\omega\in\Omega$ where $\omega_{2n}=1, \omega_{2n+1}=3$. Note that in statement of Theorems~\ref{t:posLEawayDiscreteSet} and \ref{t:posLEawayFiniteSet1}, the fixed point is only there to produce a constant potential $V_\omega(n)$. Thus, Theorem~\ref{t:posLEmarkov} can still be applied to obtain the finiteness of $\CZ_f$. However, for this model, we can provide more information. It actually follows from F\"urstenberg's Theorem that the Lyapunov exponent is positive away from these two energies $\{0,\lambda\}$. This shows that in this particular case $\CZ_f=\{ 0 , \lambda\}$.

\end{document}